\documentclass[a4paper]{amsart}
\usepackage{amsmath}
\usepackage{amssymb}%This gives us symbols like \square
\usepackage{verbatim}%This gives us the \begin{comment} environment.
\usepackage{amsthm}%This gives us Definition style of theorem.
\usepackage{mathrsfs}%This gives us \mathscr#

\usepackage{graphicx}
\usepackage{mathtools}
\usepackage[colorlinks,pdfdisplaydoctitle,linkcolor=blue,citecolor=blue]{hyperref}
\usepackage{caption}
\usepackage{subcaption}
\usepackage{url}
\usepackage{epstopdf}
\usepackage{pgf,tikz}
\usepackage{bbm}

\usepackage{extarrows}
\usetikzlibrary{arrows}
\usepackage[ruled,lined]{algorithm2e}

\usepackage{ dsfont }
\usepackage{listings}
\usepackage{bbm}
\usepackage{color}
\usepackage{ dsfont }
\usepackage{enumerate}
\usepackage{todonotes}
\DeclareMathOperator*{\argmin}{arg\,min}
\renewcommand{\Delta}{\triangle}

\definecolor{darkblue}{rgb}{0,0,0.7}

\definecolor{darkgreen}{rgb}{0.01,0.75,0.24}

\def \Ee[#1]{\mathcal{E}^{\text{{#1}}}}

\def\pa[#1,#2]{\frac{\partial {#1}}{\partial {#2}} }

\def\idom[#1,#2,#3]{\int_{#1}\hspace{1pt} {#2} \hspace{1pt} \text{d}{#3}}
\def\res[#1,#2]{\left.{#1}\right|_{#2}}

\def\var[#1,#2]{\langle \delta \mathcal{E}^{\text{{#1}}}({#2}),v\rangle}
\def\vars[#1,#2,#3]{\langle \delta^2\mathcal{E}^{\text{{#1}}}({#2})v,{#3}\rangle}
\def\vard[#1,#2,#3,#4]{\langle \delta\mathcal{E}^{\text{{#1}}}({#2})-\delta\mathcal{E}^{\text{{#3}}}({#4}),v\rangle}

\def\E{\mathbb{E}}

\newcommand{\Cov}[1]{\mathrm{Cov}\left[#1\right]}

\newcommand{\be}{\begin{equation}}
\newcommand{\en}{\end{equation}}
\newcommand{\ben}{\begin{equation*}}
\newcommand{\enn}{\end{equation*}}
\newcommand{\bea}{\begin{aligned}}
\newcommand{\ena}{\end{aligned}}

\def\ba#1\ena{\begin{align}#1\end{align}}
\def\ban#1\enan{\begin{align*}#1\end{align*}}

\theoremstyle{plain}
\newtheorem{thm}{Theorem}[section]

\newtheorem{Lem}[thm]{Lemma}

\newtheorem{ass}{Assumption}

\newtheorem{rem}[thm]{Remark}

\numberwithin{equation}{section}

\newcommand{\X}{\mathcal{X}}
\newcommand{\Y}{\mathcal Y}
\newcommand{\V}{\mathcal V}
\newcommand{\op}{F}
\newcommand{\norm}[2]{\left\Vert #1\right\Vert_{#2}}
\newcommand{\ip}[2]{\left\langle #1, #2\right\rangle}
\newcommand{\taylor}[2]{\op\left(#1\right) + \op'\left[#1\right]\left(#2-#1\right)}
\newcommand{\udag}{u^\dagger}
\newcommand{\ydag}{y^\dagger}
\newcommand{\cS}[1]{\mathcal S \left(#1; W\right)}

\newcommand{\err}{\mathrm{err}}
\newcommand{\Ctc}{C_{\mathrm{tc}}}
\newcommand{\Cdec}{C_{\mathrm{dec}}}

\begin{document}

%\title[Dynamic iRGNM in the presence of white noise]{Dynamic iteratively regularized Gauss--Newton method in the presence of white noise}
\title[A Dynamic variant of iRGNM with sequential data]{On a dynamic variant of the iteratively regularized Gauss-Newton method with sequential data}
\author[N. K. Chada] {Neil K. Chada}
\address{Applied Mathematics and Computational Science Program, King Abdullah University of Science and Technology, Thuwal, 23955, KSA}
\email{neilchada123@gmail.com}

\author[M. A. Iglesias] {Marco A. Iglesias}
\address{School of Mathematical Sciences, University of Nottingham, Nottingham, NG72RD, UK}
\email{marco.iglesias@nottingham.ac.uk}

\author[S. Lu] {Shuai Lu}
\address{School of Mathematical Sciences, Fudan University, 200433 Shanghai, China}
\email{slu@fudan.edu.cn}

\author[F. Werner] {Frank Werner}
\address{Institut f\"{u}r Mathematik, University of Wuerzburg, Emil--Fischer--Str. 30, 97074 W\"{u}rzburg}
\email{frank.werner@mathematik.uni-wuerzburg.de}

\subjclass{94A12, 86A22, 60G35, 62M99.}
\keywords{Inverse problems, regularization theory, Gauss--Newton method, convergence rates, dynamical process}

\begin{abstract}
For numerous parameter and state estimation problems, assimilating new data as they become available can help produce accurate and fast inference of unknown quantities. While most existing algorithms for solving those kind of ill-posed inverse problems can only be used with a single instance of the observed data, in this work we propose a new framework that enables existing algorithms to invert multiple instances of data in a sequential fashion. Specifically we will work with the well-known iteratively regularized Gauss--Newton method (IRGNM), a variational methodology for solving nonlinear inverse problems. We develop a theory of convergence analysis for a proposed dynamic IRGNM algorithm in the presence of Gaussian white noise. We combine this algorithm with the classical IRGNM to deliver a practical (hybrid) algorithm that can invert data sequentially while producing fast estimates. Our work includes the proof of well-definedness of the proposed iterative scheme, as well as various error bounds that rely on standard assumptions for nonlinear inverse problems. We use several numerical experiments to verify our theoretical findings, and to highlight the benefits of incorporating sequential data. The context of the numerical experiments comprises various parameter identification problems including a Darcy flow PDE, and that of electrical impedance tomography.
\end{abstract}

\maketitle

\section{Introduction}\label{se1}

A common problem in numerous scientific disciplines is the estimation of some unknown function $u^{\dagger} \in \mathcal{X}$, from observations
$y^{\dagger} \in \mathcal{Y}$ of the form
\begin{equation}
\label{eq:ip}
y^{\dagger} = F(\udag),
\end{equation}
where we assume that $F : D(F) \subset \X \to \Y$ is a nonlinear continuous operator acting between two Hilbert spaces $\X$ and $\Y$ with domain of definition $D(F)$. Due to the unavoidable presence of observational noise in real applications, the idealised equation (\ref{eq:ip}) must be replaced by
\begin{equation}
\label{eq:ip_noise}
y^{\delta} = F(\udag) + \sigma \xi,
\end{equation}
where $\xi$ could be a uniformly bounded noise or some other random noise.
Problems associated with \eqref{eq:ip} or \eqref{eq:ip_noise} are commonly referred to as inverse problems \cite{AMS10,AT87}, concerned with the estimation of some unobservable parameter or quantity of interest. Such examples of particular applications include, but not limited to, geophysical sciences, medical imaging and numerical weather prediction \cite{LB02,MW06,ORL08}.

%\todo[inline]{Add references here: \cite{BHM09,BNO97,QJ00,QJ11,JZ13}.}

Given observed data $y^{\delta} \in \Y$, a well-known regularization method to recover $\udag$ from $y^{\delta}$ is the iteratively regularized Gauss--Newton method (IRGNM) \cite{KNS08,SK11}, proposed originally by Bakushinskii \cite{ABB92}. At each iteration, the IRGNM solves a variational problem of the form
\begin{equation}
\label{eq:var_gn}
\hat u_{n+1} := \argmin_{u \in \mathcal{X}}\left[\norm{F\left(\hat u_n\right) + F'\left[\hat u_n\right] \left(u-\hat u_n\right)  - y^{\delta}}{\Y}^2  + \alpha_n \norm{u-\hat u_0}{\X}^2\right],
\end{equation}
where $\hat u_0 \in \X$ is some initial guess, $F'\left[u\right]$ is the  Fr\'{e}chet (or some other) derivative of $F$ at $u$, and $\{\alpha_n\}_{i=1}^{\infty}$ is a sequence of regularization parameters  chosen such that
$$
\alpha_0 \leq 1, \quad \quad \alpha_n \searrow 0, \quad \quad 1 \leq \frac{\alpha_n}{\alpha_{n+1}} \leq C_{\mathrm{dec}}, \quad \mathrm{for \ all } \ n \in \mathbb{N},
$$
for some constant $C_{\mathrm{dec}}$. Typically one uses $\alpha_n=\alpha_0C^{-n}_{\mathrm{dec}}$. Alternatively, we can express the minimization procedure of \eqref{eq:var_gn} in terms of the first order optimality condition as
\begin{equation}
\label{eq:gn_up}
\hat u_{n+1} = \hat u_n - (F'\left[\hat u_n\right]^*F'\left[\hat u_n\right] + \alpha_n \text{id}{\X})^{-1}\left(F'\left[\hat u_n\right]^*\left(F\left(\hat u_n\right)-y^{\delta}\right) + \alpha_n\left(\hat u_n-\hat u_0\right)\right),
\end{equation}
with the adjoint $F'\left[u\right]^* : \Y \to \X$ of $F'\left[u\right] : \X \to \Y$.
%The IRGNM has been well-studied and applied in the context of different inverse problems.
Convergence (rate) analysis for the classical IRGNM can be found in \cite{BNO97,QJ00,Kaltenbacher97} and extension towards the random noise or the Banach space setting can be found in \cite{BHM09,QJ11,JZ13} and references therein.

In this study, we assume that we have \textit{sequential noisy observations} of the form
\begin{equation}\label{eq:model_dyn}
Y_n = \op\left(\udag\right) + \sigma \xi_n,\quad  n=1,2,\ldots,
\end{equation}
where $\xi_n : \Y \to L^2 \left(\Omega, \mathcal A,\mathbb P\right)$ are independent Hilbert space processes (bounded linear functionals) with an underlying probability space $\left(\Omega, \mathcal A, \mathbb P\right)$ and $\E{\xi_n} = 0$, $\norm{\Cov{\xi_n}}{\Y\to \Y} \leq 1$.
Note that the model \eqref{eq:model_dyn} means, that for every $g \in \Y$, the quantity $\ip{\xi}{g} := \xi\left(g\right)$ is a real-valued random variable. However, in general it does not hold $\xi_n \in \Y$, and hence the observations $Y_n$ in \eqref{eq:model_dyn} do not belong to $\Y$. This implies that the model \eqref{eq:model_dyn} has to be understood in a weak sense, meaning that for every $g \in \Y$, the random variable $\ip{Y_n}{g}$ can be observed.

Observations of the form \eqref{eq:model_dyn} are available in nearly all practical applications, but usually not treated as such. Instead, sequential observations of the form \eqref{eq:model_dyn} are used to generate a final datum
\begin{align}\label{eq:avg_noise}
Z_N = N^{-1} \sum_{i=1}^NY_n = \op\left(\udag\right)+ \frac{\sigma}{N} \sum_{n=1}^N \xi_n,
\end{align}
as the average of the (first) $N$ sequential observations. The rationale behind is that the covariance operator of $Z_N$ satisfies
\[
\Cov{Z_N} = \Cov{\frac{\sigma}{N} \sum_{n=1}^N \xi_n} = \frac{\sigma^2}{N^2} \sum_{n=1}^N \Cov{\xi_n},
\]
and hence the noise level of $Z_N$ is $\frac{\sigma}{\sqrt{N}}$ instead of $\sigma$ for each of the observations $Y_n$ in \eqref{eq:model_dyn}. In our situation, where we assume that all the available data (i.e. $Y_n$ in \eqref{eq:model_dyn} or $Z_N$ in \eqref{eq:avg_noise} are a.s. not elements in $\Y$), the above classic IRGNM (cIRGNM) is not directly available. However, noticing that $\norm{\cdot - y}{\Y}^2$ in \eqref{eq:var_gn} is in finite dimensions just the negative log-likelihood functional of the normal distribution, it seems reasonable to replace $\norm{\cdot - y}{\Y}^2$ by
\begin{equation}\label{eq:S}
\mathcal{S}(g;Z_N) := \frac12 \norm{g}{\Y}^2 - \ip{g}{Z_N},\qquad g \in \Y,
\end{equation}
as this is the infinite-dimensional negative log-likelihood in the Cameron-Martin-Girsanov sense, cf. \cite{w03}. This leads to the following method modification of the cIRGNM in case of random noise:
\begin{align}
\hat{u}_{n+1} := \argmin_{u \in \X} \bigg[\mathcal S \left(\taylor{\hat{u}_n}{u}; Z_N\right) + \alpha_n \norm{u - \hat{u}_0}{\X}^2\bigg]. \label{eq_IRGNM}
\end{align}
Note that all terms in \eqref{eq_IRGNM} - especially the term $\ip{Z_N}{\taylor{\hat{u}_n}{u}}$ - are well-defined, since we have $\taylor{\hat{u}_n}{u} \in \Y$ for all $u \in \X$. This method has been proposed and analyzed in \cite{HW13}, and allows for further generalizations including different noise models or general convex penalty terms.

In this work we take a different focus motivated by many practically problems, for which one does not want to begin the reconstruction procedure until the (final) measurement $Y_N$ has been collected (so that $Z_N$ can be computed). Instead, it would be preferable to start the reconstruction immediately after obtaining $Y_1$ and update our estimate of the unknown on the fly as the new observations $Y_2, Y_3, ...Y_{N}$ become available. This motivation aligns with the aim of online algorithms for (linear) inverse problems which have been recently attracted much attention to solve filtering and data assimilation problems \cite{DLC2018,ILS13,ILLS17}.
%To this end, following \cite{ILS13} we rewrite the sequential noisy observations (\ref{eq:model_dyn}) into an artificial dynamic form
%\begin{align}\label{eq:model_dyn2}
%\left\{
%\begin{array}{l}
%u_{n+1}  =u_{n} \\
%Y_{n+1} = \op\left(u_{n+1}\right) + \sigma \xi_{n+1},
%\end{array}
%\right.
%\end{align}
%for $n=0,1\ldots$ and $u_0 = \udag$. Such an artificial dynamic contains a steady state equation associated with the unknown variable $\udag$ and the other observation equation with sequential observation $\{Y_n\}_{n=1,\ldots}$.
%Meanwhile, inserting the offline averaging setting (\ref{eq:avg_noise}) into the artificial dynamic (\ref{eq:model_dyn2}), we can generate an online averaging observation at each artificial time index $n+1$. For instance, we could rewrite the artificial dynamic \eqref{eq:model_dyn2} into an augmented form, such that
%\begin{align}\label{eq:model_dyn2_modified}
%\left\{
%\begin{array}{l}
%u_{n+1}  =u_{n}, \\
%Y_{n+1} = \op\left(u_{n+1}\right) + \sigma \xi_{n+1}, \\
%Z_{n+1} : =  \frac{1}{n+1}(nZ_{n}+Y_{n+1}) = \frac{1}{n+1}\sum_{i=1}^{n+1} Y_i,
%\end{array}
%\right.
%\end{align}
%for $n=0,1\ldots$, $u_0 = \udag$ and $Z_0=0$.
%Both artificial dynamics \eqref{eq:model_dyn2} and \eqref{eq:model_dyn2_modified} then allows us to implement some online reconstruction algorithms to recover the unknown steady state variable $u_{n+1}$ from the sequential observation $\{Y_{i}\}_{i=1}^{n+1}$ or $\{Z_{i}\}_{i=1}^{n+1}$.

In the context of the IRGNM, we propose to sequentially incorporate data by the following scheme which we call the \textbf{dynamic iteratively regularized Gauss-Newton method (dIRGNM)}:
\begin{align}
\hat{u}_{n+1} := \argmin_{\hat{u} \in \X} \bigg[\mathcal S\left(\taylor{\hat{u}_n}{\hat{u}}; Z_n\right) + \alpha_n \norm{\hat{u} - \hat{u}_0}{\X}^2\bigg]. \label{eq_DIRGNM}
\end{align}
Note that this algorithm can be started as soon as $Y_1$ (and hence $Z_1$) is available, i.e. right after the first set of observations are collected. We emphasize that the main difference between \eqref{eq_IRGNM} and \eqref{eq_DIRGNM} is the index $n$ in the used data $Z_n$ (compared to $Z_N$ in \eqref{eq_IRGNM}). However, this ensures that the data $Y_n$ (and hence $Z_n$) that is currently available are assimilated sequentially via \eqref{eq:avg_noise} into the algorithm \eqref{eq_DIRGNM}. Despite of such a subtle modification, we show that the proposed scheme will allows us to immediately benefit from the decreasing uncertainty which will, in turn, lead to faster computations of the unknown without compromising accuracy.

\subsection{Literature overview}

In many real-world application areas, it is common to have experimental settings that allow us to sequentially acquire multiple observations of the physical process under consideration (e.g. by repeating the experiment). The classical approach for solving this kind of inverse problem is to first produce the average of those observations, and use this average with a standard regularization method to infer the unknown quantity/parameter of interest. A class of methods for solving ill-posed inverse problems is the so-called variational regularization which includes the well-known Tikhonov regularization as well as various other methods such as Landweber iteration, steepest descent and $\nu$-methods \cite{EHN96,KNS08,LP13}.

The analysis of the convergence of most existing iterative methods, including those cited above, assume that observed data remain the same throughout the iterative procedure. However, exploring sequential variants of these methods in which data are updated as they become available can bring substantial benefits in practical settings. The focus on the IRGNM is particularly relevant since, for data assimilation problems, the Gauss-Newton method has been shown to have striking similarities with Kalman filtering methodologies that sequentially update parameters and states of dynamical processes \cite{BMB94,BC93,CCS21,CT21,GLN07}. 

The extensive and successful use of Kalman filter methods for large-scale data assimilation applications such as ocean and weather forecasting \cite{CB12,LLJ15,MW06}, has prompted a body of work aimed at importing and adapting those methodologies for solving-ill posed inverse problems. In \cite{DLC2018,ILLS17}, for example, regularization theory was used to analyze convergence of data assimilation algorithms, such as the Kalman filter, 3DVAR and 4DVAR in the context of solving linear inverse problems. These works have shown that using multiple instances of noisy observations lead to more robust and stable algorithms when a scaling regularization parameter is appropriately tuned. In the nonlinear case, however, whether the convergence of filtering methods, such as the ensemble Kalman filter \cite{GE09,GE94} and extended Kalman filter, can be improved by using multiple instances of data is still an open problem. Our work on the dIRGNM, in addition to providing practical algorithms that can invert data sequentially, will also pave the way towards understanding the dynamic behavior of data assimilation algorithms for nonlinear inverse problems.

\subsection{Aim of the paper}

\textcolor{black}{Our primary focus and contribution from this work is the development and understanding of the dIRGNM, which, as stated earlier, is a modified version of the IRGNM that enable us to sequentially invert observed data. We propose two particular forms of a dynamic IRGNM, the first is given above in \eqref{eq_DIRGNM} which is intended for our analysis with infinitely many observation. The second form, which we refer to as, the hybrid iterated regularized Gauss-Newton method (hIRGNM) combines the classical (cIRGNM) with the dIRGNM in the practical case when finitely many observations are available. The motivation behind the hybrid scheme is to obtain improved performance by initially running the dIRGNM for various but finitely many observations, followed by running the cIRGNM with the average of all acquired observations. Based on standard assumptions for nonlinear inverse problems, we prove well-definedness for both algorithms. In addition, we derive appropriate error bounds and convergence rates. In order to prevent from data over-fitting, our analysis also includes recommended choices for the parameter $\alpha_n$. We employ two PDE-constrained parameter identification problems in order to numerically test the convergence results of the proposed dIRGNM and hIRGNM, as well as to demonstrate their computational advantages over the cIRGNM.}

%We emphasis that to the best of the authors knowledge, the results derived in this work are new and have not been considered before. Specifically
% the theoretical understanding of nonlinear statistical inverse problem methodologies in the context of an artificial dynamics process, of the form \eqref{eq:model_dyn2} or \eqref{eq:model_dyn2_modified},
%  rather than one single instance of the data.

\subsection{Outline}
The outline of this paper is as follows. In Section \ref{sec:analysis} we provide the necessary background and material
related to the dIRGNM and assumptions, in order to carry out our analysis where we derive generic error bounds. This will lead into Section \ref{sec:error} where we discuss
 and present convergence analysis, with error bounds, of each of the various algorithms introduced which include the dIRGNM and the hIRGNM. We also present and discuss the implementation of each method.
 In order to verify such results we present numerical experiments in Section \ref{sec:num}, where we provide tests on three PDE-constrained parameter identification problems motivated from practical applications including the characterization of geological properties of the subsurface as well as medical imaging.
Finally in Section \ref{sec:conc} we conclude our findings, and present potential new directions of research.

%\section{Preliminaries}
%\label{sec:prelim}
%
%\todo[inline]{Maybe we do not need this section anymore, as the artifical dynamics will not be used anywhere.}
%

\section{Standing assumptions and error analysis}
\label{sec:analysis}

In this section, we provide error bounds for IRGNM \eqref{eq_IRGNM} in the general data model and introduce the assumptions needed. Note that - due to the only difference in the used data - the same bounds also apply for \eqref{eq_DIRGNM}. To treat both cases at the same time, let us denote by $W \in \{Z_N,Z_n,Y_n\}$ the available data, define
\begin{equation}\label{eq:IRGNM0}
\mathcal{J}[u,\hat{u}_n,\hat{u}_0,\alpha_{n},W] := \mathcal S \left(\taylor{\hat{u}_n}{u}; W\right) + \alpha_n \norm{u - \hat{u}_0}{\X}^2,	
\end{equation}
and consider
\begin{equation}\label{eq:IRGNM}
\hat{u}_{n+1} := \argmin_{u \in \X} \mathcal{J}[u,\hat{u}_n,\hat{u}_0,\alpha_{n},W].
\end{equation}
If $W = Z_N$, this equals \eqref{eq_IRGNM}, and if $W = Z_n$, then this equals \eqref{eq_DIRGNM}. Furthermore, the following analysis will also illustrate why the naive choice $W = Y_n$ does not allow for an assimilation of the sequential data \eqref{eq:model_dyn} and will not provide a convergent algorithm unless the noise vanishes.

Our analysis here closely follows the general approach to error bounds for variational regularization methods described in \cite{HW22}. %For brevity we additionally introduce
%\begin{align*}
%\mathcal{T}(g;y^{\dagger})&= \frac{1}{2}\|g - y^{\dagger}\|^2_{\mathcal{Y}}.
%\end{align*}

\subsection{Preliminary error decomposition}
Let us assume that the $n$-th iterate $\hat{u}_n \in D(F)$ is well defined. As a first step, we aim to provide an error bound for $\hat{u}_{n+1}$ defined by \eqref{eq:IRGNM}. The minimality condition of \eqref{eq:IRGNM} implies
\begin{multline}
\alpha_n \left[ \norm{\hat{u}_{n+1} - \hat{u}_0}{\X}^2 - \norm{\udag - \hat{u}_0}{\X}^2\right] \\\leq \cS{\taylor{\hat{u}_n}{\udag}} - \cS{\taylor{\hat{u}_n}{\hat{u}_{n+1}}}. \label{eq_min_analysis_1}
\end{multline}
Introducing the effective noise level
\[
\err\left(g\right):=  \frac{1}{2}\|g - y^{\dagger}\|^2_{\mathcal{Y}} - \left(\cS{g} - \cS{\ydag}\right), \qquad g \in \Y,
\]
we rewrite the right-hand side of (\ref{eq_min_analysis_1}) by
\begin{align*}
&\cS{\taylor{\hat{u}_n}{\udag}} - \cS{\taylor{\hat{u}_n}{\hat{u}_{n+1}}}\\
=& \left(\cS{\taylor{\hat{u}_n}{\udag}} - \cS{\ydag}\right)\\
&- \left(\cS{\taylor{\hat{u}_n}{\hat{u}_{n+1}}} -\cS{\ydag}\right)\\
=& \frac12 \norm{\taylor{\hat{u}_n}{\udag} - \ydag}{\Y}^2 - \err\left(\taylor{\hat{u}_n}{\udag}\right) \\
& - \frac12 \norm{\taylor{\hat{u}_n}{\hat{u}_{n+1}} - \ydag}{\Y}^2 + \err\left(\taylor{\hat{u}_n}{\hat{u}_{n+1}}\right).
\end{align*}
Note that in each of the settings $W \in \left\{Z_N, Z_n, Y_n\right\}$ we can derive an explicit formulation for $\err(g)$, namely
\begin{equation}\label{eq:error}
\err\left(g\right) = \begin{cases} \frac{\sigma}{N} \sum_{i=1}^N \ip{\xi_i}{g-\ydag} & \text{if }W = Z_N,\\
\frac{\sigma}{n} \sum_{i=1}^n \ip{\xi_i}{g-\ydag} & \text{if }W = Z_n,\\
\sigma \ip{\xi_n}{g-\ydag} & \text{if }W = Y_n,
\end{cases}
\end{equation}
for $g \in \Y$.
If we introduce 
\begin{equation}\label{eq:error2}
\lambda_n = \begin{cases} \frac{\sigma}{N}& \text{if }W = Z_N,\\
	\frac{\sigma}{n} & \text{if }W = Z_n,\\
	\sigma & \text{if }W = Y_n,
\end{cases} \qquad\text{and}\qquad \Xi_n = \begin{cases} \sum_{i=1}^N \xi_i  \stackrel{\mathcal D}{=} \sqrt{N} \xi_1 & \text{if }W = Z_N,\\
\sum_{i=1}^n \xi_i  \stackrel{\mathcal D}{=} \sqrt{n} \xi_1 & \text{if }W = Z_n,\\
\xi_n& \text{if }W = Y_n,
\end{cases}
\end{equation}
where we used that the noises $\xi_i$ are independently identical distributed,
then this gives
\begin{align}\label{eq:errn}
\err_n &:= \err\left(\taylor{\hat{u}_n}{\hat{u}_{n+1}}\right) - \err\left(\taylor{\hat{u}_n}{\udag}\right), \\ \nonumber
&=\lambda_n \langle \Xi_n,F(\hat{u}_n) + F'[\hat{u}_n](\hat{u}_{n+1} - \hat{u}_n)) - \langle \Xi_n, \taylor{\hat{u}_n}{\udag} \rangle, \\ \nonumber
&= \lambda_n \langle \Xi_n, F'[\hat{u}_n](\hat{u}_{n+1} - u^{\dagger}\rangle,
\end{align}
and we obtain by (\ref{eq_min_analysis_1}) that
\begin{multline}
\alpha_n \left[ \norm{\hat{u}_{n+1} - \hat{u}_0}{\X}^2 - \norm{\udag - \hat{u}_0}{\X}^2\right]  + \frac12  \norm{\taylor{\hat{u}_n}{\hat{u}_{n+1}} - \ydag}{\Y}^2 \\\leq \err_n + \frac12 \norm{\taylor{\hat{u}_n}{\udag} - \ydag}{\Y}^2. \label{eq_min_analysis_2}
\end{multline}
To proceed further, we need the following variational source condition, which has been first formulated in \cite{hkps07} and has become a standard assumption in the analysis of variational regularization methods. In many situations it turns out that variational source conditions are necessary and sufficient for convergence rates \cite{hw17}. Note that - as typical for source conditions in general - the smoothness of $\udag$ is therein measured relative to the smoothing properties of $F$.

\begin{ass}[Variational source condition]\label{ass:vsc}
There exists a concave index function $\varphi$ (i.e. $\varphi(0) = 0$ and $\varphi$ and monotonically increasing) such that for all $u \in D(F)$ it holds
\begin{equation}\label{eq:vsc}
\norm{u-\udag}{\X}^2 \leq \norm{u-\hat{u}_0}{\X}^2 - \norm{\udag - \hat{u}_0}{\X}^2 + \varphi \left(\frac12 \norm{\op(u) - \op\left(\udag\right)}{\Y}^2\right).
\end{equation}
\end{ass}

Plugging Assumption \ref{ass:vsc} into (\ref{eq_min_analysis_2}) with $u = \hat{u}_{n+1}$ yields
\begin{multline}\label{eq:aux}
\alpha_n \norm{\hat{u}_{n+1} - \udag}{\X}^2 + \frac12  \norm{\taylor{\hat{u}_n}{\hat{u}_{n+1}} - \ydag}{\Y}^2 \\\leq \err_n + \alpha_n\varphi \left(\frac12 \norm{\op\left(\hat{u}_{n+1}\right) - \op\left(\udag\right)}{\Y}^2\right) + \frac12 \norm{\taylor{\hat{u}_n}{\udag} - \ydag}{\Y}^2.
\end{multline}
In order to further treat the nonlinearity, we employ the
following assumption.
\begin{ass}[Tangential cone condition]\label{ass:tcc}
There exists a constant $\Ctc \geq 1$ and $\eta > 0$ sufficiently small such that
\begin{align*}
&\frac{1}{\Ctc} \norm{\op(v) - \ydag}{\Y}^2 - \eta \norm{\op(u) - \ydag}{\Y}^2 \\
\leq & \norm{\taylor{u}{v} - \ydag}{\Y}^2\\
\leq & \Ctc \norm{\op(v) - \ydag}{\Y}^2 + \eta \norm{\op(u) - \ydag}{\Y}^2.§
\end{align*}
\end{ass}
\begin{rem}
This tangential cone condition follows from the standard tangential cone condition with some $\Ctc$, see \cite[Lemma 5.2]{HW13}. If $\varphi \geqslant \sqrt{t}$ as $t \to 0$, it can - using the techniques from \cite{FW15} - be replaced by a Lipschitz-type assumption.
\end{rem}

The tangential cone condition gives for the second term on the left-hand side of \eqref{eq:aux} that
\[
\frac{1}{\Ctc} \norm{\op\left(\hat{u}_{n+1}\right) - \ydag}{\Y}^2 - \eta \norm{\op\left(\hat{u}_n\right) - \ydag}{\Y}^2 \leq \norm{\taylor{\hat{u}_n}{\hat{u}_{n+1}} - \ydag}{\Y}^2,
\]
and for the third term on the right-hand side, with (\ref{eq:ip}), that
\begin{align*}
\norm{\taylor{\hat{u}_n}{\udag}-\ydag}{\Y}^2 &\leq \Ctc \norm{\op\left(\udag\right) - \ydag}{\Y}^2 + \eta \norm{\op\left(\hat{u}_n\right) - \ydag}{\Y}^2 \\&= \eta \norm{\op\left(\hat{u}_n\right) - \ydag}{\Y}^2.
\end{align*}
Inserting above two inequalities into (\ref{eq:aux}) we obtain the recursive error estimate
\begin{multline}
\alpha_n \norm{\hat{u}_{n+1} - \udag}{\X}^2 + \frac{1}{2\Ctc} \norm{\op\left(\hat{u}_{n+1}\right) - \ydag}{\Y}^2 \\\leq \err_n + \alpha_n\varphi \left(\frac12 \norm{\op\left(\hat{u}_{n+1}\right) - \op\left(\udag\right)}{\Y}^2\right) + \eta \norm{\op\left(\hat{u}_n\right) - \ydag}{\Y}^2. \label{eq_min_analysis_3}
\end{multline}
We then abbreviate
\begin{align*}
d_n &:= \norm{\hat{u}_n - \udag}{\X}^2,\\
t_n &:= \frac12\norm{\op\left(\hat{u}_n\right) - \ydag}{\Y}^2,
\end{align*}
and estimate
\[
\varphi\left(t\right) - \frac{t}{\alpha} \leq \sup_{\tau \geq 0} \left[\varphi\left(\tau\right) - \frac{\tau}{\alpha}\right] = \left(-\varphi\right)^* \left(-\frac{1}{\alpha}\right) =: \Psi\left(\alpha\right),
\]
for $t = \frac 12 t_{n+1}$ with the Fenchel conjugate $\left(-\varphi\right)^*$ of the convex function $-\varphi$. Therewith, we have proven the following:
\begin{Lem}[Preliminary error estimate]
\label{lem:free}
Let Assumptions \ref{ass:vsc} and \ref{ass:tcc} hold and assume that $\hat{u}_n \in D(F)$ is well defined. Then the next iterate $\hat{u}_{n+1}$ defined by \eqref{eq:IRGNM} satisfies a.s. a preliminary error decomposition
\begin{equation}\label{eq:aux2}
	\alpha_n d_{n+1} + \frac{1}{2\Ctc} t_{n+1} \leq \err_n + \alpha_n \Psi \left(2 \Ctc \alpha_n\right) + 2 \eta t_n,
\end{equation}
with
\[
\err_n := \lambda_n\ip{\Xi_n}{\op'\left[\hat{u}_n\right]\left(\hat{u}_{n+1} - \udag\right)},
\]
with $\lambda_n$ and $\Xi_n$ as in \eqref{eq:error2}.
\end{Lem}

The recursive error estimate \eqref{eq:aux2} is similar to those obtained for the standard IRGNM, see e.g. \cite{HW13}. Before we continue, let us discuss this error estimate briefly. In case of noise-free observations $Y_n = \op\left(\udag\right)$, we have $\xi_n = 0$ and hence $\err_n = 0$. This shows that the first term in \eqref{eq:aux2} is in fact a (preliminary) \textbf{propagated data noise error}. Correspondingly, if $F$ is linear, we have $\eta = 0$ and hence the third term in \eqref{eq:aux2} is an upper bound for the \textbf{nonlinearity error}. The remaining second term $\alpha_n \Psi \left(2 \Ctc \alpha_n\right)$ in \eqref{eq:aux2} is a bound for the \textbf{approximation error}, which can clearly be made arbitrarily small by letting $\alpha_n \searrow 0$. In the following we will discuss the individual error terms.

\subsubsection{The approximation error}

To simplify the approximation error bound \\ $\alpha_n \Psi \left(2 \Ctc \alpha_n\right)$, we need an additional restriction on the source function $\varphi$:
\begin{ass}[Saturation of source functions]\label{ass:varphi}
	For the function $\varphi$ from the variational source condition \eqref{eq:vsc} there exists an $\epsilon > 0$ such that $\varphi^{1+\epsilon}$ is concave.
\end{ass}

Assumption \ref{ass:varphi} implies that
\[
\varphi\left(C \tau\right)^{1+\epsilon} = \varphi\left(C \tau + \left(1-C\right)0\right)^{1+\epsilon} \leq C \varphi\left(\tau\right)^{1+\epsilon} + \left(1-C\right) \varphi \left(0\right)^{1+\epsilon} = C \varphi\left(\tau\right)^{1+\epsilon},
\]
for all $\tau \geq 0$ and $C \geq 1$, and thus the monotonicity of $\varphi$ yields
\[
\varphi\left(C \tau\right) \leq \max\left\{1, C^\frac{1}{1+\epsilon}\right\} \varphi\left(\tau\right) \qquad\text{for all}\qquad C, \tau \geq 0.
\]
For the approximation error function $\Psi$, this implies
\begin{equation}\label{eq:Psi_C}
\Psi \left(C \alpha\right) = \sup_{\tau \geq 0} \left[\varphi\left(\tau\right) - \frac{\tau}{C \alpha}\right] = \sup_{s \geq 0} \left[\varphi\left(C^{\frac{1+\epsilon}{\epsilon}}s\right) - \frac{C^{\frac{1}{\epsilon}}s}{\alpha}\right] \leq \max\left\{1, C^{\frac{1}{\epsilon}}\right\} \Psi \left(\alpha\right).
\end{equation}
Consequently, under Assumption \ref{ass:varphi} we can simplify \eqref{eq:aux2} to
\[
\alpha_n d_{n+1} + \frac{1}{2\Ctc} t_{n+1} \leq \err_n + \left(2 \Ctc \right)^{\frac{1}{\varepsilon}}\alpha_n \Psi \left(\alpha_n\right) + 2 \eta t_n.
\]

\subsubsection{The nonlinearity error}

The nonlinearity error term $2\eta t_n$ can in principle be bounded by the other two error contributions by the help of the following abstract lemma:
\begin{Lem}\label{lem:nl}
	Let $a_n,b_n$ be two sequences such that
	\begin{equation}\label{eq:an}
		a_{n+1} \leq C \eta a_n + b_n \qquad\text{for all}\qquad n \in \mathbb N_0,
	\end{equation}
	with a constants $C > 0, 0 < \eta < C^{-1}$. If $a_0 \leq \frac{\eta}{1-\eta C} b_0$, then this implies
	\[
	a_{n+1} \leq \gamma b_n \qquad\text{for all}\qquad n \in \mathbb N_0,
	\]
	with $\gamma := \frac{1}{1-\eta C}$.
\end{Lem}
\begin{proof}
	We perform an induction over $n$. For $n =0$ we have
	\[
	a_1 \leq C a_0 + b_0  \leq \frac{\eta C}{1-\eta C} b_0 + b_0 = \left(\gamma -1\right) b_0 + b_0 = \gamma b_0,
	\]
	by assumption. For the induction step, we insert the induction hypothesis into \eqref{eq:an} and obtain
	\[
	a_{n+1} \leq C \eta a_n + b_n \leq C \eta \gamma b_n + b_n = \left(1+C \eta \gamma\right)b_n,
	\]
	and thus the claim is shown as soon as we prove
	\[
	\left(1+C \eta \gamma\right) \leq \gamma.
	\]
	But the latter is true if and only if
	\[
	\gamma \geq \frac{1}{1-C \eta},
	\]
	which holds by the definition of $\gamma$.
\end{proof}

We will illustrate this exemplarily in case of noise-free observations $Y_n = \op \left(\udag\right)$ at hand, i.e. $\xi_n = 0$ and hence $\err_n = 0$. In this case, \eqref{eq:aux2} reduces to
\begin{equation}\label{eq:aux3}
\alpha_n d_{n+1} + \frac{1}{2\Ctc} t_{n+1} \leq \left(2\Ctc\right)^{\frac{1}{\varepsilon}}\alpha_n \Psi \left(\alpha_n\right) + 2 \eta t_n.
\end{equation}
Applying Lemma \ref{lem:nl} to this inequality after neglecting the $\alpha_{n}d_{n+1}$ term on the left-hand side, we obtain - assuming that $t_0$ and $\eta$ are sufficiently small - the bound
\begin{equation}\label{eq:nl_error}
t_n  \leq \frac{\left(2\Ctc\right)^{1+\frac{1}{\varepsilon}}}{1-4\Ctc\eta}\alpha_{n-1} \Psi \left(\alpha_{n-1}\right),
\end{equation}
for the nonlinearity error. To derive a concrete (noise-free) convergence rate from this, we need to pose an additional assumption on the decay rate of the $\alpha_n$ as mentioned in Section \ref{se1}.
\begin{ass}[Regularization parameters]\label{ass:regpar}
	The regularization parameters $\alpha_n$ are chosen such that
	\[
    \alpha_0 \leq 1, \quad  \alpha_n \searrow 0, \quad  1 \leq \frac{\alpha_n}{\alpha_{n+1}} \leq \Cdec \qquad\text{for all}\qquad n \in \mathbb N.
	\]
\end{ass}
This assumption allows us to compare $\alpha_{n-1}$ (and $\Psi\left(\alpha_{n-1}\right)$) with $\alpha_n$ (and $\Psi\left(\alpha_n\right)$). We derive from \eqref{eq:nl_error} that
\begin{align*}
t_n & \leq \frac{\left(2\Ctc\Cdec\right)^{1+\frac{1}{\varepsilon}}}{1-4\Ctc\eta}\alpha_{n} \Psi \left(\alpha_{n}\right), \nonumber \\
d_n & \leq \left(1+ \frac{4 \Ctc}{1-4\Ctc\eta}\right) \left(2\Ctc\Cdec\right)^{\frac{1}{\epsilon}} \Psi \left(\alpha_{n}\right), 
\end{align*}
and by plugging these into \eqref{eq:aux3}, we immediately obtain:
\begin{thm}[Convergence rates for noise-free data]
Let Assumptions \ref{ass:vsc}-\ref{ass:regpar} hold and suppose that each $\hat{u}_n \in D(F)$ is well-defined. If $\eta$ and $t_0$ are sufficiently small, we have the convergence rates
\[
d_{n} = \mathcal O \left(\Psi \left(\alpha_n\right)\right)\qquad\text{and}\qquad t_{n} = \mathcal O \left(\alpha_n\Psi \left(\alpha_n\right)\right),
\]
as $n \to \infty$.
\end{thm}

\subsubsection{The propagated data noise error}

The propagated data noise error in \eqref{eq:aux2} is somewhat more difficult to handle, as $\err_n$ depends (implicitly) on $\hat{u}_{n+1}$, which is why \eqref{eq:aux2} should only be considered as a preliminary error estimate. To derive a more helpful bound without implicit dependencies, we have to bound $\err_n$ and factorize over $\hat{u}_{n+1}$.  Therefore as seen in Lemma \ref{lem:free},  $\err_n$ is always of the form $\err_n = \lambda_n \ip{\Xi_n}{g-\ydag}$ with a scalar $\lambda_n$ and some Hilbert space process $\Xi_n$. If we now take a Gelfand triple $\left(\V,\Y,\V'\right)$, where the embedding $\iota : \V\hookrightarrow \Y$ is a Hilbert-Schmidt operator, then this allows for
\begin{equation}
\label{eq:tmp}
\err_n \leq \lambda_n \norm{\Xi_n}{\V'} \norm{F'[\hat{u}_n](\hat{u}_{n+1} - u^{\dagger})}{\V},
\end{equation}
with $\lambda_n$ and $\Xi_n$ as in \eqref{eq:error2}, and
\begin{align*}
 \E\norm{\Xi_n}{\V'}^2  &= \norm{ \E \Xi_n}{\V'}^2  + \mathrm{trace}( \iota \mathrm{Cov}(\Xi_n) \iota^*) \\
 &\leq \norm{ \E \Xi_n}{\V'}^2  + \mathrm{trace}( \iota \iota^*) < \infty,
 \end{align*}
which follows  $\norm{\Xi_n}{\V'} \leq \infty$ a.s..
To bound the second term in \eqref{eq:tmp}, we employ the following assumption.
\begin{ass}[$\V$-smoothingness]\label{ass:smoothing}
There exists a parameter $\theta \in \left(0,1\right)$ and a constant $C_\theta \geq 1$ such that
\[
\norm{\op'\left[\hat{u}_n\right]\left(u-\udag\right)}{\V} \leq C_\theta  \norm{\op'\left[\hat{u}_n\right]\left(u-\udag\right)}{\Y}^\theta \norm{u-\udag}{\X}^{1-\theta},
\]
for all $u \in \X$.
\end{ass}
\begin{rem}	
This assumption is e.g. satisfied if $\op'\left[\hat{u}_n\right]$ maps Lipschitz continuously into a smoother space than $\V$ which obeys a classical interpolation inequality, see e.g. \cite[Rem. 2.6]{HW22}. Thus, Assumption \ref{ass:smoothing} characterizes in fact the smoothing properties of the forward operator $F$ in terms of its derivative $F'$.
\end{rem}

Together with Young's inequality with $\varepsilon > 0$, i.e.
\begin{equation}\label{eq:young}
ab \leq \varepsilon a^r + \frac{1}{r'} \left(\frac{1}{r\varepsilon}\right)^{\frac{r'}{r}} b^{r'},
\end{equation}
for $a,b \geq 0$ and $r,r' \in \left[1,\infty\right]$ such that $\frac{1}{r} + \frac{1}{r'} = 1$, the $\V$-smoothingness yields
\begin{align*}
	\err_n  &\leq C_\theta \lambda_n \norm{\Xi_n}{\V'} \norm{\op'\left[\hat{u}_n\right]\left(\hat{u}_{n+1} -\udag\right)}{\Y}^\theta \norm{\hat{u}_{n+1} - \udag}{\X}^{1-\theta}\\
	&\stackrel{\eqref{eq:young}}{\leq}  \frac12 \alpha_n d_{n+1}+  \left(\frac{\theta+1}{2}\right) \left(\frac{1-\theta}{2}\right)^{\frac{1-\theta}{1+\theta}} \left(\frac12 \alpha_n\right)^{\frac{\theta-1}{\theta+1}} C_\theta^{\frac{2}{1+\theta}} (\lambda_n \norm{\Xi_n}{\V'})^{\frac{2}{1+\theta}} \\   & \qquad\times \norm{\op'\left[\hat{u}_n\right]\left(\hat{u}_{n+1} -\udag\right)}{\Y}^{\frac{2\theta}{1+\theta}}\\
	&\stackrel{\eqref{eq:young}}{\leq}   \frac12 \alpha_n d_{n+1} +\varepsilon  \norm{\op'\left[\hat{u}_n\right]\left(\hat{u}_{n+1} -\udag\right)}{\Y}^2 +2 ^{-(1+\theta)}\left(\frac{\theta}{\varepsilon}\right)^\theta \left(1-\theta \right)^{1-\theta}  C_\theta^2 \alpha_n ^{\theta-1} \lambda_n^2 \norm{\Xi_n}{\V'}^{2},
\end{align*}
where in the first inequality we used $\varepsilon = 1$, $r = \frac{2}{1-\theta}$ and $r' = \frac{2}{1+\theta}$, and in the second inequality we used $r = \frac{1+\theta}{\theta}$ and $r' = 1+\theta$. Now the term $\norm{\op'\left[\hat{u}_n\right]\left(\hat{u}_{n+1} -\udag\right)}{\Y}^2$ can again be handled by the tangential cone condition, since
\begin{align*}
	\norm{\op'\left[\hat{u}_n\right]\left(\hat{u}_{n+1} -\udag\right)}{\Y}^2 \leq &2 \norm{\taylor{\hat{u}_n}{\hat{u}_{n+1}} - \ydag}{\Y}^2 \\
	&+ 2 \norm{\taylor{\hat{u}_n}{\udag}-\ydag}{\Y}^2\\
	\leq & 2 \Ctc \norm{\op \left(\hat{u}_{n+1}\right) - \ydag}{\Y}^2 + 4 \eta \norm{\op \left(\hat{u}_n\right)- \ydag}{\Y}^2.
\end{align*}
Plugging this in, we obtain the error estimate
\[
\err_n \leq 4 \varepsilon \Ctc t_{n+1} + 8 \varepsilon \eta t_n + C \left(\theta, \varepsilon\right) \lambda_n^2 \alpha_n^{\theta-1}  \norm{\Xi_n}{\V'}^{2}+ \frac12 \alpha_n d_{n+1},
\]
with an arbitrary constant $\varepsilon > 0$ and the constant
\[
C \left(\theta,\varepsilon\right) =  2 ^{-(1+\theta)}\left(\frac{\theta}{\varepsilon}\right)^\theta \left(1-\theta \right)^{1-\theta}  C_\theta^2  .
\]
Choosing $\varepsilon$ as the solution of $4 \varepsilon \Ctc = \frac{1}{4 \Ctc}$, i.e. $\varepsilon := \left(16 \Ctc^2\right)^{-1}$ and combining the above estimate with \eqref{eq:aux2}, we obtain the following.
\begin{Lem}[Total error estimate]
Let Assumptions \ref{ass:vsc}, \ref{ass:tcc} and \ref{ass:smoothing} hold and suppose that $\hat{u}_n \in D(F)$ is well-defined. If $\eta$ and $t_0$ are sufficiently small, then the total error estimate
\begin{align}\label{eq:ree}
	\frac{1}{2}\alpha_n d_{n+1} + \frac{1}{4\Ctc} t_{n+1} &\leq \alpha_n \Psi \left(2 \Ctc \alpha_n\right) + 2\left(1+\frac{1}{4 \Ctc^2}\right) \eta t_n \\&+ C \left(\theta, \frac{1}{16 \Ctc^2}\right) \lambda_n^2  \alpha_n^{\theta-1}  \norm{\Xi_n}{\V'}^{2}, \nonumber
\end{align}
with $\lambda_n$ and $\Xi_n$ as in \eqref{eq:error2} according to the specific data model holds true almost surely. If additionally Assumption \ref{ass:varphi} hold, then the total estimate obeys
\begin{align}\label{eq:ree2}
	\frac{1}{2}\alpha_n d_{n+1} + \frac{1}{4\Ctc} t_{n+1} &\leq \left(2 \Ctc\right)^{\frac{1}{\epsilon}}\alpha_n \Psi \left(\alpha_n\right) + 2\left(1+\frac{1}{4 \Ctc^2}\right) \eta t_n \\&+ C \left(\theta, \frac{1}{16 \Ctc^2}\right) \lambda_n^2  \alpha_n^{\theta-1}  \norm{\Xi_n}{\V'}^{2}, \nonumber
\end{align}
almost surely.
\end{Lem}
The still contained nonlinearity error on the right-hand side of \eqref{eq:ree} can now be handled similar to the noise free-case. Applying Lemma \ref{lem:nl} to \eqref{eq:ree} after neglecting the $\alpha_n d_{n+1}$ term on the left hand side, we obtain in view of \eqref{eq:Psi_C} and Assumption \ref{ass:regpar} - for sufficiently small $t_0$ and $\eta$ - the bound
\[
t_{n+1} \leq \gamma_{\mathrm{nl}} \left(\alpha_n \Psi\left(2 \Ctc\alpha_n\right) + C \left(\theta, \frac{1}{16 \Ctc^2}\right) \lambda_n^2 \alpha_n^{\theta-1}  \norm{\Xi_n}{\V'}^{2}\right),
\]
for all $n$ with
\[
\gamma_{\mathrm{nl}} := \frac{4\Ctc}{1-8\eta \Ctc \left(1+ \frac{1}{4 \Ctc^2}\right)}.
\]
Plugging this into \eqref{eq:ree} and dropping the $t_{n+1}$ term on the left-hand side, we get after division by $\alpha_n$ the following result by adjusting the iteration number accordingly.
\begin{Lem}[Final recursive error estimate]\label{lem:finalerror}
Let Assumptions \ref{ass:vsc}-\ref{ass:smoothing} hold and assume that $\hat{u}_n \in D(F)$ is well-defined. If $t_0$ and $\eta > 0$ are sufficiently small, then we have the error bound	
\begin{align}
d_{n} & \leq C_1 \left(2 \Ctc\right)^{\frac{1}{\epsilon}}  \Psi \left(\alpha_{n-1}\right) + C_2 \lambda_n^2 \alpha_{n-1}^{\theta-2}  \norm{\Xi_n}{\V'}^{2} \label{eq:ree2b} \\
& \leq C_1 \left(2 \Ctc \Cdec\right)^{\frac{1}{\epsilon}} \Psi \left(\alpha_{n}\right) + C_2 \lambda_n^2 \alpha_{n-1}^{\theta-2}  \norm{\Xi_n}{\V'}^{2}, \label{eq:ree3b}
\end{align}
a.s. with constants
\begin{align*}
C_1 & := 2\left(1 + 2 \left(1+\frac{1}{4 \Ctc^2}\right) \eta \gamma_{\mathrm{nl}}\right),\\
C_2 &:= 2 \left(1+2 \left(1+\frac{1}{4 \Ctc^2}\right) \eta \gamma_{\mathrm{nl}}\right) C \left(\theta, \frac{1}{16 \Ctc^2}\right).
\end{align*}
\end{Lem}
The above error estimate \eqref{eq:ree2b} plays an important role in the following analysis and discussion.

\subsection{Well-definedness of the method}

As a first application, we are now in position to prove the well-definedness of the dIRGNM in \eqref{eq:IRGNM}. Recall therefore that we have $W = Y_n$ and hence $\lambda_n = \sigma $ and $\Xi_n = \xi_n$ in this case.
\begin{thm}\label{thm:well_def}
Let Assumptions \ref{ass:vsc}-\ref{ass:regpar} hold and let $\hat{u}_0 \in D(F)$. Assume that there exists an open ball $B_r\left(\udag\right) \subset D(F)$ around $\udag$ in $D(F)$ and that $\norm{\xi_n}{\V'}$ in \eqref{eq:model_dyn} satisfies the deviation inequality
\begin{equation}\label{eq:die}
\mathbb P \left[\left|\norm{\xi_n}{\V'} - \E{\norm{\xi_n}{\V'}}\right| \geq \delta\right] \leq 2\exp\left(- c \delta\right),
\end{equation}
for all $\delta > 0$ with some constant $c > 0$. Suppose furthermore that $t_0$ and $\eta > 0$ are sufficiently small to allow the application of Lemma \ref{lem:nl}. Then, if both $\alpha_{n} > 0$ and $\sigma^2 \alpha_{n}^{\theta-2}$ are sufficiently small, then we have also $\hat{u}_{n+1} \in B_r \left(\udag\right) \subset D(F)$ with probability at least
\[
1- \exp\left(-c\frac{\alpha_n^{\frac{\theta}{2}-1}}{\sigma} \left(\frac{r}{C}-\left(\sqrt{\Psi \left(\alpha_{n}\right)} + \frac{\sigma}{\alpha_n^{\frac{\theta}{2} - 1}} \E{\norm{\xi_n}{\V'}}\right)\right)\right).
\]
\end{thm}
\begin{proof}
Let $X := \norm{\xi_n}{\V'}$. By Jensen's inequality, we have $|\E{X}|^2 \leq \E{|X|^2}$.
As $\E{|X|^2} < \infty$ by the Hilbert-Schmidt property of the embedding $\iota : \V \to \Y$, this shows that $\E{X} < \infty$.
Now suppose that $\hat u_{n} \in D(F)$. Then by Lemma \ref{lem:finalerror}, the error estimate \eqref{eq:ree2b} holds true. Thus, due to $\sqrt{a +b} \leq \sqrt{a} + \sqrt{b}$ and adjusting the iteration number appropriately, we have
\begin{align*}
\norm{ \hat u_{n+1} - \udag}{\X} &\leq C \left[\sqrt{\Psi \left(\alpha_{n}\right)} + \frac{\sigma}{\alpha_n^{\frac{\theta}{2} - 1}} X\right]\\
& \leq C \left[\sqrt{\Psi \left(\alpha_{n}\right)} + \frac{\sigma}{\alpha_n^{\frac{\theta}{2} - 1}} \E{X}\right] + C\frac{\sigma}{\alpha_n^{\frac{\theta}{2} - 1}}\left|X - \E{X}\right|,
\end{align*}
a.s. with some $C > 0$. Now we assume that both $\alpha_n>0$ and $\sigma^2\alpha_n^{\theta-2}$ are sufficiently small such that the first term in above inequality is smaller than $r$ and thus by \eqref{eq:die} we have
\begin{align*}
\mathbb P \left[ \norm{ \hat u_{n+1} - \udag}{\X} \leq r \right] & \geq \mathbb P \left[C\frac{\sigma}{\alpha_n^{\frac{\theta}{2} - 1}}\left|X - \E{X}\right| \leq r-C \left(\sqrt{\Psi \left(\alpha_{n}\right)} + \frac{\sigma}{\alpha_n^{\frac{\theta}{2} - 1}} \E{X}\right) \right] \\
& = 1- \mathbb P \left[\left|X - \E{X}\right| \geq \frac{\alpha_n^{\frac{\theta}{2}-1}}{\sigma} \left(\frac{r}{C}-\left(\sqrt{\Psi \left(\alpha_{n}\right)} + \frac{\sigma}{\alpha_n^{\frac{\theta}{2} - 1}} \E{X}\right)\right)\right] \\
& \geq 1- \exp\left(-c\frac{\alpha_n^{\frac{\theta}{2}-1}}{\sigma} \left(\frac{r}{C}-\left(\sqrt{\Psi \left(\alpha_{n}\right)} + \frac{\sigma}{\alpha_n^{\frac{\theta}{2} - 1}} \E{X}\right)\right)\right),
\end{align*}
which proves the claim.
\end{proof}

\begin{rem}
We provide some extended discussion below.
	\begin{itemize}
		\item A deviation inequality of the form \eqref{eq:die} is e.g. satisfied for Gaussian white noise $\xi_1$, see \cite[Thm. 2.1.20]{gn16}.
		\item The assumption that both $\alpha_{n} > 0$ and $\sigma^2 \alpha_{n}^{\theta-2}$ are sufficiently small is natural, as both terms should tend to $0$ anyway during the iteration.
		\item Note that the stated probability in Theorem \ref{thm:well_def} tends to $1$ as $\alpha_{n} > 0$ and $\sigma^2 \alpha_{n}^{\theta-2}$ tend to $0$. As a corollary, Theorem \ref{thm:well_def} implies that all iterates from a certain $n \in \mathbb N$ on will be well-defined with overwhelming probability under reasonable assumptions, or differently phrased that all iterates will be well-defined with overwhelming probability if the initial parameters $\alpha_0$, $\hat u_0$ are chosen carefully.
        \item Theorem \ref{thm:well_def} shows that well-definedness probability of \eqref{eq:IRGNM} with $W=Z_n$ is larger than that with $W=Y_n$, which highlights the focus of our current work.
	\end{itemize}
\end{rem}

\section{Error bounds for different observation models}
\label{sec:error}
In this section, we provide discussion on the error bounds or asymptotical behavior for the different observation models $W \in \left\{Z_N, Z_n, Y_n\right\}$ based on the recursive error estimate in Lemma \ref{lem:finalerror}. We also provide an in-depth description of our proposed algorithms, for which we will test in the following section after this. Precise algorithmic forms will be given, related to each observation
model.

\subsection{Error bound for the classical IRGNM}\label{sec:single}

Let us start by analyzing the classical IRGNM, i.e. \eqref{eq:IRGNM} with $W = Z_N$. According to \eqref{eq:error2} we have $\lambda_n = \sigma / N$ and $\Xi_n = \sum_{i=1}^N \xi_i  \stackrel{\mathcal D}{=} \sqrt{N} \xi_1 $, so that \eqref{eq:ree2} specializes to
\begin{equation}\label{eq:ree3}
	d_{n} \leq C_1 (2\Ctc\Cdec)^{\frac{1}{\epsilon}} \Psi(\alpha_{n}) + C_2\norm{\xi_1}{\V'}^{2} \frac{\sigma^2}{N\alpha_{n}^{2-\theta}},
\end{equation}
where $\norm{\xi_1}{\V'}^{2}$ can be handled as a (random) constant.

To determine an optimal regularization parameter $\alpha_N$ for (\ref{eq:ree3}), let us informally search for the infimal value
\[
\inf_{\alpha > 0} \left[ \Psi \left(\alpha\right) +   \frac{\sigma^2}{ \alpha^{2-\theta}} \right] = \inf_{\alpha > 0} \left[\left(-\varphi\right)^* \left(- \frac{1}{\alpha}\right) +  \frac{\sigma^2}{ \alpha^{2-\theta}} \right].
\]
If we set
\[
\Psi_\theta(t) := \left(-\varphi\right)^*\left(-t^{\frac{1}{2-\theta}}\right),
\]
then we can compute this infimum by means of Fenchel duality as
\begin{align*}
\inf_{\alpha > 0} \left[ \Psi \left(\alpha\right) +  \frac{\sigma^2}{ \alpha^{2-\theta}} \right] & = \inf_{\alpha > 0} \left[ \Psi_\theta \left(\frac{1}{\alpha^{2-\theta}}\right) +  \frac{\sigma^2}{ \alpha^{2-\theta}} \right]	\\
& = \inf_{\alpha' > 0} \left[ \Psi_\theta \left(\alpha'\right) +  \alpha'\sigma^{2}\right]\\
& = - \sup_{\alpha'>0} \left[- \Psi_\theta  \left(\alpha'\right) -  \alpha'\sigma^{2} \right]\\
& = \left(-\Psi_\theta\right)^* \left(- \sigma^{2}\right),
\end{align*}
and by the equality condition in Young's inequality, this infimum is attained for $\frac{1}{\alpha^{2-\theta}} \in \partial \left(-\Psi_\theta\right)^* \left(- \sigma^{2}\right)$.

Now we define the stopping criterion as
\begin{equation}\label{eq:nstar}
	n_* = \min \left\{n \in \mathbb N ~ \,\, \Big| \,\, \frac{1}{\alpha_n^{2-\theta}} \geq \partial \left(-\Psi_\theta\right)^* \left(- \frac{\sigma^{2}}{N}\right) \right\}.
\end{equation}
and obtain the following
\begin{thm}\label{thm:fixed}
	Let Assumptions \ref{ass:vsc}-\ref{ass:smoothing} hold. Assume that all iterates $\hat{u}_n$ of the classical IRGNM \eqref{eq_IRGNM} are well defined with probability larger than $1-\delta$ and that $\norm{F\left(\hat{u}_0\right) - \ydag}{\Y}$ and $\eta$ are sufficiently small.
	Then for the single fixed observation model we obtain
	\[
	\norm{\hat{u}_{n_*} - \udag}{\X}^2  = \mathcal O \left( \left(-\Psi_\theta\right)^* \left(- \frac{\sigma^{2}}{N}\right)\right),
	\]
	with probability larger than $1-\delta$ and $n_*$ chosen by (\ref{eq:nstar}).
\end{thm}

The proof is straightforward by the proposed parameter choice rule and we skip it here. In particular, Theorem \ref{thm:fixed} can be viewed as an extended convergence rate result compared with the cIRGNM in the deterministic setting \cite{KNS08}.

\subsection{Asymptotical analysis for infinitely many (averaged) observations}\label{se:average}
This part focuses on the most interesting case with infinitely many (averaged) observations, i.e. an infinite sequence of observations $Y_n$ as in \eqref{eq:model_dyn}.

Note that in case of $W = Y_n$, i.e. $\lambda_n = \sigma$ and $\Xi_n = \xi_n$, the recursive error bound \eqref{eq:ree3b} contains the terms $\sigma^2$ and $\norm{\xi_n}{\V'}^{2}$, which will in general not tend to $0$ this shows that no convergence (and hence no assimilation) can be expected from the corresponding scheme.
Meanwhile, if we consider \eqref{eq_DIRGNM}, where $\lambda_n = \sigma/n$ and $\Xi_n = \sum_{i=1}^n \xi_n \stackrel{\mathcal D}{=} \sqrt{n} \xi_1$, the recursive error bound \eqref{eq:ree3b} becomes
\begin{align}
d_{n} & \leq C_1 \left(2 \Ctc \Cdec\right)^{\frac{1}{\epsilon}} \Psi \left(\alpha_{n}\right) + C_2\norm{\xi_1}{\V'}^{2} \frac{\sigma^2}{n \alpha_n^{2-\theta}}, \label{eq:ree4}
\end{align}
where $\norm{\xi_1}{\V'}^{2}$ can be handled as a (random) bounded constant.

Thus as a central result, we obtain the following:
\begin{thm}\label{thm:infiniteobs}
	Let Assumptions \ref{ass:vsc}, \ref{ass:tcc} and \ref{ass:smoothing} hold, suppose that $\norm{F\left(\hat{u}_0\right) - \ydag}{\Y}$ and $\eta$ are sufficiently small, and assume that all iterates $\hat{u}_n$ are well-defined a.s.. If $\alpha_n$ is chosen such that
	\[
	\alpha_n \searrow 0\qquad\text{and}\qquad n \alpha_n^{2-\theta} \nearrow \infty,
	\]
	then for averaged observations we have $d_n \to 0$ as $n \to \infty$ such that dIRGNM converges a.s. for infinitely many averaged observations.

If additionally Assumptions \ref{ass:varphi}-\ref{ass:regpar} hold and we choose the regularization parameter $\alpha_n$ such that
\[
\frac{1}{\alpha_n^{2-\theta}} \in \partial \left(-\Psi_\theta\right)^* \left(- \frac{\sigma^2}{n}\right),
\]
then there holds
\begin{align*}
\norm{\hat{u}_{n} - \udag}{\X}^2  = \mathcal O \left( \left(-\Psi_\theta\right)^* \left(- \frac{\sigma^{2}}{n}\right)\right).
\end{align*}
\end{thm}
\begin{proof}
	As all iterates are well-defined by assumption, the first result now follows immediately noticing both terms in (\ref{eq:ree4}), by adopting to the proposed parameter choice rule, vanish when $n\rightarrow \infty$.

Concerning the second results, the proposed parameter choice rule then allows us to obtain
\begin{align*}
d_{n} \leq  C \left(-\Psi_\theta\right)^* \left(- \frac{\sigma^{2}}{n}\right),
\end{align*}
which proves the claim.
%We end the proof noticing that
%\begin{align*}
%\E d_{n} & = \sum_{k=1}^{\infty} \P [E_k\backslash E_{k-1}] \E [d_{n} | E_k \backslash E_{k-1}] \\
% & \leq \sum_{k=1}^{\infty} \P [E_k \backslash E_{k-1}] \max_{E_k} d_{n} \\
% & \leq \frac{C}{c^2} \left(\sum_{k=1}^{\infty} k \exp(-2(k-1)) \right) \left(-\Psi_\theta\right)^* \left(- \frac{\sigma^{2}}{n}\right)\\
% & \leq C \left(-\Psi_\theta\right)^* \left(- \frac{\sigma^{2}}{n}\right)
%\end{align*}
%with a bounded constant $\sum_{k=1}^{\infty} k \exp(-2(k-1))$.
\end{proof}

\begin{rem}\label{rem:dIRGNM}
\begin{enumerate}
	\item The assumption that all iterates are well-defined a.s. is reasonable in view of Theorem \ref{thm:well_def} and can be interpreted as a conditioning on some event with overwhelming probability. To derive \textit{overall} rates of convergence in expectation, one would have to specify what is considered as the reconstruction if $\hat u_n$ is no longer well-defined.
\item Theorem \ref{thm:infiniteobs} yields a qualitative result showing that by \eqref{eq:IRGNM} we can obtain a vanishing asymptotical behavior for the dIRGNM \eqref{eq_DIRGNM} by choosing the regularization parameter appropriately. Though the index $\theta$ might be unknown, we can slightly modify the condition on $\alpha_n$ such that
\[
\alpha_n \searrow 0\qquad\text{and}\qquad n \alpha_n^{2} \nearrow \infty,
\]
are sufficient to guarantee the same result. A natural choice would be $\alpha_n \sim n^{-\beta}$, i.e., for $\beta\in (1/2,1-\theta/2)$ and we will examine the numerical performance in Section \ref{sec:num} for different choices of $\beta$.
\item For H\"older-type source conditions, we have $\varphi\left(t\right) = ct^\nu$ with some $0 \leq \nu < 1$ and $c > 0$. Straight-forward computations show
\[
\left(-\varphi\right)^* \left(-s\right) \sim s^{\frac{\nu}{\nu-1}}, \qquad s > 0
\]
and hence
\[
\Psi \left(\alpha_n\right) \sim \alpha_n^{\frac{\nu}{1-\nu}}.
\]
Then by choosing $\alpha_n \sim n^{-\frac{1-\nu}{2-\nu-\theta(1-\nu)}}$, we obtain an asymptotical decaying rate
\begin{align*}
\norm{\hat{u}_{n} - \udag}{\X}^2 = \mathcal{O}(n^{-\frac{\nu}{2-\nu - \theta (1-\nu)}}).
\end{align*}
\end{enumerate}
\end{rem}

\subsection{Analysis for finitely many averaged observations}

Let us now consider the case that we have access to finite $N \in \mathbb N$ sequential observations. Our aim is to use the dIRGNM iteration. Clearly, the iteration should be stopped after the $N$th iteration, as no further data is available then, and additional iterations should be avoided. In this case, \eqref{eq:ree4} holds true for all $n \leq N$. As a consequence of the above considerations, we obtain the following result:
\begin{thm}\label{thm:hIRGNM}
Let Assumptions \ref{ass:vsc}, \ref{ass:tcc}, \ref{ass:varphi} and \ref{ass:smoothing} hold. Assume that all iterates $\hat{u}_n$ are well defined with probability larger than $1-\delta$ and that $\norm{F\left(\hat{u}_0\right) - \ydag}{\Y}$ and $\eta$ are sufficiently small. Choose the regularization parameters $\alpha_1, ..., \alpha_{N-1}$ arbitrary and $\alpha_{N}$ such that
\[
\frac{1}{\alpha_N^{2-\theta}} \in \partial \left(-\Psi_\theta\right)^* \left(- \frac{\sigma^2}{N}\right).
\]
Then we obtain the final estimate
\[
\norm{\hat{u}_{N} - \udag}{\X}^2  = \mathcal O \left( \left(-\Psi_\theta\right)^* \left(-\frac{\sigma^2}{N} \right)\right),
\]
with probability larger than $1-\delta$.
\end{thm}
%The proof of above theorem is analogously as that of Theorems \ref{thm:fixed} \& \ref{thm:infiniteobs} and we skip the details.
Note that the above result has the same convergence rate as in Theorem \ref{thm:fixed} if we would first collect all data, and then run the cIRGNM once on the averaged data. The advantage of dIRGNM with online outputs can be clearly observed by comparing these two results.
Such advantage also allows us to design a hIRGNM to first run the dIRGNM for the sequential finitely many (averaged) observation $Z_n$ with $n\leq N$ and then move to the cIRGNM for the final averaged observation $Z_N$.

\subsection{Summary of Algorithms}

In this subsection we provide the summary of the cIRGNM and its two proposed variants, i.e. dIRGNM and hIRGNM.

\begin{algorithm}[h]
\caption{Classical iterated regularized Gauss-Newton method (\texttt{cIRGNM})}
\label{alg:cIRGNM}
\SetKwInOut{Input}{inputs}
 \SetKwInOut{Output}{output}
%\begin{algorithmic}[1]
     \Input{$\hat{u}_0$, $u^*$ $\alpha_0$, $M$, $W$, $C_{dec}>1$}
%     $$
%Z_N = N^{-1} \sum_{i=1}^NY_i
%$$

    %\Output{$U_{dIRGN}=\hat{u}_{N+1}$.}
%        Set $s_{0}=0$ \\
     %\KwInput{Input}: $u_0$,  $\alpha_0$\\

\For{$n=1,\ldots,M$}{
Compute
$$\hat{u}_{n} := \argmin_{\hat{u} \in \X} \Bigg[\mathcal S \left(\taylor{\hat{u}_{n-1}}{u}; W\right) + \alpha_n \norm{u - u^*}{\X}^2\Bigg],$$
with $\mathcal S $ defined in (\ref{eq:S}) and $\alpha_{n}=\alpha_{0}C_{dec}^{-n}$.
%$$\mathcal{J}[u,\hat{u}_{n-1},\hat{u}_0,\alpha_{n},W]$$
}
    \Output{$\hat{u}_{M}$.}
\end{algorithm}
The  cIRGNM with generic observations $W$ is displayed in Algorithm \ref{alg:cIRGNM} where we have used the standard choice of regularization parameter $\alpha_{n}$. Usually, one can use $\hat{u}_0 = u^*$ to start the iteration but we keep them differently as we need to do so in order to define the hIRGNM later.
For the purpose of monitoring performance, we select a maximum number of $M$ iterations. However, we recognise that in practice this algorithm needs to be stopped, for example, via the discrepancy principle.

\begin{algorithm}[h]
\caption{Dynamic iterated regularized Gauss-Newton method (\texttt{dIRGNM})}\label{alg:dIRGNM}
\SetKwInOut{Input}{inputs}
 \SetKwInOut{Output}{outputs}
%\begin{algorithmic}[1]
     \Input{$\hat{u}_0$, $\alpha_0$, $N$, $\{Y_i\}_{i=1}^N$, $\beta>0$}
%        Set $s_{0}=0$ \\
     %\KwInput{Input}: $u_0$,  $\alpha_0$\\

\For{$n=1,\ldots,N$}{
(1) Collect data $Y_{n}.$\\
(2) Compute
$$\hat{u}_{n} := \argmin_{\hat{u} \in \X} \Bigg[\mathcal S \left(\taylor{\hat{u}_{n-1}}{u}; Z_{n}\right) + \alpha_n \norm{u - \hat{u}_0}{\X}^2\Bigg],$$
where
$$Z_n = n^{-1} \sum_{i=1}^nY_i, \quad \text{and}\quad \alpha_{n}=\alpha_{0}n^{-\beta}.$$
}
%\end{algorithmic}
    \Output{$U_{dIRGNM}^{N}=\hat{u}_{N}$ and $\alpha_{N}$.}
\end{algorithm}

\begin{algorithm}[h]
\caption{Hybrid iterated regularized Gauss-Newton method (\texttt{hIRGNM})}\label{alg:hIRGNM}
\SetKwInOut{Input}{input}
 \SetKwInOut{Output}{output}
%\begin{algorithmic}[1]
     \Input{$\hat{u}_0$, $\alpha_0$, $N$, $M$, $\{Y_i\}_{i=1}^N$, $C_{dec}>1$, $\beta>0$.}
\textbf{First Part}: Compute
    $$(U_{dIRGNM}^{N}, \alpha_{N})=\texttt{dIRGNM}(\hat{u}_{0},\alpha_{0},N,\{Y_i\}_{i=1}^N,\beta).$$

\textbf{Second Part}: Set $\hat{w}_{0}=U_{dIRGNM}^{N}$ and $\tilde{\alpha}_{0}=\alpha_{N}$.
        Compute
     $$\hat{w}_{n}=\texttt{cIRGNM}(\hat{w}_{0}, \hat{u}_{0}, \tilde{\alpha}_{0},M,Z_{N},C_{dec})$$
      with
     $$
Z_N = N^{-1} \sum_{i=1}^NY_i.
$$

\Output{$\hat{w}_{M}$.}
%\end{algorithmic}
\end{algorithm}

The proposed dIRGNM is summarised in Algorithm \ref{alg:dIRGNM}. We recall that in contrast to the cIRGNM in which the observations are fixed throughout the entire algorithm, the dIRGNM allows us to use observations $Y_{1}, Y_{2},\dots, Y_{n},\dots$ as they become available. More specifically, at each iteration $n$, we use $Z_{n}$ i.e. the average of the $n$ available observations, in order to produce the estimate $\hat{u}_{n}$. While the previous section ensures the asymptotic convergence of the dIRGNM, in practical settings we have only access to limited number of experiments. Therefore, we propose the hybrid version shown in Algorithm \ref{alg:hIRGNM}. The first part of this hybrid IRGNM consists of applying the dIRGNM with $N$ iterations. For the second part we use the  cIRGNM using the final estimate of the dIRGNM as starting point, as well as the average of the $N$ measurements collected upon completion of the dIRGNM. Furthermore, for the second part we choose the regularization parameter $\tilde{\alpha}_{n}=\alpha_{N}C_{dec}^{-n}$ where $\alpha_{N}$ is the final value computed with the dIRGNM.

As discussed in Section \ref{se1}, within the classical setting we would have to wait until all observations are acquired, and use the cIRGNM with the average of all these observations (i.e. with $W=Z_{N}$). However, the numerical experiments from the following section show that the hybrid version can offer significant computational advantages. Indeed, by the time all measurements have been collected and assimilated with the dIRGNM encoded in the hybrid version, the estimate of the unknown already shows good levels of accuracy. Consequently, convergence of the second part of the hIRGNM is then achieved in much fewer iterations than those required by the  cIRGNM. For problems where an iteration of the dIRGNM can be computed within the time-scale of measurement acquisition, faster estimates can be obtained using the hybrid algorithm compared to the classical one.

For all the algorithms we adopt the standard practice of starting the iteration using the same element, $\hat{u}_{0}$, that appears in the stabilization term of the cost functional (\ref{eq:IRGNM0}) that we minimize at each iteration of these algorithms. However, it is worth emphasizing, that for the second part of Algorithm \ref{alg:hIRGNM}, we initialise the iterations using the estimate from dIRGNM while keeping the same initial guess, $\hat{u}_{0}$, in the stabilization term.

\section{Numerical Experiments}
\label{sec:num}
%\todo[inline]{We should discuss our assumptions for the examples we consider here. For the first example we get e.g. that $F$ satisfies the tangential cone condition (reference), and that $F'$ is $H^2$-smoothing by classical regularity theory for elliptic PDEs. The variational source condition can be guaranteed by combining the two facts above and computing an $s$ such that $\udag \in H^s$. For the other examples this is more involved (and maybe not known or true), but we should discuss this.}

In this section, we provide three numerical examples verifying the theoretical finding of current work. Our focus mainly concentrates on the convergence of dIRGNM for infinitely many (averaged) observation, i.e. Theorem \ref{thm:infiniteobs}, and the comparison between hIRGNM and cIRGNM when the same finitely many observation is given, i.e. Theorems \ref{thm:fixed} and \ref{thm:hIRGNM}.

\subsection{Example 1}
In the first benchmark example, the unknown solution $u$ is the potential coefficient of the following PDE
\begin{align}
\label{Exeq1}
-\Delta p + u p &=f, \quad \text{in}~\Omega, \\
p&=g, \quad \text{on}~ \partial \Omega, \nonumber
\end{align}
where $\Omega \subset \mathbb{R}^2$ is a bounded domain with Lipschitz boundary $\partial \Omega$, $f \in L^2(\Omega)$ and $g \in H^{3/2}(\Omega)$. We define the parameter-to-measurements operator $F:  L^2(\Omega) \rightarrow  L^2(\Omega)$ via $p=F(u)$, where $p$ is the unique solution of (\ref{Exeq1}). %The explicit form of the Fr\'echet derivative and its adjoint can be found in \cite[Example 4.2]{HNS95} and we skip these details here. 

Note that this operator obeys the tangential cone condition as shown in \cite{}. Thus Assumption \ref{ass:tcc} is satisfied. To treat white noise, we choose $\V = H^a \left(\Omega\right)$ with $a > 1 = d/2$ to ensure $\norm{\xi_1}{\V'} < \infty$ a.s. . Furthermore, the Fr\'echet derivative $v=F'[u]h$ for $u\in L^{2}(\Omega), h\in L^{2}(\Omega)$ can - as shown in \cite[Example 4.2]{HNS95} - be expressed as the solution to
\begin{align}\label{Exeq2}
-\Delta v + uv &=-h\op(u), \quad \text{in} ~\Omega, \\
v&=0, \qquad \qquad\text{on}~ \partial \Omega. \nonumber
\end{align}
Note that the weak form of \eqref{Exeq2} has unique solution $v\in H^2(\Omega) \cup H^1_0(\Omega)$. This representation now allows us to verify Assumption \ref{ass:smoothing} whenever $a < 2$: By means of elliptic regularity theory, the operator $F'[u] : L^2 \left(\Omega\right) \to H^2 \left(\Omega\right)$ is bounded (in fact a homomorphism), and thus it follows from \cite[Rem. 2.6]{HW22} that Assumption \ref{ass:smoothing} is satisfied with $\theta = \frac{a}{2}$ and $C_\theta = \norm{F'[u]}{L^2 \to H^2}$. Finally, we can also verify Assumption \ref{ass:vsc} similar to \cite[Ex. 2.2]{HW22} by using the tangential cone condition. Precisely, if $\udag \in H^s \left(\Omega\right)$ for some $s > 0$, then \eqref{eq:vsc} holds true with $\varphi(\lambda) = C \lambda^{\frac{s}{s+2}}$ with some constant $C > 0$.

Our aim is to obtain the optimality conditions for the minimization procedure in Algorithms \ref{alg:cIRGNM}-\ref{alg:hIRGNM}. We note that the cost functionals in all these algorithms only vary in the measurements that they employ. Hence, here we focus only on the generic form of the minimization given in (\ref{eq:IRGNM}) and which, for the example under consideration, can be written as the minimizer of
\begin{align}
\nonumber
\mathcal{Q}(u,v) :=  &\frac12 \norm{\op(\hat{u}_{n})+v}{L^{2}(\Omega)}^2 - \ip{W}{\op(\hat{u}_{n})+v}_{L^{2}(\Omega)} + \frac{\alpha_n}{2} \norm{u - \hat{u}_0}{L^{2}(\Omega)}^2,
\end{align}
where $v$ satisfies the constraint
\begin{align}\label{Exeq3}
-\Delta v + \hat{u}_{n}v &=(\hat{u}_{n}-u)\op(u_{n}), \quad \text{in} ~\Omega, \\
v&=0, \qquad \qquad \qquad\quad \text{on}~ \partial \Omega. \nonumber
 \end{align}
Let us define the Lagrangian $\mathcal{L}:V\times V\times L^{2}(\Omega)\to \mathbb{R}$:
\begin{align}
\mathcal{L}(v,\lambda,u):= \mathcal{Q}(u,v) +\ip{-\Delta v + \hat{u}_{n}v -(\hat{u}_{n}-u)\op(\hat{u}_{n})}{\lambda}_{L^{2}(\Omega)},
\end{align}
which we now employ to solve the unconstrained optimization problem. To this end, we derive expression for the optimality conditions:
\begin{align}
D_{v}\mathcal{L}(v,\lambda,u)\tilde{v}=0,\label{Exeq4}\\
D_{\lambda}\mathcal{L}(v,\lambda,u)\tilde{\lambda}=0,\label{Exeq5}\\
D_{u}\mathcal{L}(v,\lambda,u)h=0,\label{Exeq6}
\end{align}
for all $(\tilde{v},\tilde{\lambda},h)\in V\times V\times L^{2}(\Omega)$. It follows trivially that the condition (\ref{Exeq4}) yields directly the constraint (\ref{Exeq3}). Furthermore, note that
\begin{align}\nonumber
D_{v}\mathcal{L}(v,\lambda,u)\tilde{v}= &\ip{\op(\hat{u}_{n})+v-W}{\tilde{v}}_{L^{2}(\Omega)}   +\ip{-\Delta \tilde{v} + \hat{u}_{n}\tilde{v} }{\lambda}_{L^{2}(\Omega)},
\end{align}
which, after integrating by parts and applying boundary conditions yields
\begin{align}
D_{v}\mathcal{L}(v,\lambda,u)\tilde{v}= &\ip{\op(\hat{u}_{n})+v-Y_{n}-\Delta \lambda + \hat{u}_{n}\lambda }{\tilde{v}}_{L^{2}(\Omega)}.
\end{align}
Hence, (\ref{Exeq5}) is equivalent to the following adjoint equation for $\lambda\in V$
\begin{align}\label{Exeq7}
-\Delta \lambda + \hat{u}_{n}\lambda &=W-\op(\hat{u}_{n})-v,
\end{align}
with homogeneous Dirichlet boundary conditions. Finally, it is easy to see that (\ref{Exeq6}) is equivalent to
\begin{align}\label{Exeq8}
u& =\hat{u}_0 -\frac{1}{\alpha_{n}}\lambda \op(\hat{u}_{n}).
%D_{u}\mathcal{L}(v,\lambda,u)h= & \ip{\alpha_{n}(u-u_0)+\lambda p }{h}
\end{align}
We use the previous equation in (\ref{Exeq3}) which we then combine with (\ref{Exeq7}) to obtain the linear system on $(\lambda, v)$ given by
\begin{align}\label{Exeq9}
\left(\begin{array}{cc}
-\Delta+\hat{u}_{n}  & I\\
-\alpha_{n}^{-1} \big(\op(\hat{u}_{n})\big)^2  & -\Delta +\hat{u}_{n}\end{array}\right)\left(\begin{array}{c}\lambda  \\
v \end{array}\right)&=\left(\begin{array}{c}W-\op(\hat{u}_{n}) \\
(\hat{u}_{n}-\hat{u}_0)\op(\hat{u}_{n}) \end{array}\right),
\end{align}
where $I$ denotes the identity in $L^{2}(\Omega)$. At a given iteration level $n$, we solve (\ref{Exeq9}) and use $\lambda$ in (\ref{Exeq8}) to compute the update $\hat{u}_{n+1}$. Replacing with $W$ with $Z_{N}$ and $Z_{n}$ gives the corresponding updates for Algorithms \ref{alg:cIRGNM}-\ref{alg:dIRGNM}, respectively.

\subsubsection{Numerical results}

\begin{figure}[h!]
\centering
\includegraphics[scale=0.37,trim=50 0 50 0]{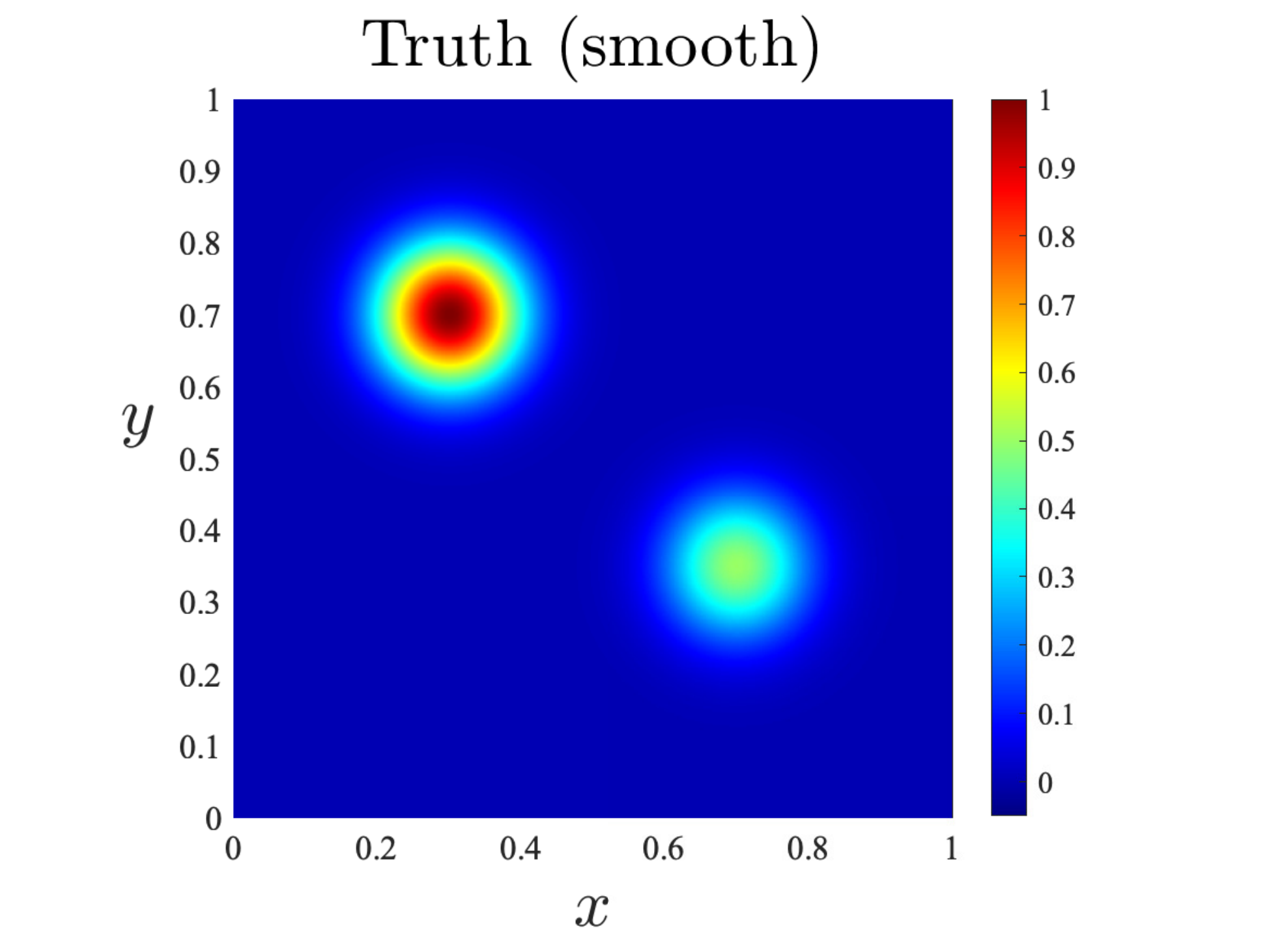}
\includegraphics[scale=0.37,trim=50 0 50 0]{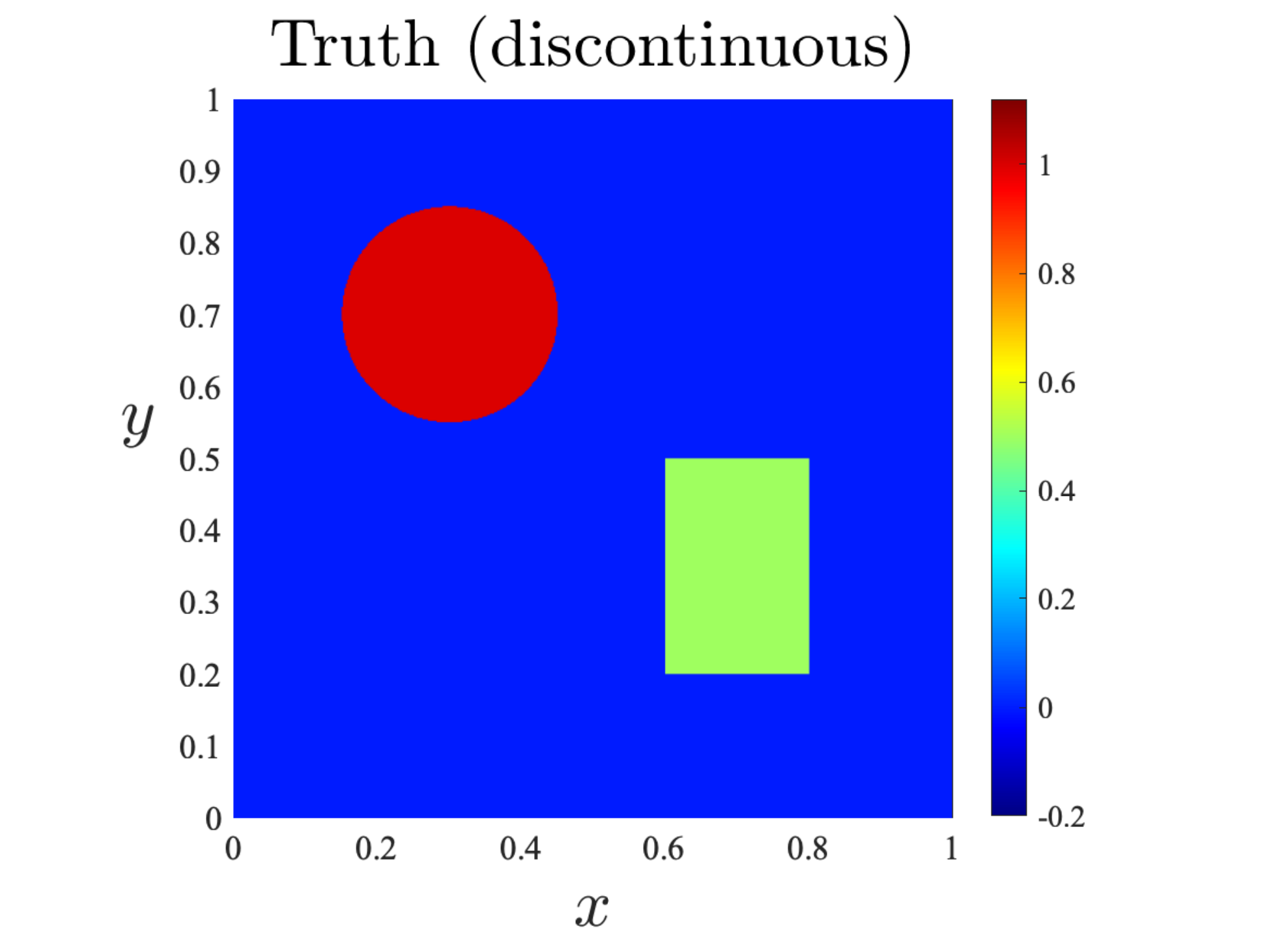}\\
 \caption{Example 1. True coefficient $u^{\dagger}(x,y)$ for the smooth (left) and discontinuous (right) cases.
 \ref{alg:cIRGNM}-\ref{alg:hIRGNM}.}
 %Bottom: Computational mesh employed for Algorithms 
    \label{FigExp1_1}
\end{figure}

We consider two experiments using a domain \newline $\Omega=[0,1]^{2}$. For the first set we consider a smooth truth defined by
$$u^{\dagger}(x,y)=\exp\Big[-100\Big( (x-0.3)^2+(y-0.7)^2\big))\Big]+\frac{1}{2}\exp\Big[-100\big( (x-0.7)^2+(y-0.35)^2\big)\Big],$$
while for the second we use
$$u^{\dagger}(x,y)=\left\{\begin{array}{cc}
  1,   & \text{if}~ (x-0.3)^2+(y-0.7)^2<0.15^2, \\
   0.5,  & \text{if}~(x,y)\in [0.6,0.8]\times[0.2,0.5],\\
   0, & \text{otherwise}.\end{array}\right.$$
In the top row of Figure \ref{FigExp1_1} we show the plots of these functions. For both cases we define $f(x,y)=(x+y)u^{\dagger}(x,y)$ and $g(x,y)=(x+y)\vert_{\partial \Omega}$, so that the noise-free data is given by $F(u^{\dagger})=p^{\dagger}(x,y)=(x+y)$. We 
specify a RHS of $f=1$ for the Darcy flow PDE.
%\todo[inline]{Discussion re. theory}

We implement Algorithms \ref{alg:cIRGNM}-\ref{alg:hIRGNM} in \texttt{MATLAB} and use pdetool toolbox to solve equation (\ref{Exeq1}) as well as the linear system (\ref{Exeq9}). We use ta mesh which consists of 7444 linear elements and 3837 nodes. Using the analytical solution, and thus avoiding inverse crimes, we evaluate the noise free observations on the nodes of the computational mesh, and produce the sequence of synthetic observations $Y_{n}$ (see e.g. (\ref{eq:model_dyn})), using a Gaussian random vector $\xi_{n}\in \mathbb{R}^{3837}$ with zero mean and standard deviation $\sigma=5\times 10^{-4}$.
For all algorithms we use $\hat{u}_{0}(x,y)=0$ (for all $(x,y)\in \Omega$) and $\alpha_{0}=10^{-3}$. For Algorithm \ref{alg:cIRGNM} and for the second part of Algorithm \ref{alg:dIRGNM} we use $C_{dec}=1.5$ in the definition of $\alpha_{n}$.

\begin{figure}[h!]
\centering
\includegraphics[scale=0.33,trim=0 0 0 0]{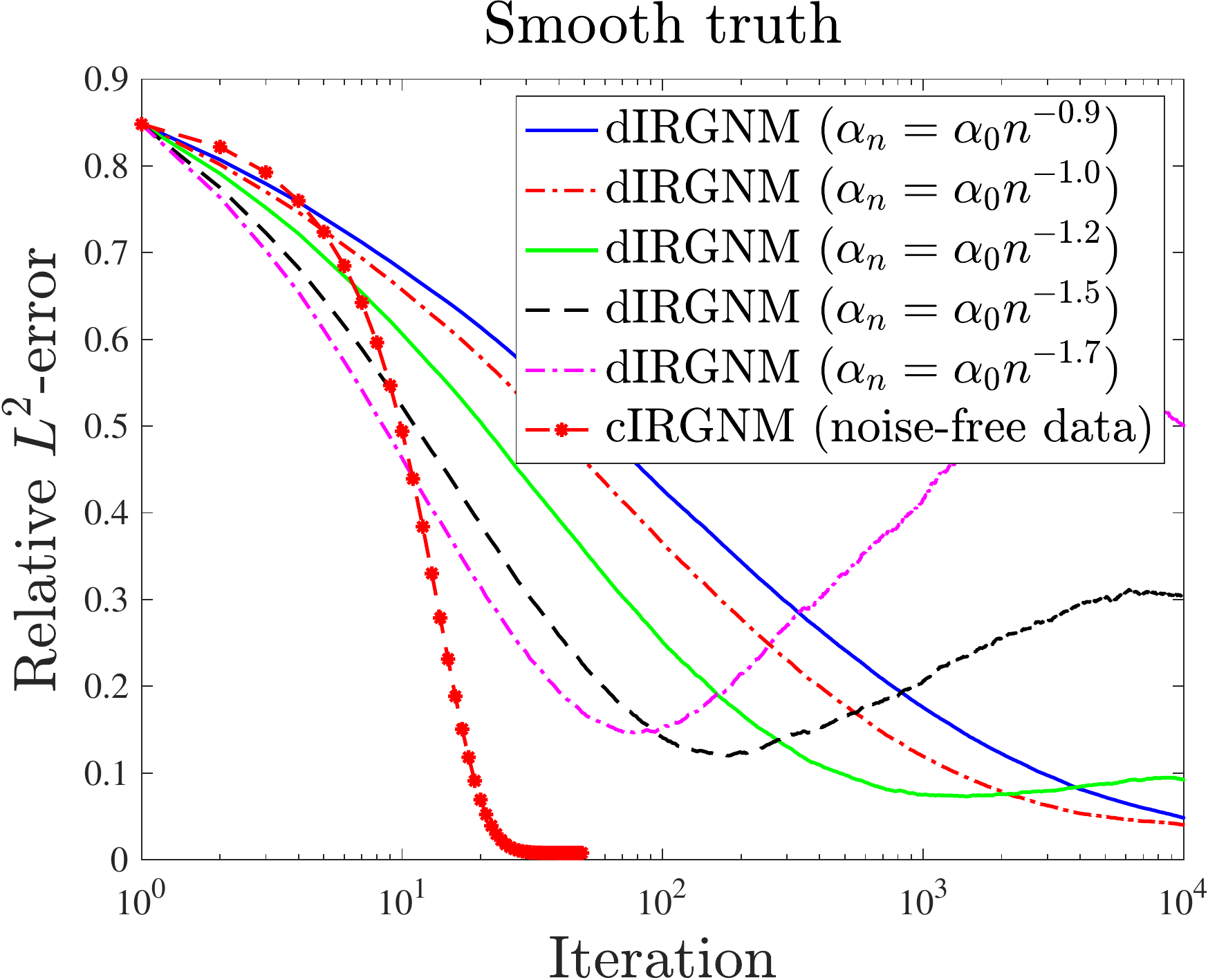}
\vspace{3mm}
\includegraphics[scale=0.33,trim=0 0 0 0]{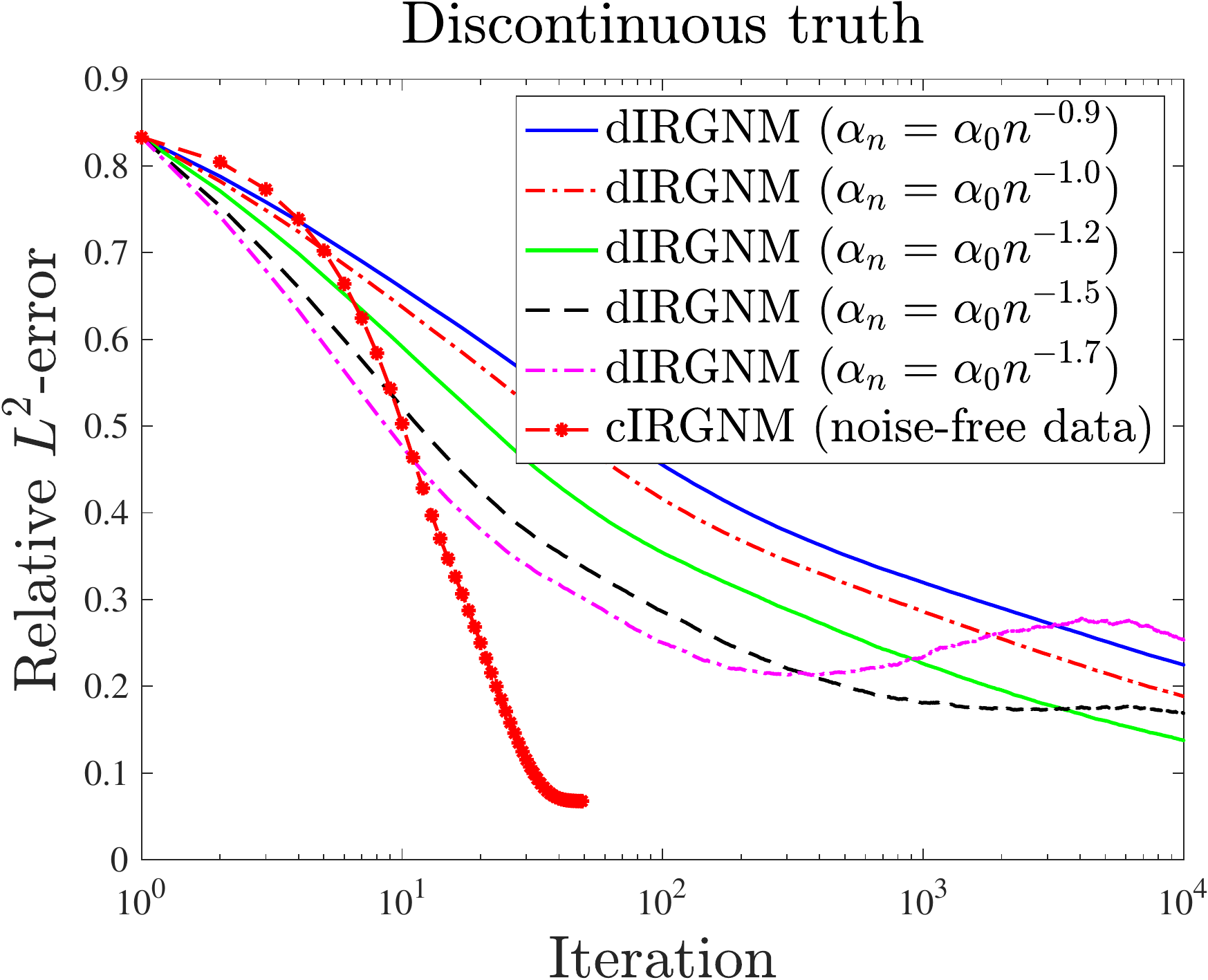}
 %\caption{Initial experiment for the iRGNM with the true underlying unknown, and the averaged noise after $n=25$ iterations.}
 \caption{Example 1. Relative $L^{2}$ errors obtained using the dIRGNM with various choices of $\beta$ for the continuous (left) and discontinuous (right) truth.}
    \label{FigExp1_4}
\end{figure}

\begin{figure}[h!]
\centering
\includegraphics[scale=0.37,trim=0 0 0 0]{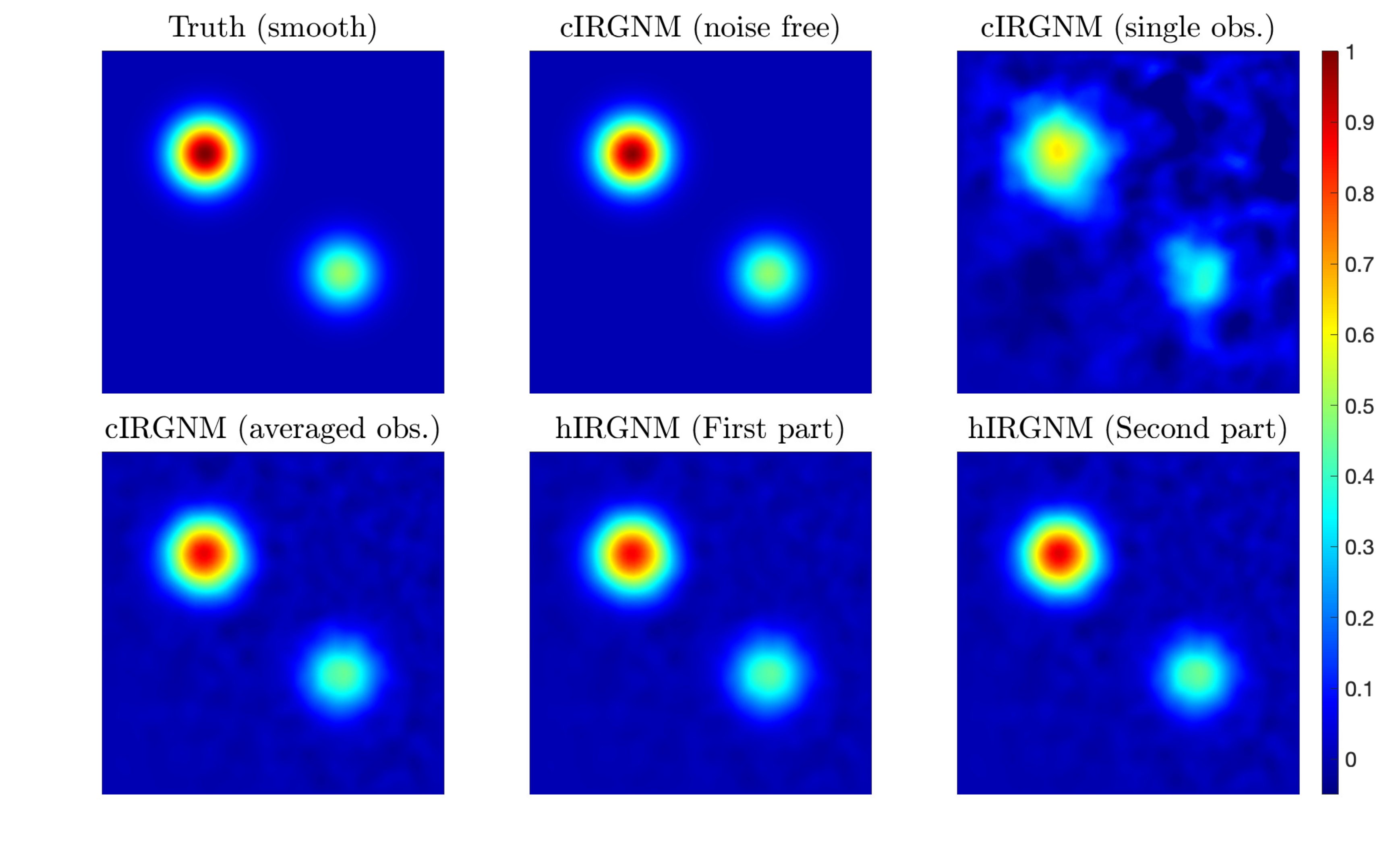}
 \caption{Example 1. (case with the smooth truth). Truth (top-left) and estimates of the unknown obtained with the cIGRNM with noise-free observations (top-middle), a single set of observations (top-right) and the average of $N=500$ observations (bottom-left). Bottom-middle and bottom-right panels show the estimates obtained from the first and second part of the hIRGNM using the same $N=500$ observations.}
    \label{FigExp1_2}
\end{figure}

\begin{figure}[h!]
\centering
\includegraphics[scale=0.37,trim=0 0 0 0]{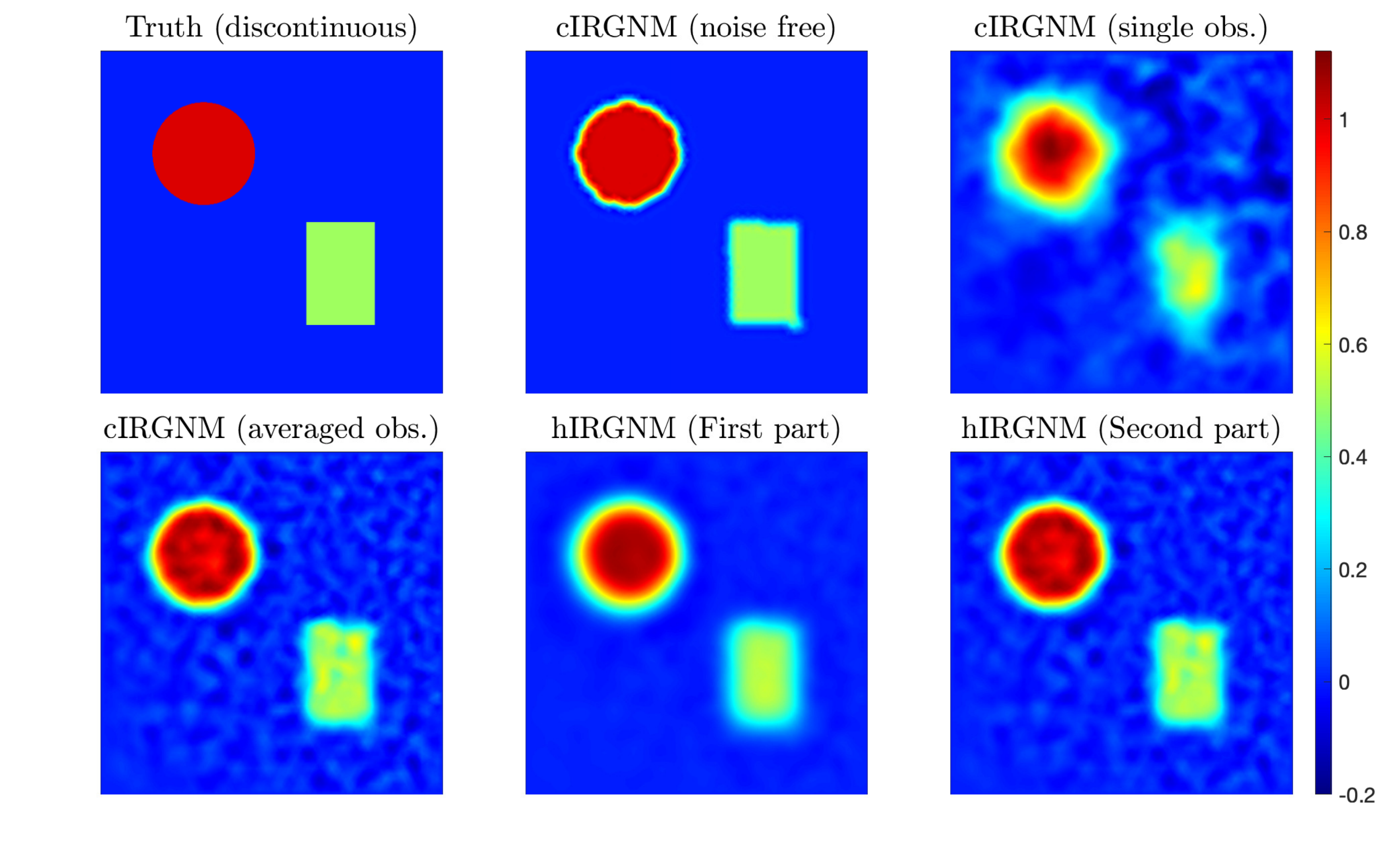}
  \caption{Example 1. (case with the discontinuous truth). Truth (top-left) and estimates of the unknown obtained with the cIGRNM with noise-free observations (top-middle), a single set of observations (top-right) and the average of $N=500$ observations (bottom-left). Bottom-middle and bottom-right panels show the estimates obtained from the first and second part of the hIRGNM using the same $N=500$ observations.}
     \label{FigExp1_3}
\end{figure}

To assess the convergence of the dIRGNM, we implement Algorithm \ref{alg:dIRGNM} with $N=10^4$ for various selections of $\beta$ in the definition of $\alpha_{n} := \alpha_0 n^{-\beta}$.
At each iteration we compute the relative error with respect to the truth defined by
$$E_{n}=\frac{\norm{\hat{u}_{n} - u^{\dagger}}{L_{2}} }{\norm{u^{\dagger}}{L_{2}} }.$$
As suggested in Item 2 of Remark \ref{rem:dIRGNM}, we shall choose $\beta\in (1/2,1-\theta/2)$ theoretically to obtain the convergence of dIRGNM. Such a remark is confirmed in Figure \ref{FigExp1_4} where relative error for various choices of $\beta$ are displayed with the smooth (resp. discontinuous) truth.
For validation purposes, in these plots we also display the relative error w.r.t the truth that we obtain from applying the cIRGNM with noise-free observations (i.e. we set $W=F(u^{\dagger})$). The estimates obtained with the noise-free case are highly accurate as we can also visually appreciate from the plots shown on the top-middle panels of Figures \ref{FigExp1_2}-\ref{FigExp1_3}. Though different choices of $\beta$ yield decaying relative error in the first hundred iterations, if $\beta>1$, we do obtain some amplified relative error when a sufficiently large number of observations are averaged.

\begin{figure}[h!]
\centering
\includegraphics[scale=0.35,trim=0 0 0 0]{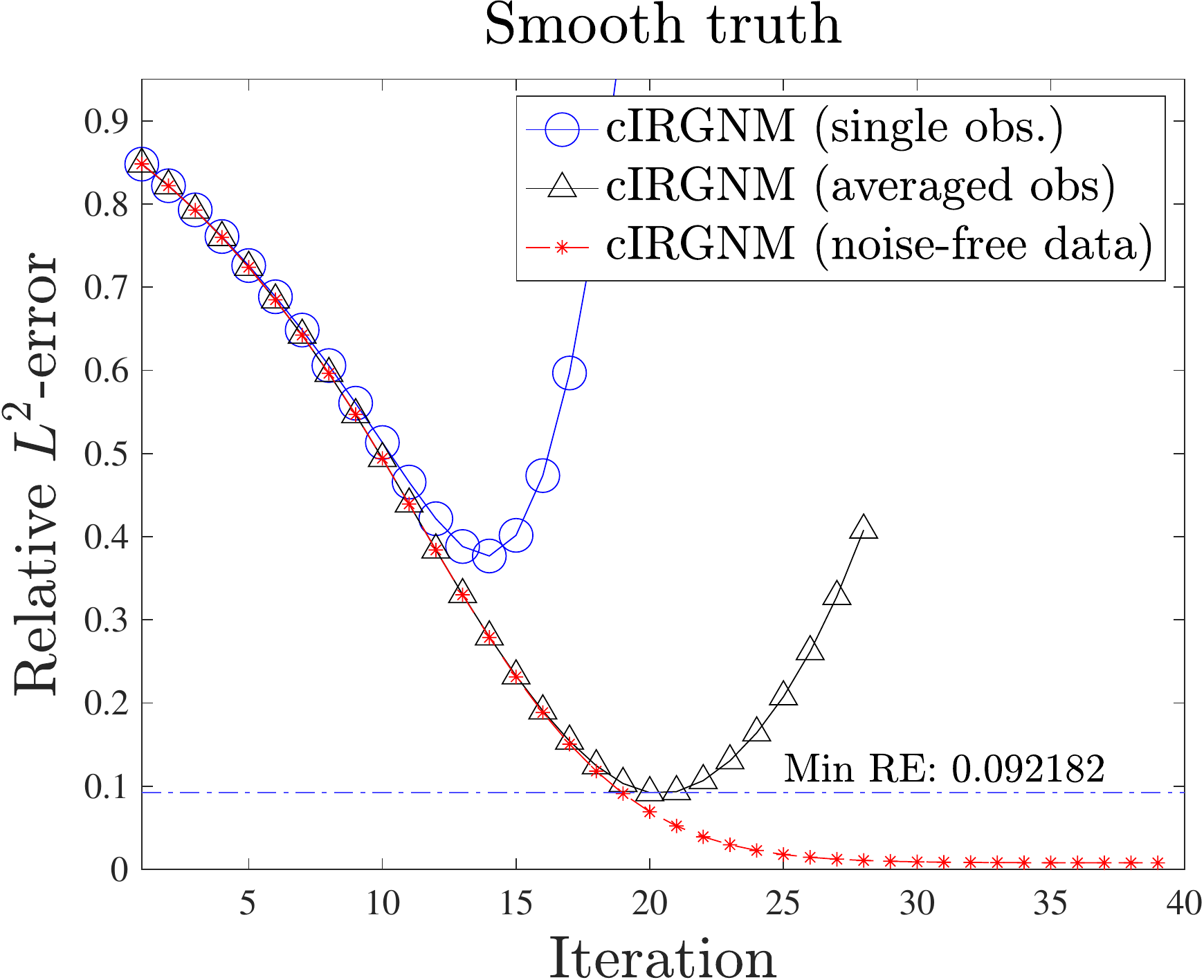}
\includegraphics[scale=0.35,trim=0 0 0 0]{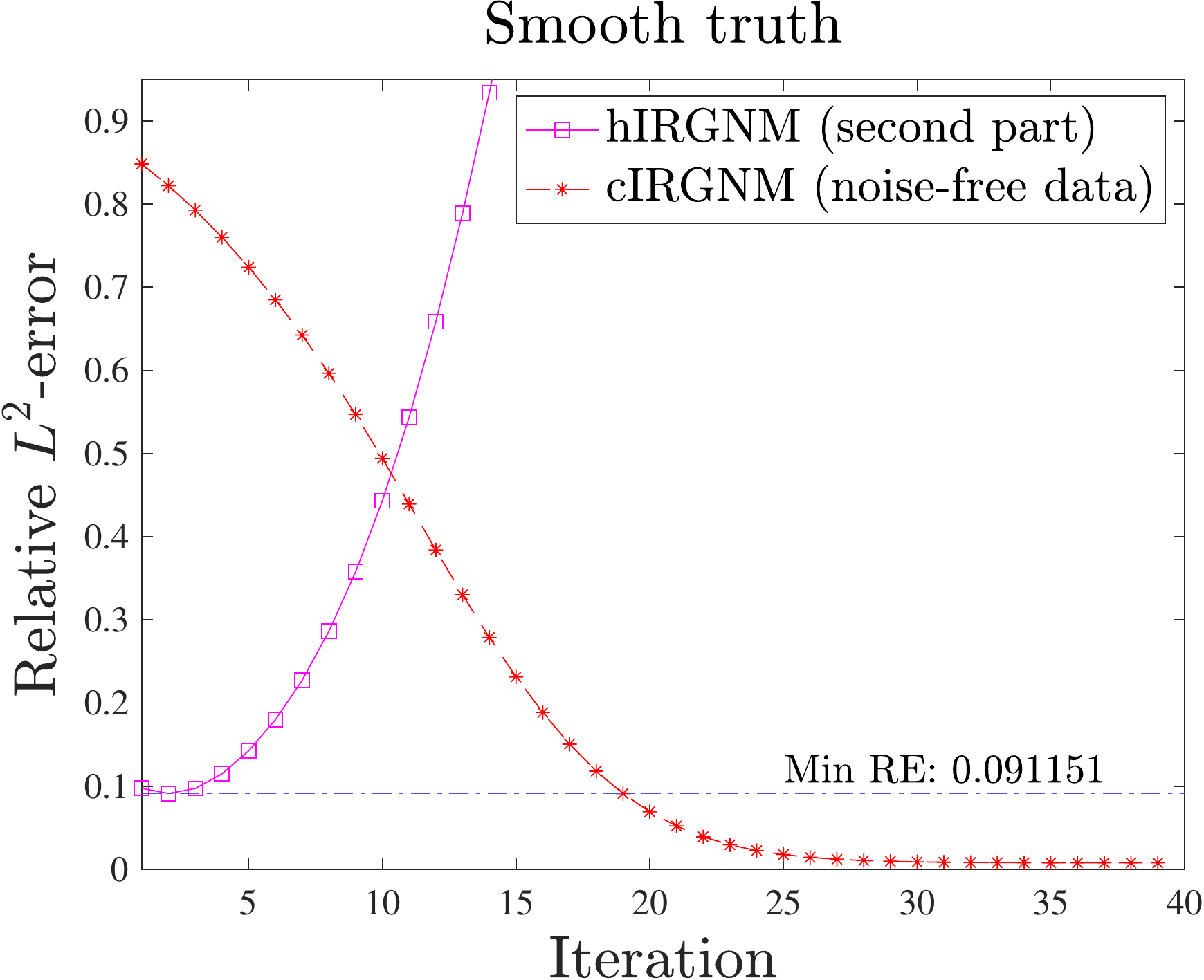}\\
\vspace{3mm}
\includegraphics[scale=0.35,trim=0 0 0 0]{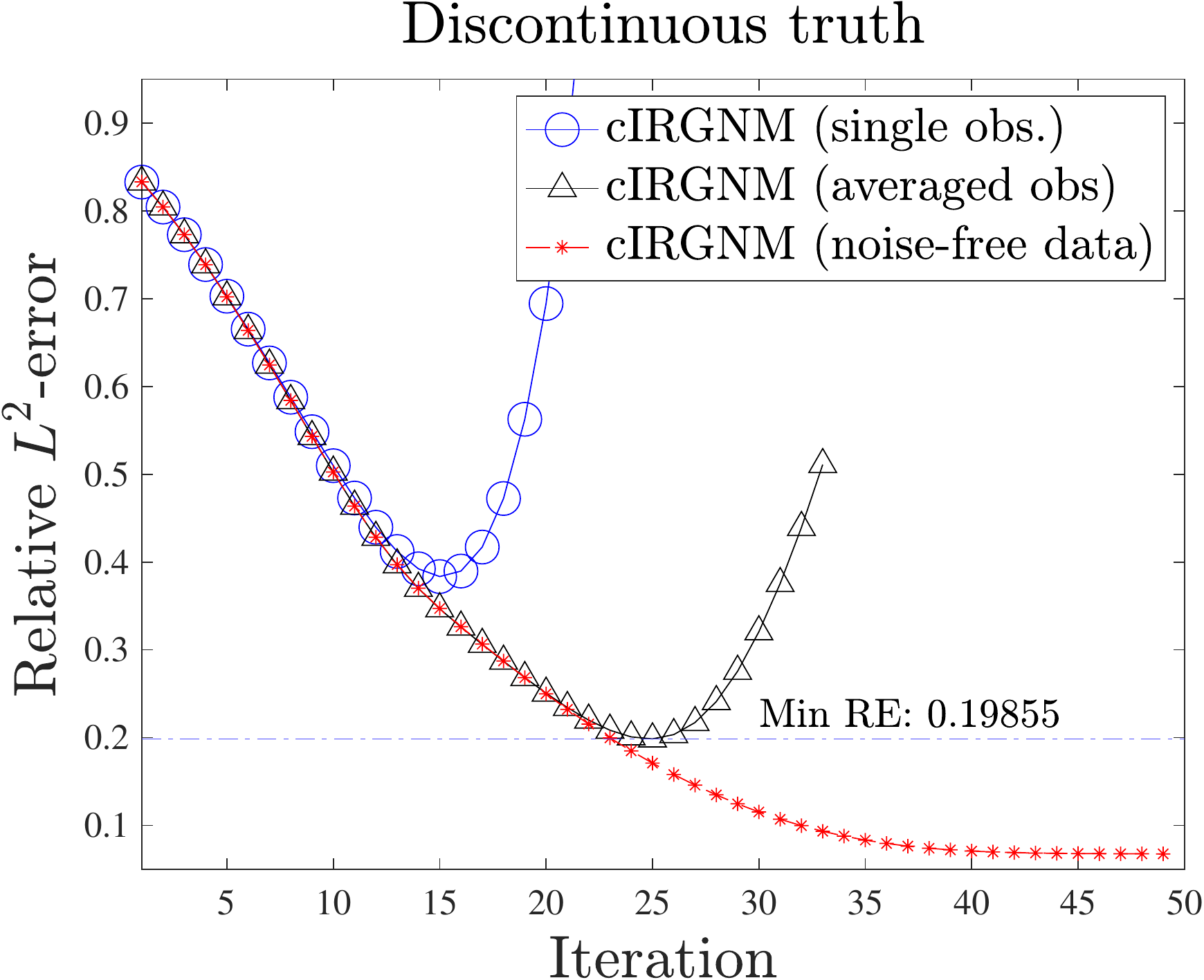}
\includegraphics[scale=0.35,trim=0 0 0 0]{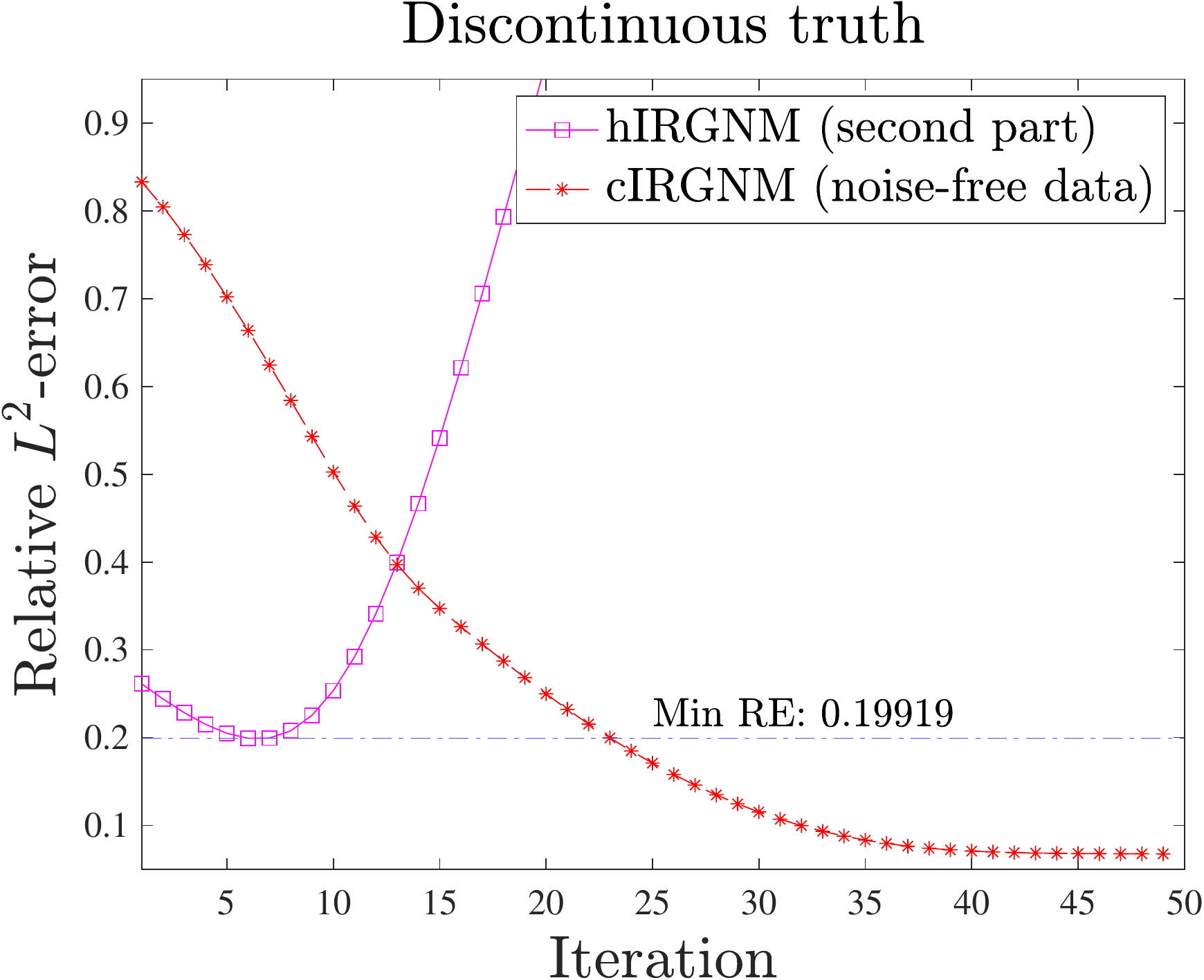}
 %\caption{Initial experiment for the iRGNM with the true underlying unknown, and the averaged noise after $n=25$ iterations.}
\caption{Example 1. Left: Relative $L^{2}$ errors for the case with the smooth (top) and discontinuous (bottom) truth obtained using the cIRGNM with three different observations: (i) noise-free, (ii) a single one and (iii) the averaged of $N=500$. Right: Relative $L^{2}$ errors obtained during the second part of the hIRGNM with the same $N=500$ observations. For comparisons the right panels also display the relative error obtained with the noise-free cIRGNM. The numerical values displayed on the left (resp. right) plots corresponds to the minimum relative error achieved via the cIRGNM with averaged measurements (resp. the second part of the hIRGNM).}
    \label{FigExp1_5}
\end{figure}

We now investigate the case where we have access to only finitely many observations and compare the performance of the (online+offline) hIGRNM  and the(offline) cIRGNM. More precisely we consider $N=500$ synthetic set of observations which we use for (i) the hIRGNM (Algorithm \ref{alg:hIRGNM}) with sequential averaged observations $W=\{Z_i\}_{i=1}^{500}$ and (ii) the cIRGNM (Algorithm \ref{alg:cIRGNM}) with $W=Z_{500}$. To further demonstrate the advantage of using multiple observations we also implement the cIRGNM with $W=Y_{1}$ which corresponds to the standard approach of aiming at reconstructing the unknown with a single data set. The relative error w.r.t the truth obtained with the cIRGNM are shown in the left panels of Figure \ref{FigExp1_5}. As we expect when noisy observations are employed, the error starts increasing after several iterations due to the ill-posedness of the inverse problem. Since the noise level of the averaged observation is smaller than that of a single observation, it comes as no surprise that the cIRGNM with $W=Z_{500}$ reaches a lower minimum value (displayed on the plots). The corresponding estimates of the unknown for which the relative error reaches the minimum are shown in the top-right and bottom-left panels of Figures \ref{FigExp1_2}-\ref{FigExp1_3}.

For the dIRGNM encoded in the first part of the hIRGNM we use $\beta=1.2$ informed by the previous experiments that suggest that we can safely go slightly above the value predicted by the theory in order to achieve faster convergence without the risk of compromising accuracy. The error from the iterations during the first part of the hIRGNM corresponds to the first $N=500$ iterations shown in Figure \ref{FigExp1_4} (for $\beta=1.2$). In the right panels of Figure \ref{FigExp1_5} we show the iterations achieved during the second part of the hIRGNM (Algorithm \ref{alg:hIRGNM}). The minimum relative error achieved by the hybrid method is shown on the corresponding plots. When we compare left and right panels, we can notice that the minimum relative error value attained by the hIRGNM is very similar to the one obtained using the cIRGNM with $W=Z_{500}$. However, we notice the second part of the hybrid method reaches the minimum after a only a few iterations. In fact, for case with a smooth truth, the second part of the hIRGNM takes only one iteration to reach the minimum value. In the bottom-middle and bottom-left panels of Figures \ref{FigExp1_2}-\ref{FigExp1_3} we show the estimated from the first and the second part (when minimum is reached) of the hIRGNM, respectively. We can visually appreciate that the (dynamic) first part of the hybrid algorithm exhibits quite a good level of accuracy.

The selected realistic value of the noise standard deviation ($\sigma=5\times 10^{-4}$) enabled us to show the advantage of inverting the average of multiple observations compared to the standard practice of inverting a single set of observations. Nonetheless, it is worth mentioning that this value of $\sigma$ is small enough so that we did not observe substantial differences when using different realizations of the noisy observations that we produced. In effect, we conducted multiple experiments (not shown) with different random selections of the $N=500$ observations and the results showed consistency across the experiments even in the case where only a single data set ($W=Y_1$) was inverted via the cIRGNM.

\subsection{Example 2.}
For our second numerical example we consider the estimation of the log-permeability of a porous medium. In particular, for the forward model, given a source function $f \in L^{\infty}(\Omega)$,
where ${\Omega} \subset \mathbb{R}$ is a Lipschitz domain, and permeability $\kappa=\exp(u) \in L^{\infty}(\Omega)$, we are interested
in solving the following Darcy flow model
\begin{align}
\label{eq:darcy}
-\nabla \cdot (\exp(u) \nabla p)&=f, \quad \in {\Omega}, \\
p&=0, \quad \in \partial {\Omega}, \nonumber
\end{align}
for the pressure $p \in H^1_0(D)$. The inverse problem associated with \eqref{eq:darcy}
is the recovery of the log-permeability $u$ given $K$ point-wise measurements of the pressure evaluated at interior testing points $\{x_{i}\}_{i=1}^{K}\subset \Omega$. In this situation it is unclear whether the tangential cone condition from Assumption \ref{ass:tcc} holds true. As our derivation of the variational source condition in the previous example makes use of the tangential cone condition, this also remains unclear here. However, Assumption \ref{ass:smoothing} could be verified similarly using the explicit form of $F'[u]$.

In this example, to treat more general setting, we consider the following weighted $L^2$ space
\begin{align}
\label{eq:darcy00}
\mathcal{H}\equiv \{u\in L^{2}(\Omega)\big\vert\quad \vert\vert  \mathcal{C}^{-1/2} u \vert \vert_{L{^2}(\Omega) } \leq \infty\},
\end{align}
where $\mathcal{C}$ is a covariance operator induced by a correlation function as follows
\begin{align}\label{eq:eit100}
\mathcal{C}[u](x)=\int_{\Omega}u(x')c(x,y')dxdx'.
\end{align}
We choose a Mat\'{e}rn correlation function given by
\begin{align}\label{matern}
c(x,x'):= c_{0}\frac{2^{1-\nu}}{\Gamma(\nu)}K_{\nu}\bigg(\frac{|x-x'|}{\ell}\bigg)\bigg(\frac{|x-x'|}{\ell}\bigg)^{\nu}.
\end{align}
where $c_{0}\in \mathbb{R}^+$ is a scaling factor, $\nu \in \mathbb{R}^+$ is a smoothness parameter, $\ell \in \mathbb{R}^+$ denotes the length-scale,
$\Gamma(\cdot)$ is the Gamma function and $K_{\nu}(\cdot)$ is the modified Bessel function of the second kind.

The forward map $F:\mathcal{H}\to \mathbb{R}^K$ is defined by $F(u)=(p(x_{1}),\dots,p(x_{K}))$ where $p$ is the solution to (\ref{eq:darcy}) evaluated..

In order to compute the minimizers in Algorithms \ref{alg:cIRGNM}-\ref{alg:hIRGNM}, we now modify our update formula (\ref{eq:gn_up}) in Section \ref{se1} based on the modified weighting of $\mathcal{C}$, which for the cIRGNM is given as
\begin{align}
\label{eq:darcy01}
\hat{u}_{n+1}-\hat{u}_{n}=(\op'[\hat{u}_{n}]^*\op'[u_{n}]+\alpha_{n}\mathcal{C}^{-1})^{-1}\Big(\op'[\hat{u}_{n}]^*(W-\op(\hat{u}_{n}))+\alpha_{n}\mathcal{C}^{-1}(\hat{u}_0-\hat{u}_{n})\Big),
\end{align}
and with suitable modifications for implementation of the hIRGNM and dIRGNM. For computational efficiency we can then use
Woodbury lemma for \eqref{eq:darcy01} yielding
\begin{align}\label{eq:darcy02}
\hat{u}_{n+1}=u_{0}+\mathcal{C}\op'(u_{n})^*(\op'(u_{n})\mathcal{C}\op'(u_{n})^*+\alpha_{n}I)^{-1}\Big(W-\op(u_{n})-\op'(u_{n})(u_0-u_{n})\Big).
\end{align}

We use \texttt{MATLAB} for the numerical implementation of Algorithms \ref{alg:cIRGNM}- \ref{alg:hIRGNM} and use bespoke solver based on a second-order centred finite difference method to numerically solve \eqref{eq:darcy}. The same scheme is used for the implementation of the Fr\'echet derivatives and the discrete adjoint equation which are derived as discussed in \cite{CRV99,ACD09,GDG14}.

\subsubsection{Numerical results}
\begin{figure}[h!]
\centering
\includegraphics[scale=0.37]{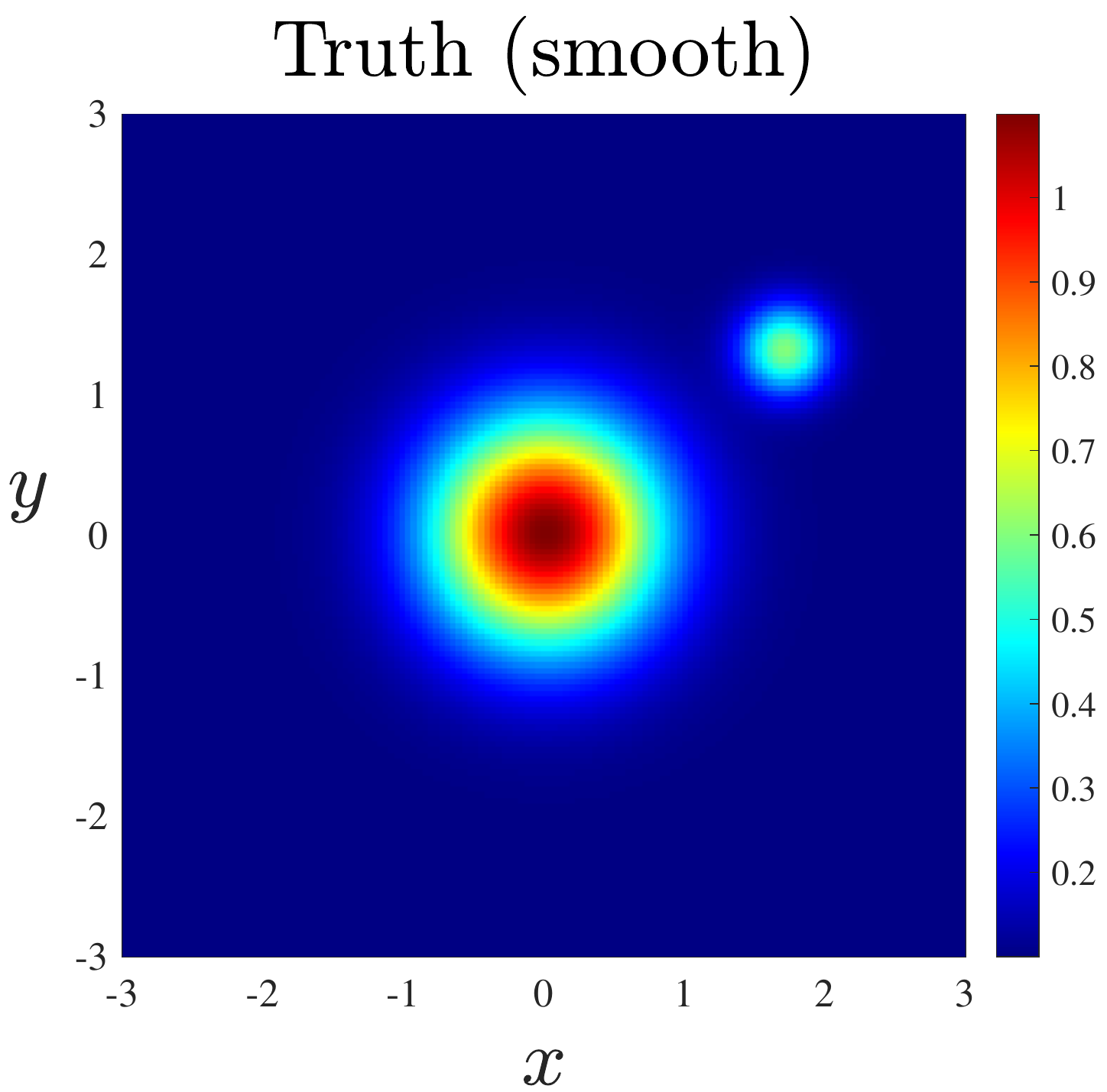}
\includegraphics[scale=0.37]{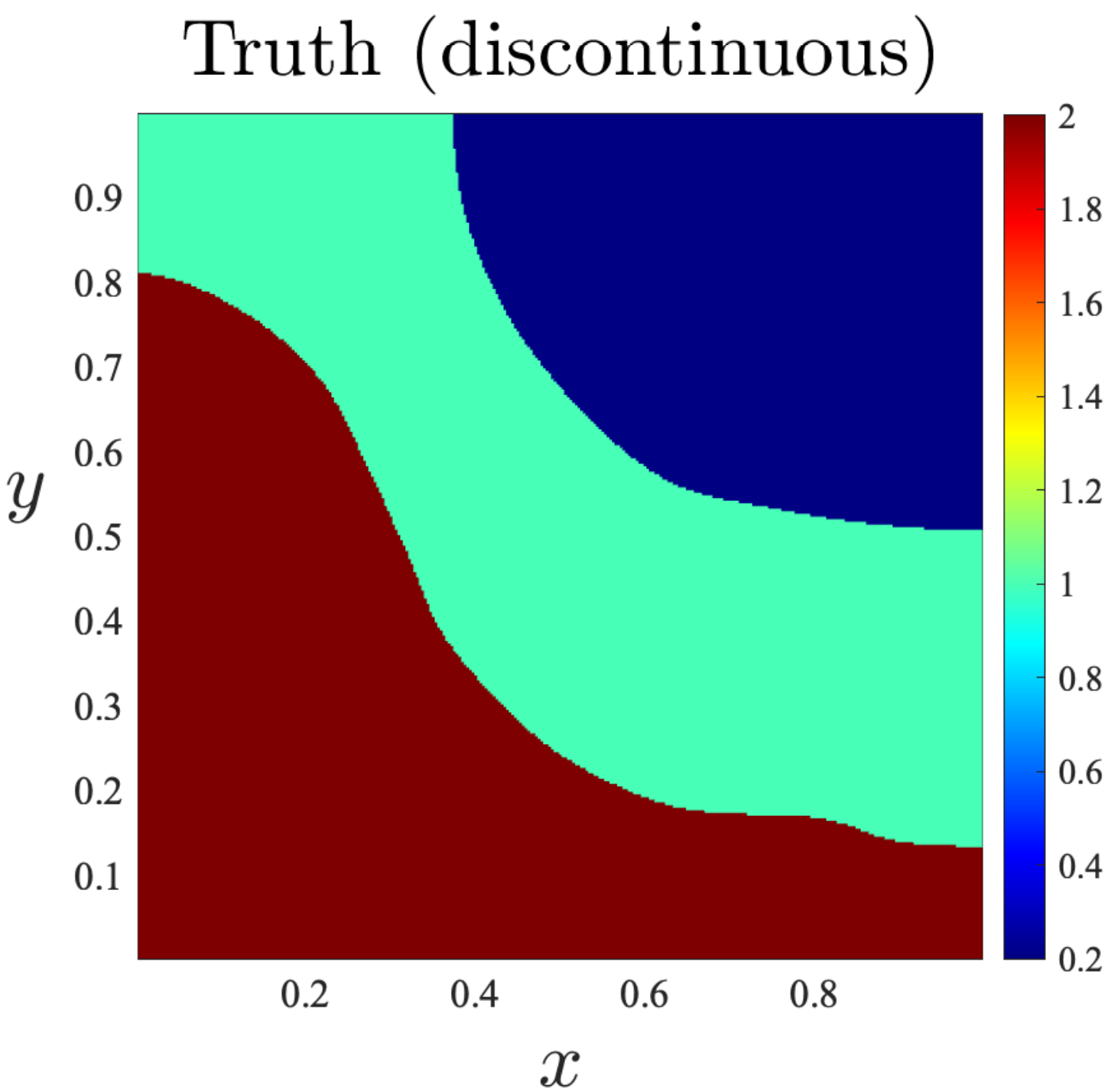}
 \caption{Example 2. True coefficient $u^{\dagger}(x,y)$ for the smooth (left) and discontinuous (right) case.}
\label{FigExp2_1}
\end{figure}

Here we consider two different domains for each numerical experiment. For the first experiment, we consider a domain of $\Omega=[-3,3]^2$
with a truth defined as
$$
u^{\dagger}(x,y)=\exp\Big[-100\Big( (x-0.3)^2+(y-0.7)^2\big))\Big]+\frac{1}{2}\exp\Big[-100\big( (x-0.7)^2+(y-0.35)^2\big)\Big],
$$
For the second experiment, we have a discontinuous truth which is defined on the domain $\Omega = [0,1]^2$, where the truth is taken to be a discontinuous-function
with some random features, which resembles a channel model \cite{ILS14}. Plots of the true permeabilities are presented in Figure \ref{FigExp2_1}. For both set of experiments we define a regular grid of $K=14\times 14$ testing points within $\Omega$.

To avoid inverse crimes we employ a mesh with $300 \times 300$ elements while a coarser mesh ($250 \times 250$) is used for the computations in Algorithms \ref{alg:cIRGNM}-\ref{alg:hIRGNM}.
As before, noisy observations $Y_{n}$ are obtained by adding Gaussian noise to the noise-free measurements, with standard deviation of $\sigma=2\times 10^{-3}$.
Furthermore, we use  $\alpha_{0}=10^{-3}$ and $C_{dec}=1.5$. In addition, we use $\hat{u}_0(x,y)=1$ (for all $(x,y)\in \Omega$) for the smooth truth case, while for the discontinuous case $\hat{u}_0$ is random sample from a Gaussian with covariance as defined in (\ref{eq:eit100}). For (\ref{matern}) we use parameters $c_{0}=1$, $\nu=3$ and $\ell=0.08$ 

\begin{figure}[h!]
\centering
\includegraphics[scale=0.31,trim=0 0 0 0]{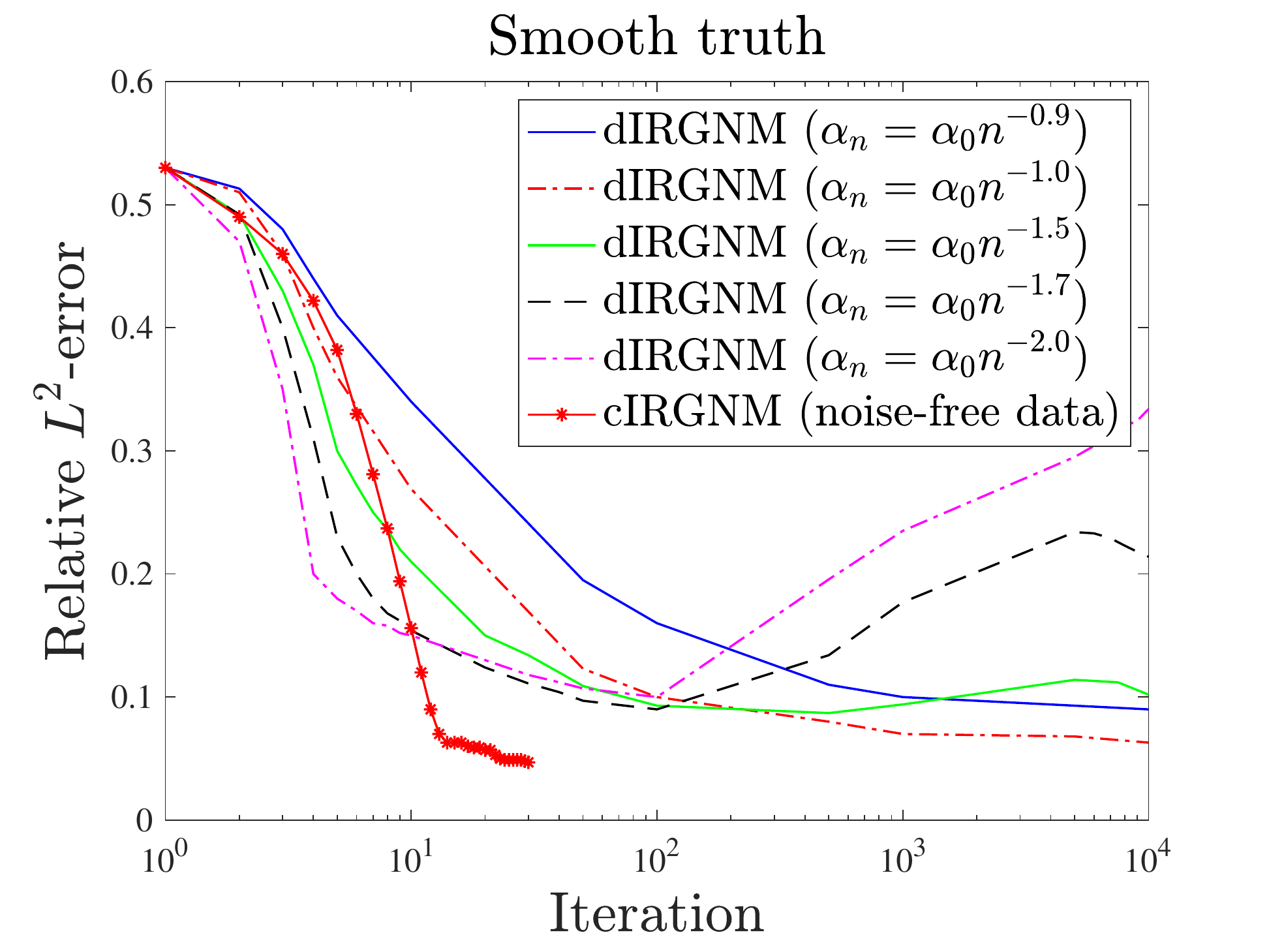}
\vspace{3mm}
\includegraphics[scale=0.31,trim=0 0 0 0]{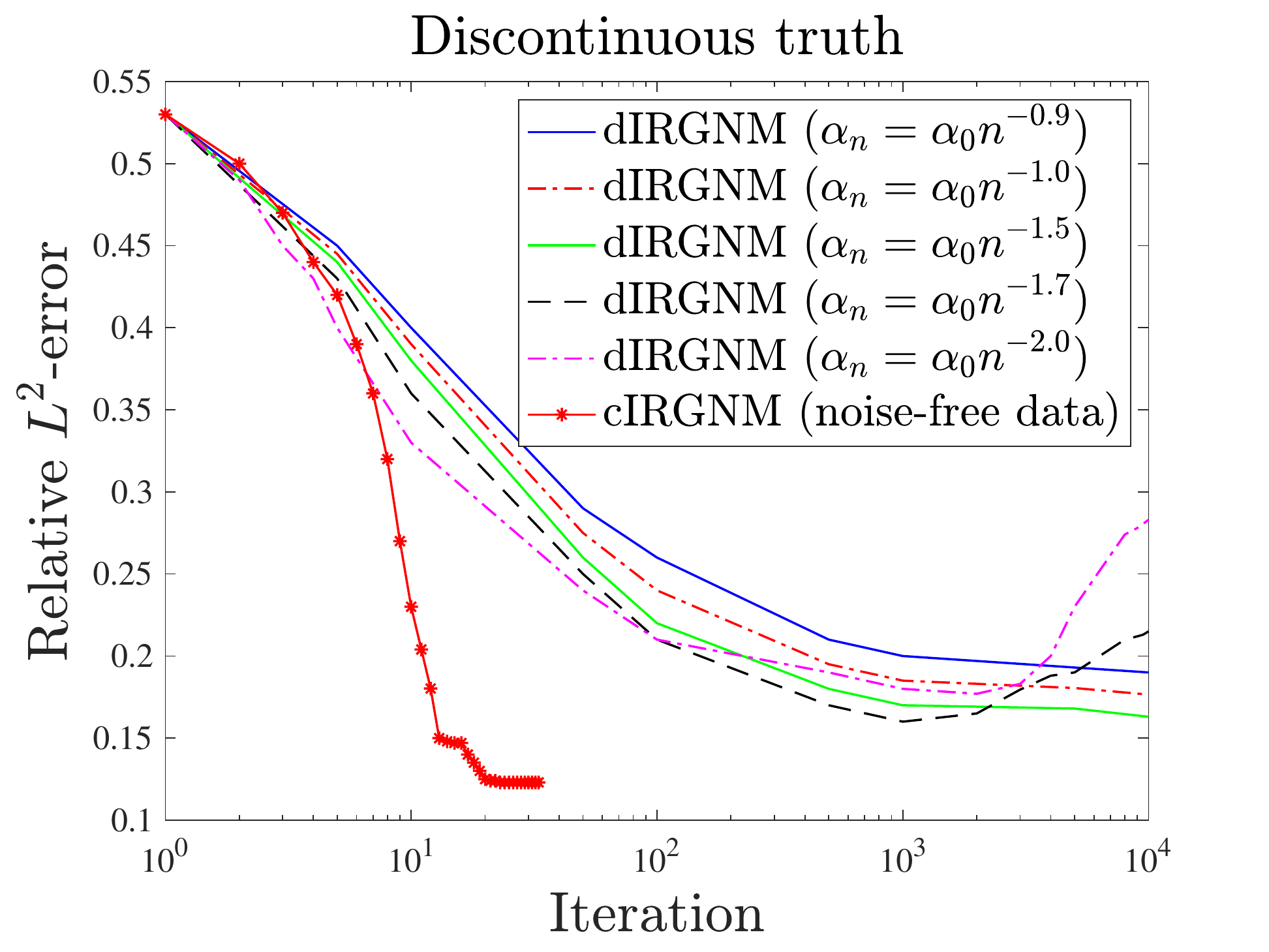}
 %\caption{Initial experiment for the iRGNM with the true underlying unknown, and the averaged noise after $n=25$ iterations.}
\caption{Example 2. Relative $L^{2}$ errors obtained using the dIRGNM with various choices of $\beta$ for the continuous (left) and discontinuous (right) truth.}
    \label{FigExp2_4}
\end{figure}

\begin{figure}[h!]
\centering
\includegraphics[scale=0.31,trim=0 0 0 0]{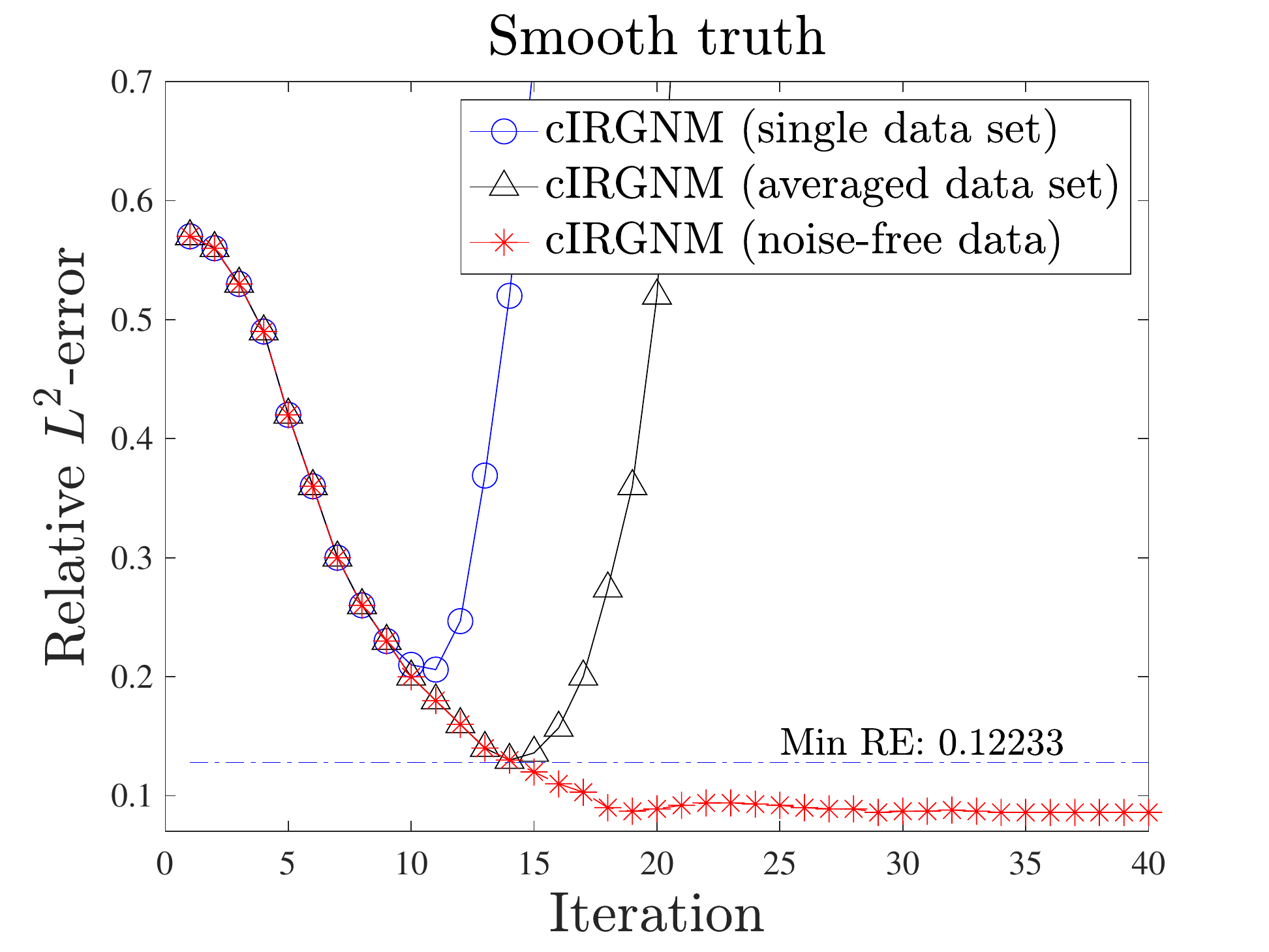}
\includegraphics[scale=0.31,trim=0 0 0 0]{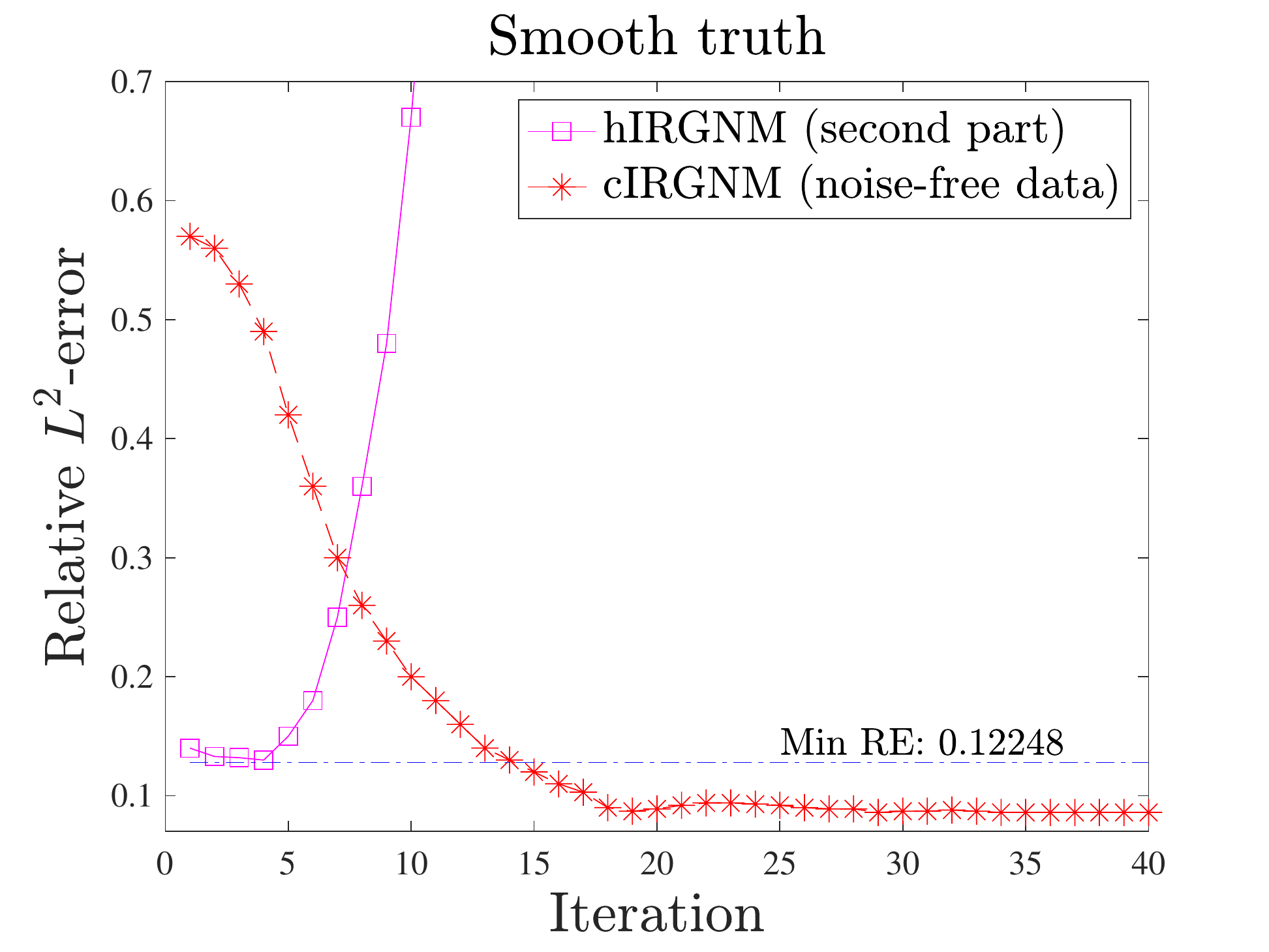}\\
\vspace{3mm}
\includegraphics[scale=0.31,trim=0 0 0 0]{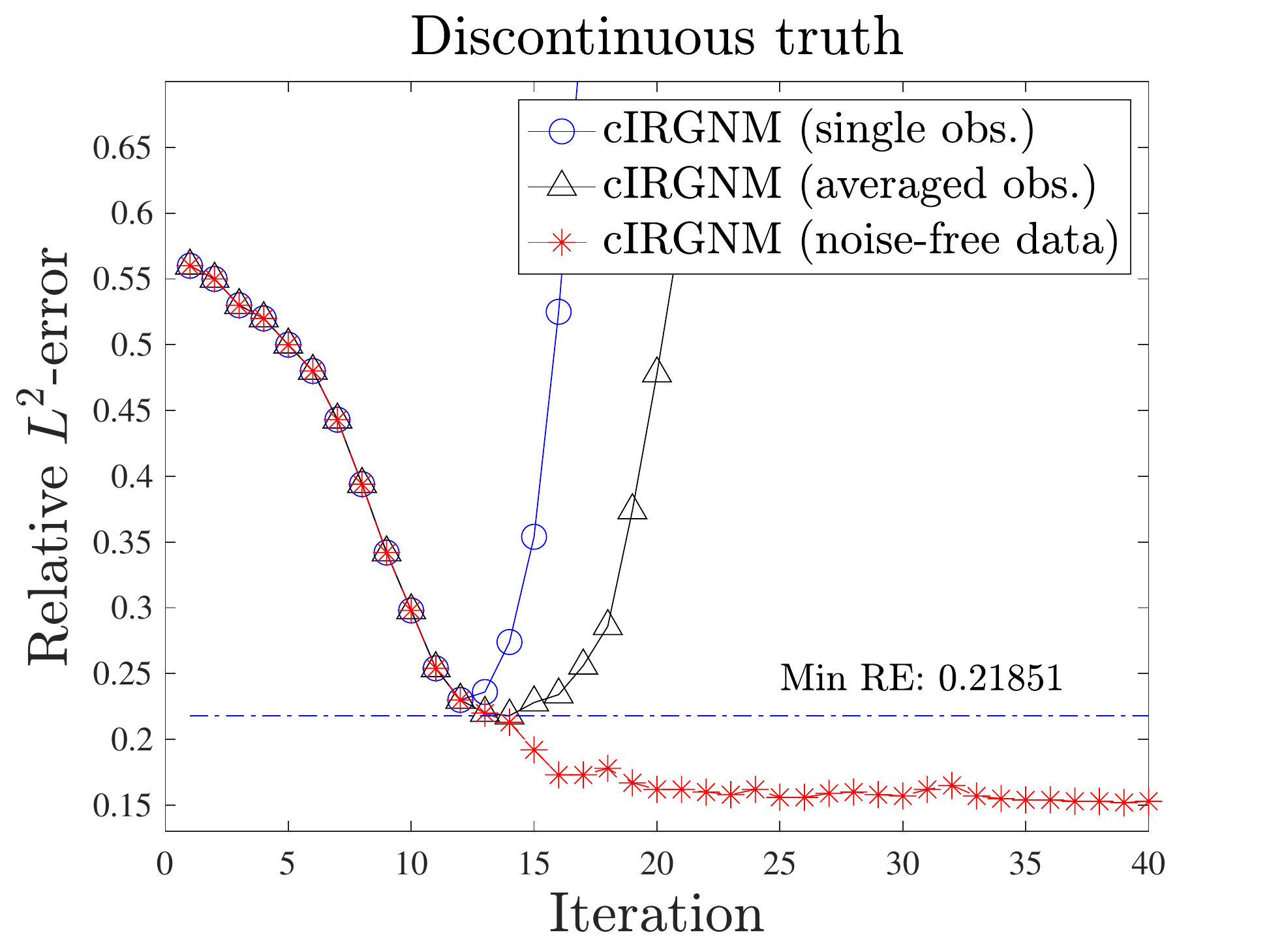}
\includegraphics[scale=0.31,trim=0 0 0 0]{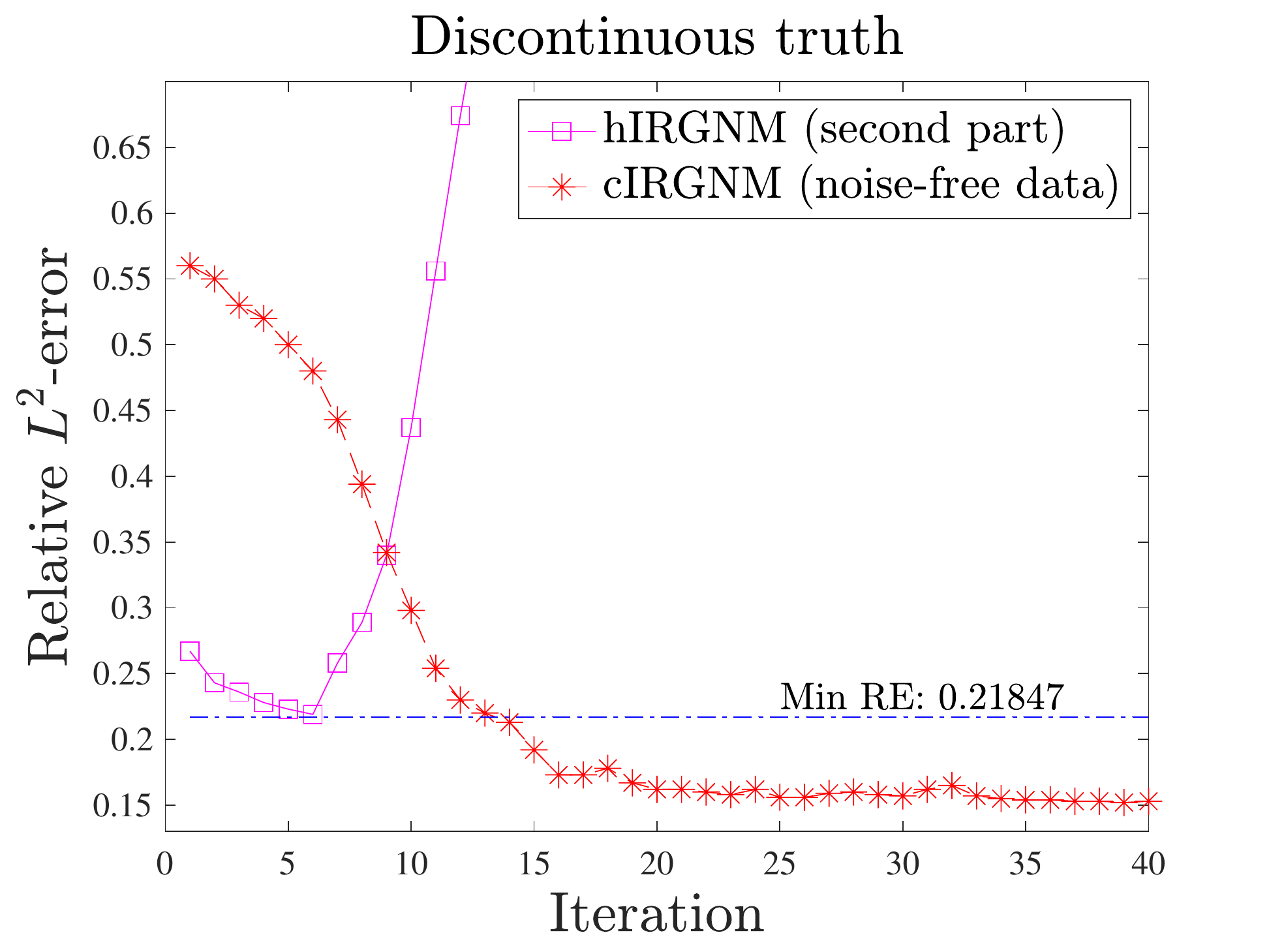}
\caption{Example 2. Left: Relative $L^{2}$ errors for the case with the smooth (top) and discontinuous (bottom) truth obtained using the cIRGNM with three different observations: (i) noise-free, (ii) a single one and (iii) the averaged of $N=500$. Right: Relative $L^{2}$ errors obtained during the second part of the hIRGNM with the same $N=500$ observations. For comparisons the right panels also display the relative error obtained with the noise-free cIRGNM. The numerical values displayed on the left (resp. right) plots corresponds to the minimum relative error achieved via the cIRGNM with averaged measurements (resp. the second part of the hIRGNM).}
    \label{FigExp2_5}
\end{figure}

\begin{figure}[h!]
\centering
\includegraphics[width=1.0\textwidth]{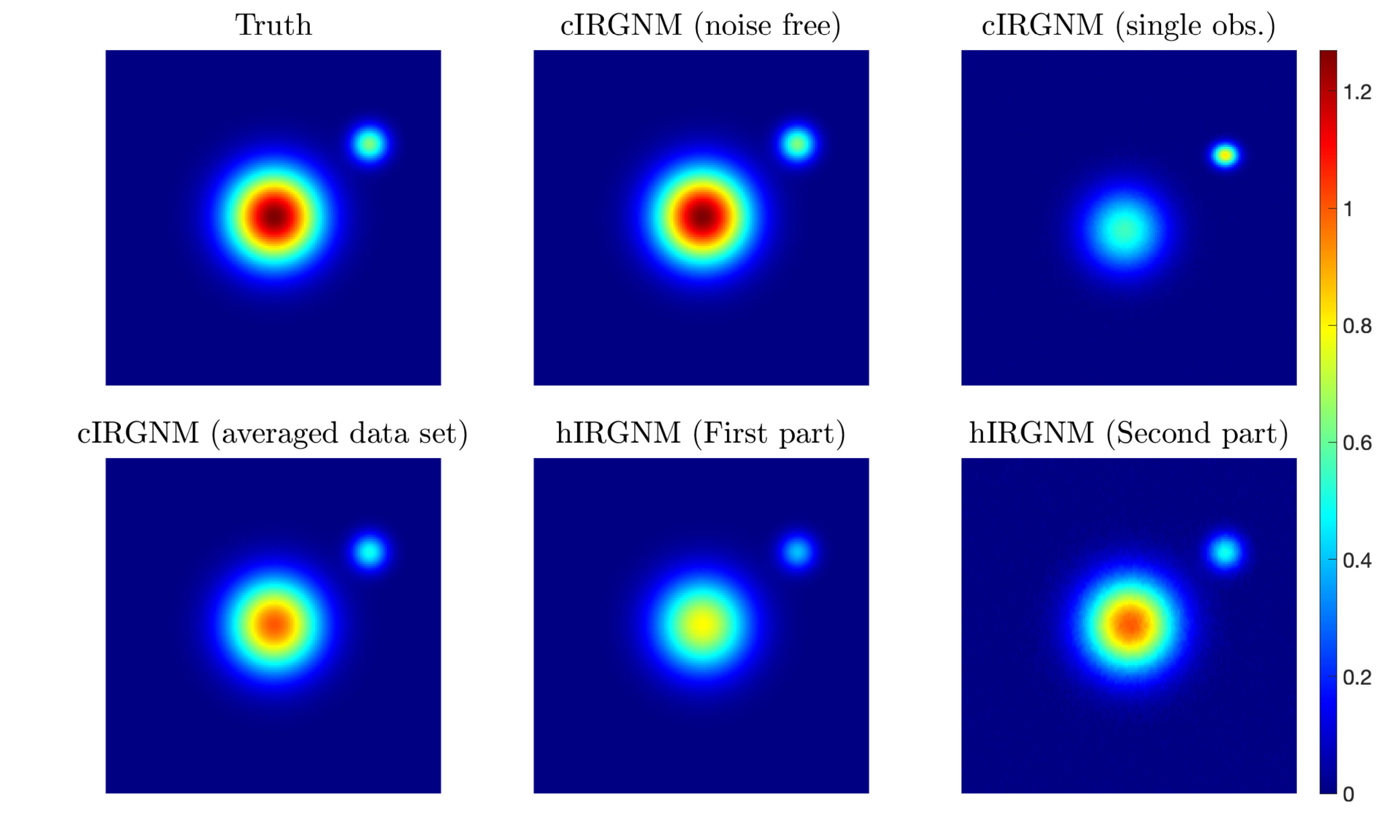}
 \caption{Example 2 (case with the smooth truth). Truth (top-left) and estimates of the unknown obtained with the cIGRNM with noise-free observations (top-middle), a single set of observations (top-right) and the average of $N=500$ observations (bottom-left). Bottom-middle and bottom-right panels show the estimates obtained from the first and second part of the hIRGNM using the same $N=500$ observations.}
 \label{FigExp2_2}
\end{figure}

\begin{figure}[h!]
\centering
\includegraphics[width=1.0\textwidth]{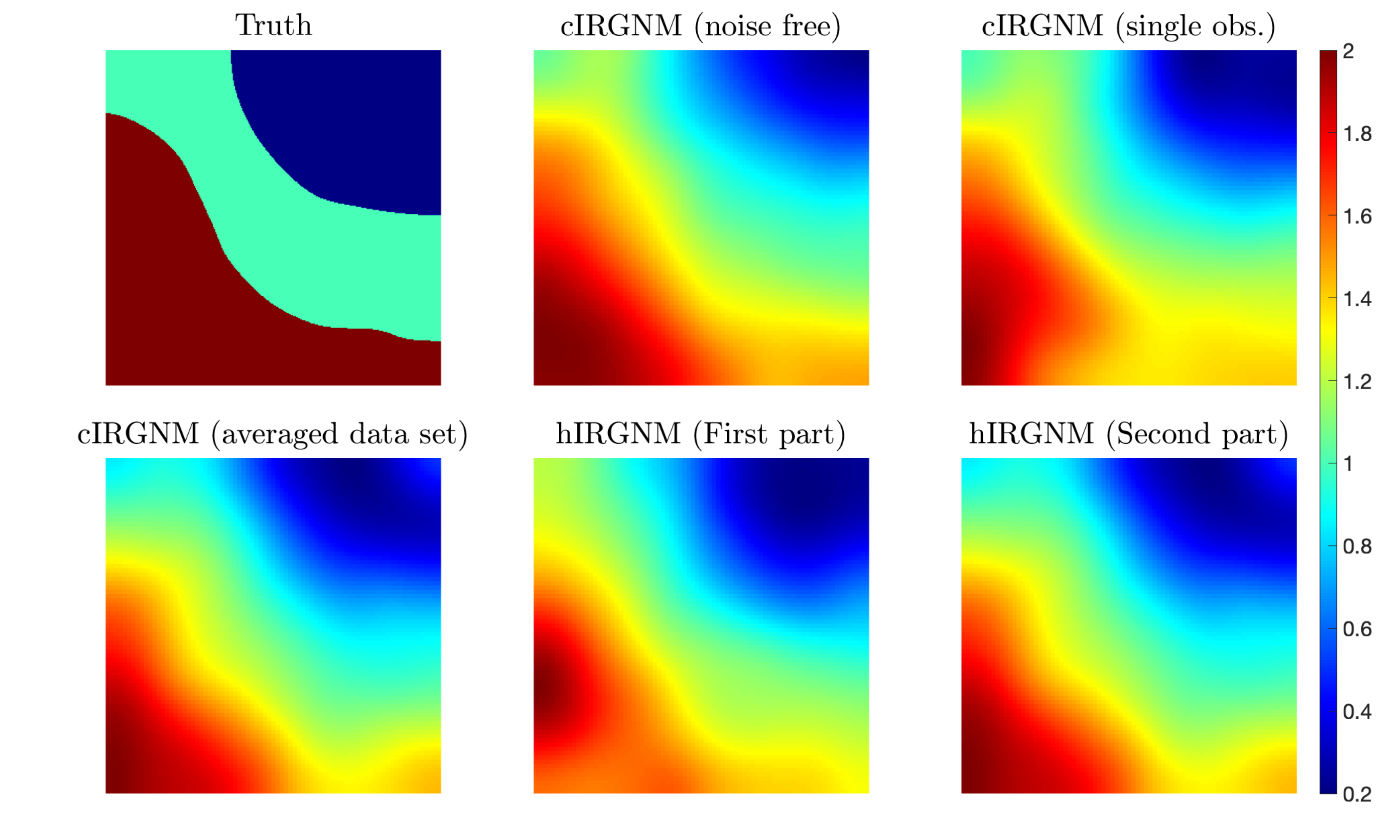}
 \caption{Example 2 (case with the discontinuous truth). Truth (top-left) and estimates of the unknown obtained with the cIGRNM with noise-free observations (top-middle), a single set of observations (top-right) and the average of $N=500$ observations (bottom-left). Bottom-middle and bottom-right panels show the estimates obtained from the first and second part of the hIRGNM using the same $N=500$ observations.}\label{FigExp2_3}
\end{figure}

Similar as in the previous example, we focus on the convergence of dIRGNM and the comparison between hIRGNM and IRGNM with the same finitely many observations. The former is validated in Figure \ref{FigExp2_4} again with different choice of $\beta$ which verifies the theoretical prediction. Namely that an ideal choice of $\beta$ is $\beta \in (1/2,1-\theta/2)$, which provides some stability without overfitting, despite it being slower to converge. Comparison between hIRGNM and cIRGNM with the same finitely many observation is presented in Figure \ref{FigExp2_5}, where one can observe that the hIRGNM ($\beta=1.5$) has already provided an accurate reconstruction in the first dIRGNM iteration and achieve the minimal relative error by just a few iteration in the second  cIRGNM iteration. Also for the smooth truth, the second part of the hIRGNM takes only two iterations to reach the minimum value, similar to the previous example. In particular the minimal relative error by hIRGNM is quite similar to those by cIRGNM. 
To visualize the reconstructed solutions, we provide them for the different algorithms in Figures \ref{FigExp2_2}-\ref{FigExp2_3}.
What we observe from these plots is that the best reconstruction is through the noise-free cIRGNM, however as we also see the worst reconstruction is related to the single observation case of the cIRGNM. Also we finally observe that the averaged observation case of the cIRGNM, matches that of the second part of the hIRGNM. Therefore
we can conclude this matches the phenomenon observed in the previous example.

%Finally we see that for the Figure \ref{FigExp2_4} the effect of the choice of $\beta>0$ related to the step size for the dIRGNM which verifies the theory suggested earlier. Namely that an ideal choice of $\beta$ is $\beta<1$, which provides some stability without overfitting, despite it being slower. This can be seen compared to other choices.

\subsection{Example 3}

The context of our final numerical examples is electrical impedance tomography (EIT) \cite{LB02}. For the forward problem we employ the complete electrode model (CEM) introduced in \cite{SCI92}. We consider a medium with physical domain denoted by $\Omega$ an electric conductivity $\kappa$. A set of electrodes $\{e_l\}_{1=1}^{m_e}$ are attached on the boundary $\partial \Omega$ with contact impedance $\{z_l\}_{l=1}^{m_e}$. The aim of the CEM is to compute the electric potential $\nu$ inside $\Omega$ as well as the voltages $\{V_l\}_{1=1}^{m_e}$ on the electrodes. The governing equations are
\begin{subequations}
\label{eq:CEM}
\begin{alignat}{4}
\nabla \cdot (\exp(u) \nabla \nu) &= 0 , \quad \ \in \ \Omega, \\
\nu + z_l\exp(u)\nabla\nu \cdot \textbf{n} &= V_l , \quad   \in \ e_l,  \ \ l=1,\ldots,m_{e},  \\
\nabla\nu \cdot \textbf{n} &=0 , \quad \   \in  \partial \Omega \textrm{\textbackslash} \cup^{m_{e}}_{l=1}e_{l},  \\
\int_{e_l}\exp(u) \nabla \nu \cdot \textbf{n} \ ds &= I_l , \quad   \in \ e_l,  \ \ l=1,\ldots,m_{e},
\end{alignat}
\end{subequations}
where $u=\log(\kappa)$, $\textbf{n}$ denotes the outward normal vector on the boundary and $I_l$ ($l=1\dots,m_e$) is the current injected through the electrode $e_l$. We require that
$$I=(I_{1},\dots,I_{m_e})\in \mathbb{R}_{0}^{m_{e}} \equiv \Big\{V\in  \mathbb{R}^{m_{e}}\Big\vert \sum_{m=1}^{m_{e}}V_{l}=0\Big\},$$
which implies conservation of charge. For $\kappa=\exp(u)\in C(\overline{\Omega})$, the weak form (\ref{eq:CEM}) has a unique solution $(\nu,V)\in H^{1}(\Omega)\times \mathbb{R}_{0}^{m_{e}}$ \cite{SCI92}.

For the inverse problem we employ $n_p$ injection patterns $\mathbf{I}_{j}=\{I_{j,k}\}_{k=1}^{m_{e}}$ ($j=1,\dots,n_{p}$), and pose the EIT problem of estimating the unknown (log) conductivity $u$ from measurements of $\mathbf{V}_{j}=\{V_{j,k}\}_{k=1}^{m_{e}}$ ($j=1,\dots,n_{p}$). The forward map $F:\mathcal{H}\to \mathbb{R}^{n_{p}m_{e}}$ is defined by $F(u)=\mathbf{V}\equiv (\mathbf{V}_{1},\dots,\mathbf{V}_{n_p})$ where, as in the second example, $\mathcal{H}$ is defined via (\ref{eq:darcy00}). 

The question whether the tangential cone condition is satisfied in this example has received considerable attention during the recent decade, see e.g. \cite{} and the references therein. However, it remains unclear whether Assumption \ref{ass:tcc} holds true. Concerning Assumptions \ref{ass:vsc} and \ref{ass:smoothing}, the same comments as in the previous example apply.

\subsubsection{Numerical results}
In this example, we implement Algorithms \ref{alg:cIRGNM}-\ref{alg:hIRGNM} in \texttt{MATLAB} using the toolbox EIDORS \cite{Al06} to solve (\ref{eq:CEM}) with the Finite Element method. Contact impedances $\{z_l\}_{l=1}^{m_e}$ are chosen with value $0.01$. We employ $m_{e}=16$ electrodes and $n_{p}=16$ injection patterns in which current of $0.1$ Amps is injected through each pair of adjacent electrodes.

Similar to the implementation of the second example, we use a discretized version of the update formula from (\ref{eq:darcy02}). The midpoint rule is applied for the discretization $\mathcal{C}$ in (\ref{eq:eit100}). The parameters for the Mat\'{e}rn correlation function (\ref{matern}) are $c_{0}=25$, $\nu=1.0$ and $\ell=0.1$. For the discretised Fr\'echet derivative $\op'(u_{n})$ we use the built-in command in EIDORS \texttt{calc\_jacobian} which yields the matrix $D_{\kappa}\mathbf{V}$. Then, via the chain rule we compute $F'[u]=D_{\kappa}\mathbf{V}\exp(u)$.

\begin{figure}[h!]
\centering
\includegraphics[scale=0.37,trim=0 0 0 0]{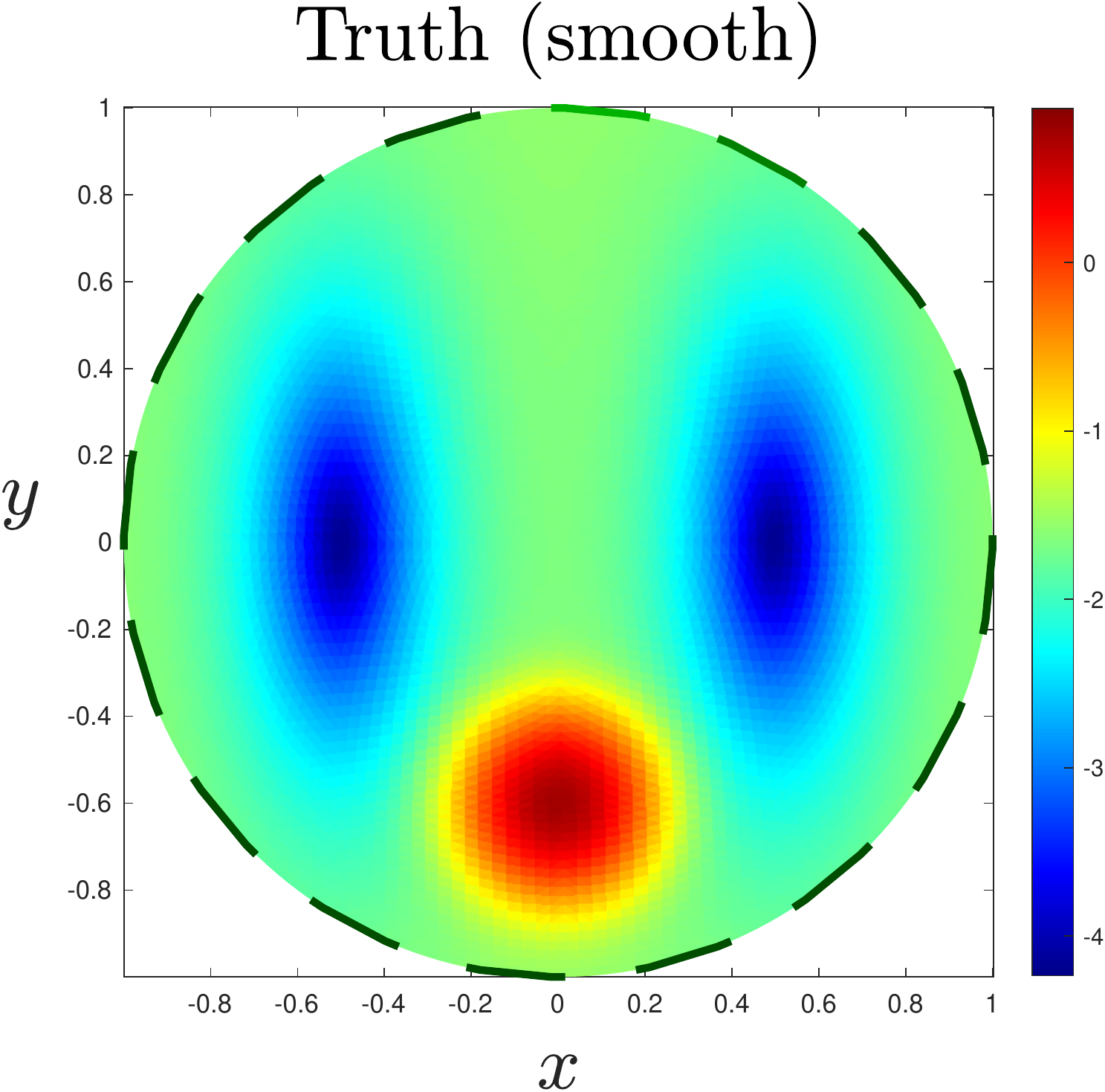}
\includegraphics[scale=0.37,trim=0 0 0 0]{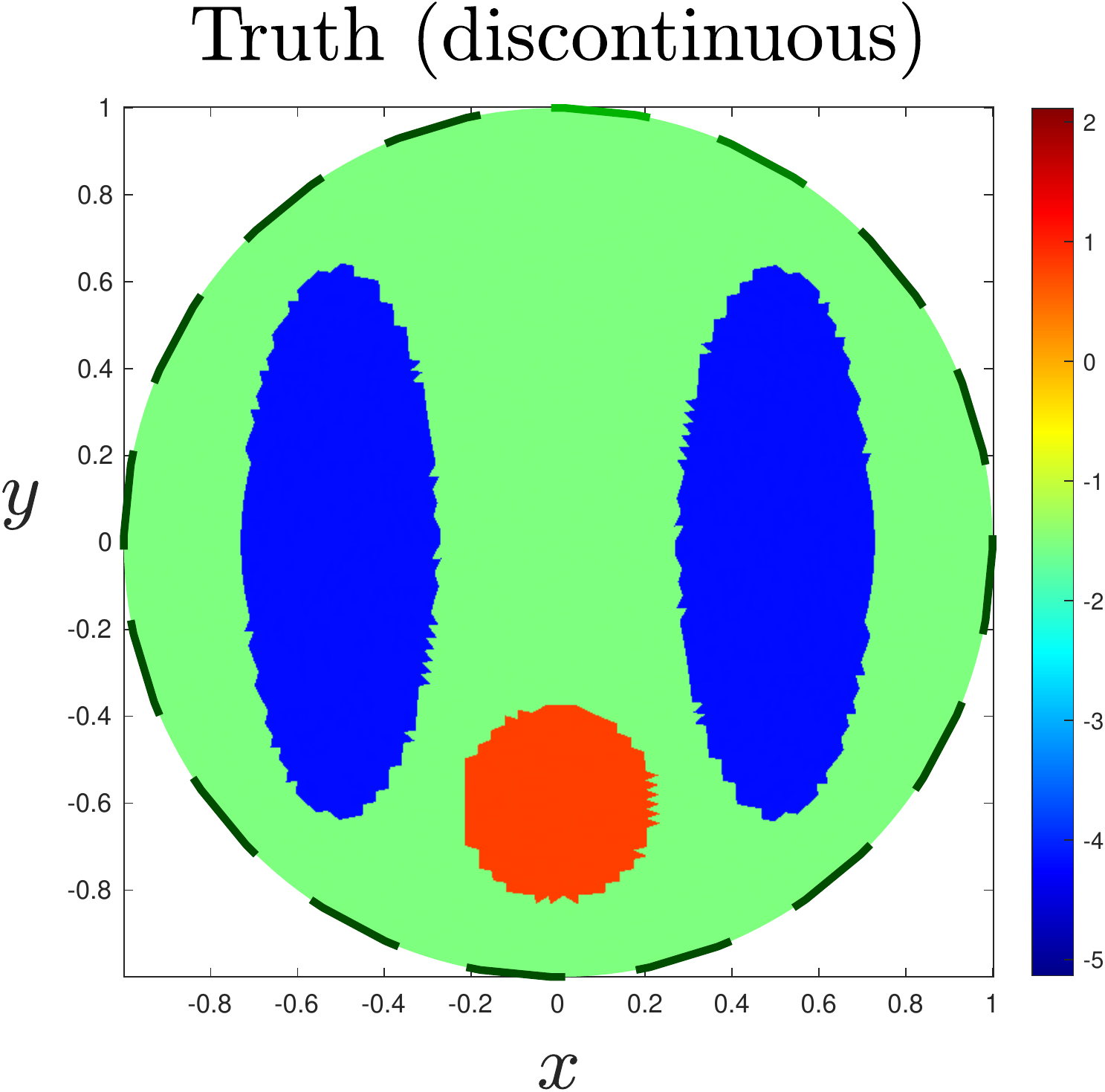}\\
 \caption{Example 3. True coefficient $u^{\dagger}(x,y)$ for the smooth (left) and discontinuous (right) case.}

    \label{FigExp3_1}
\end{figure}

We conduct two set of examples in which we use a smooth and a discontinuous truth shown in Figure \ref{FigExp3_1}. Noise free voltages are computed solving (\ref{eq:CEM}) using the truth and the collection of injection patterns. To avoid inverse crimes we employ a mesh with 9216 elements while a coarser mesh (with 7744 elements) is used for the computations in Algorithms \ref{alg:cIRGNM}-\ref{alg:hIRGNM}.
As before, noisy observations $Y_{n}$ are obtained by adding Gaussian noise to the noise-free measurements as in eq. \eqref{eq:model_dyn}. We use standard deviation of $\sigma=2.5\times 10^{-3}$. Furthermore, we use $\hat{u}_{0}(x,y)=-1$ (for all $(x,y)\in \Omega$),  $\alpha_{0}=10^{-3}$, and $C_{dec}=1.5$.

\begin{figure}[h!]
\centering
\includegraphics[scale=0.33,trim=0 0 0 0]{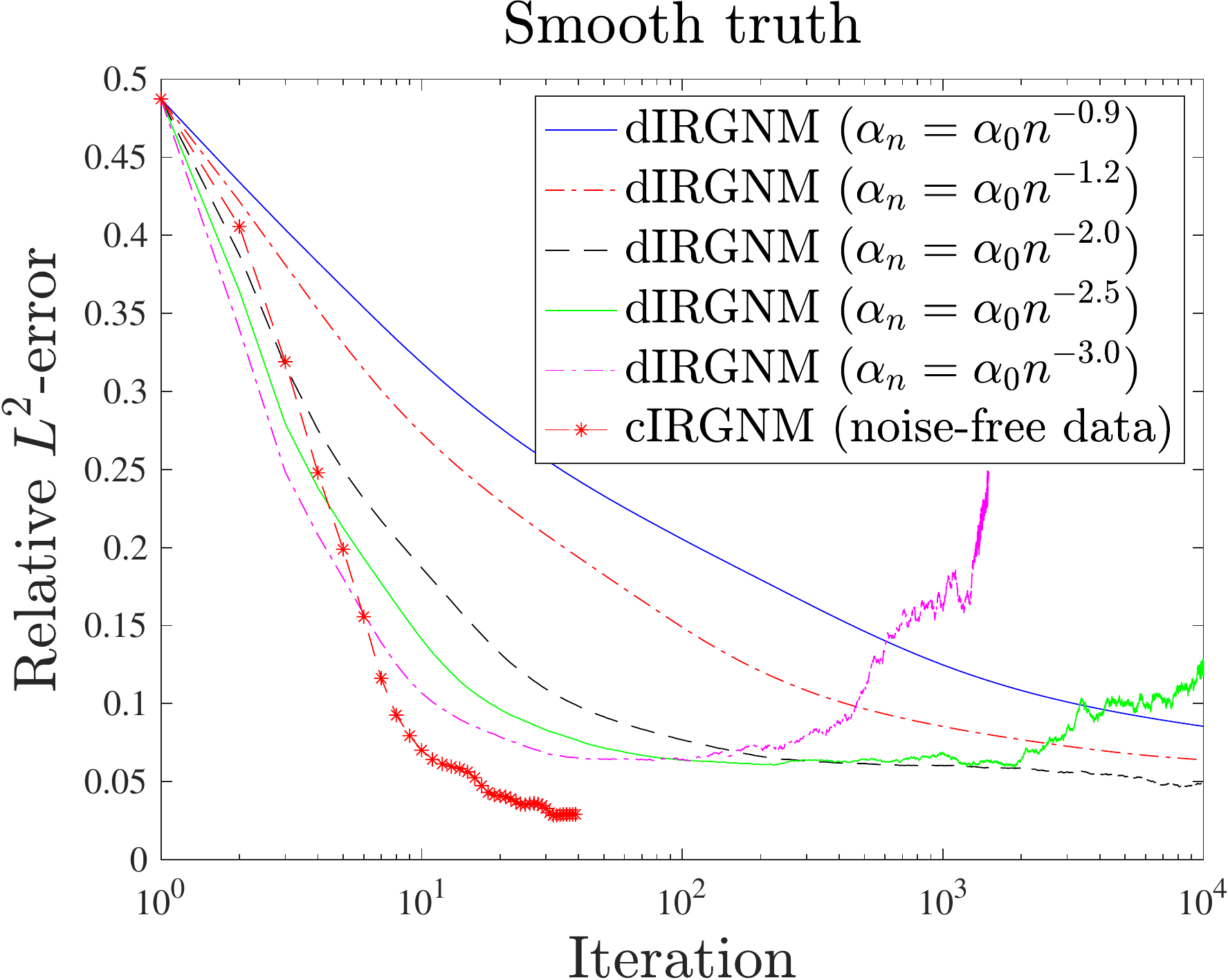}
\vspace{3mm}
\includegraphics[scale=0.33,trim=0 0 0 0]{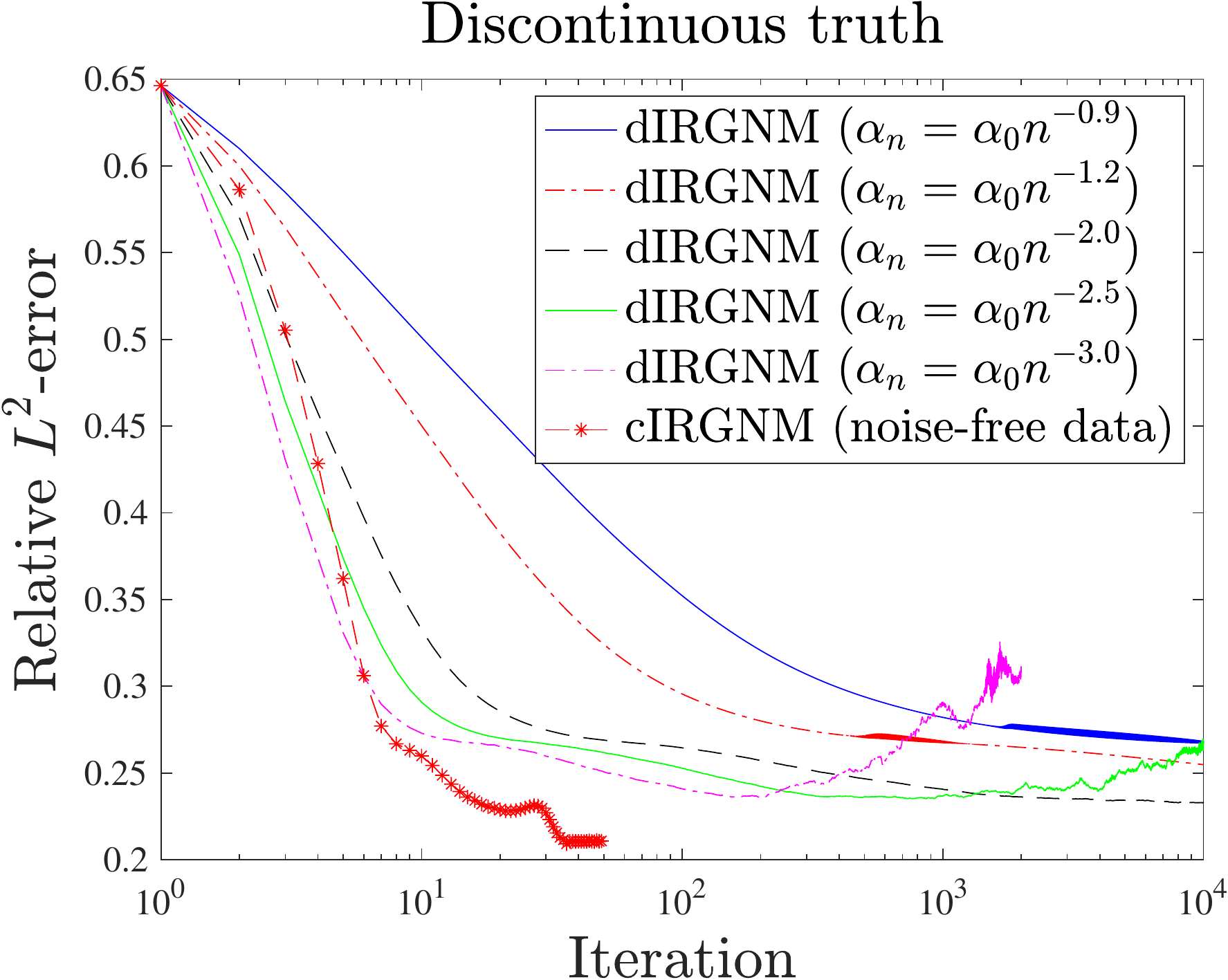}
 %\caption{Initial experiment for the iRGNM with the true underlying unknown, and the averaged noise after $n=25$ iterations.}
\caption{Example 3. Relative $L^{2}$ errors obtained using the dIRGNM with various choices of $\beta$ for the continuous (left) and discontinuous (right) truth.}
    \label{FigExp3_4}
\end{figure}
Again we focus on the convergence of dIRGNM and the comparison between hIRGNM and IRGNM with the same finitely many observation. The relative error w.r.t. the truth obtained using the dIRGNM for different choices of $\beta$ is shown in Figure \ref{FigExp3_4}. Compared with previous two examples, dIRGNM seems to be more robust with respect to the choice of $\beta$ where the amplified relative error appears more obvious when $\beta>2$. We also include the relative error obtained using the IRGNM with noise-free observations as reference. 

\begin{figure}[h!]
\centering
\includegraphics[scale=0.34,trim=0 0 0 0]{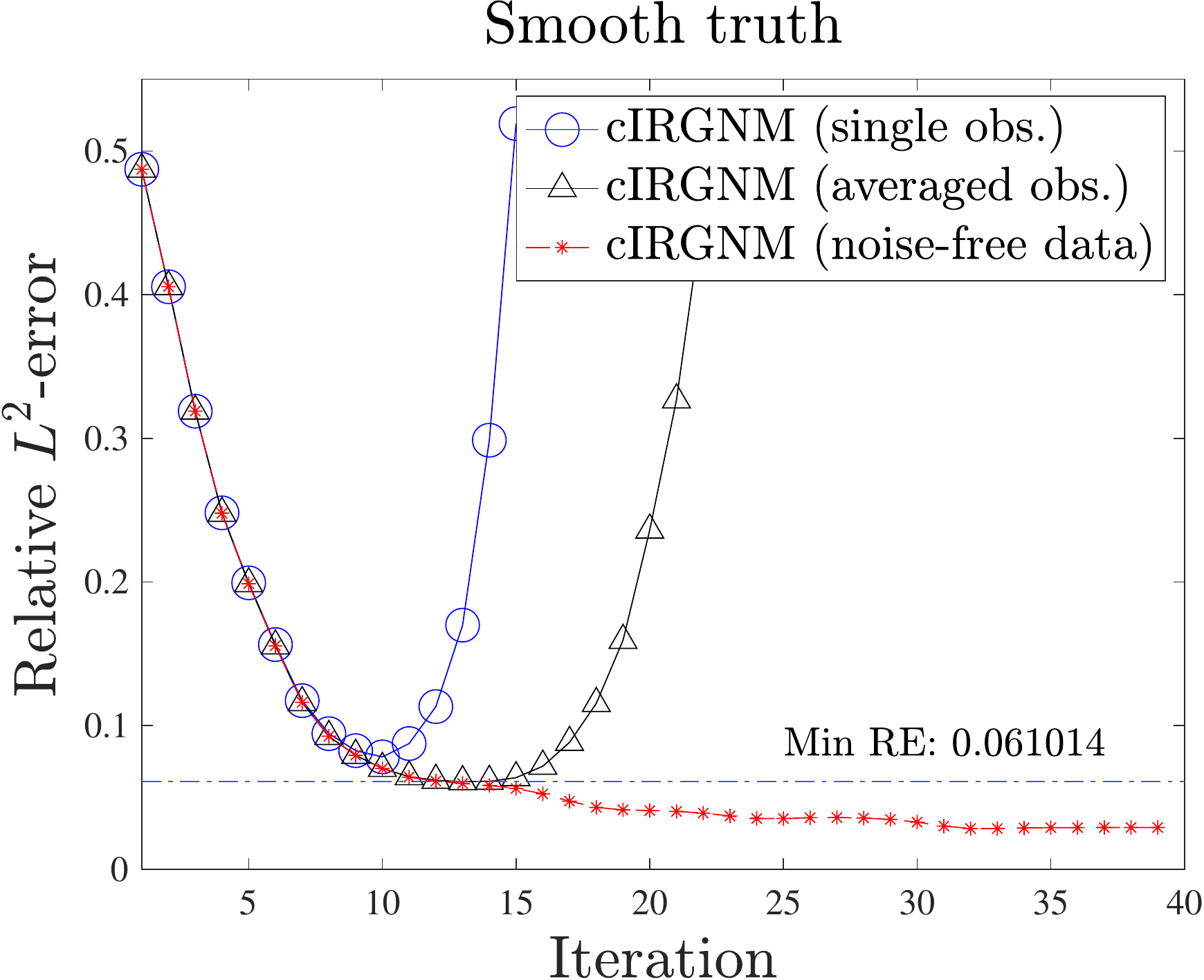}
\includegraphics[scale=0.34,trim=0 0 0 0]{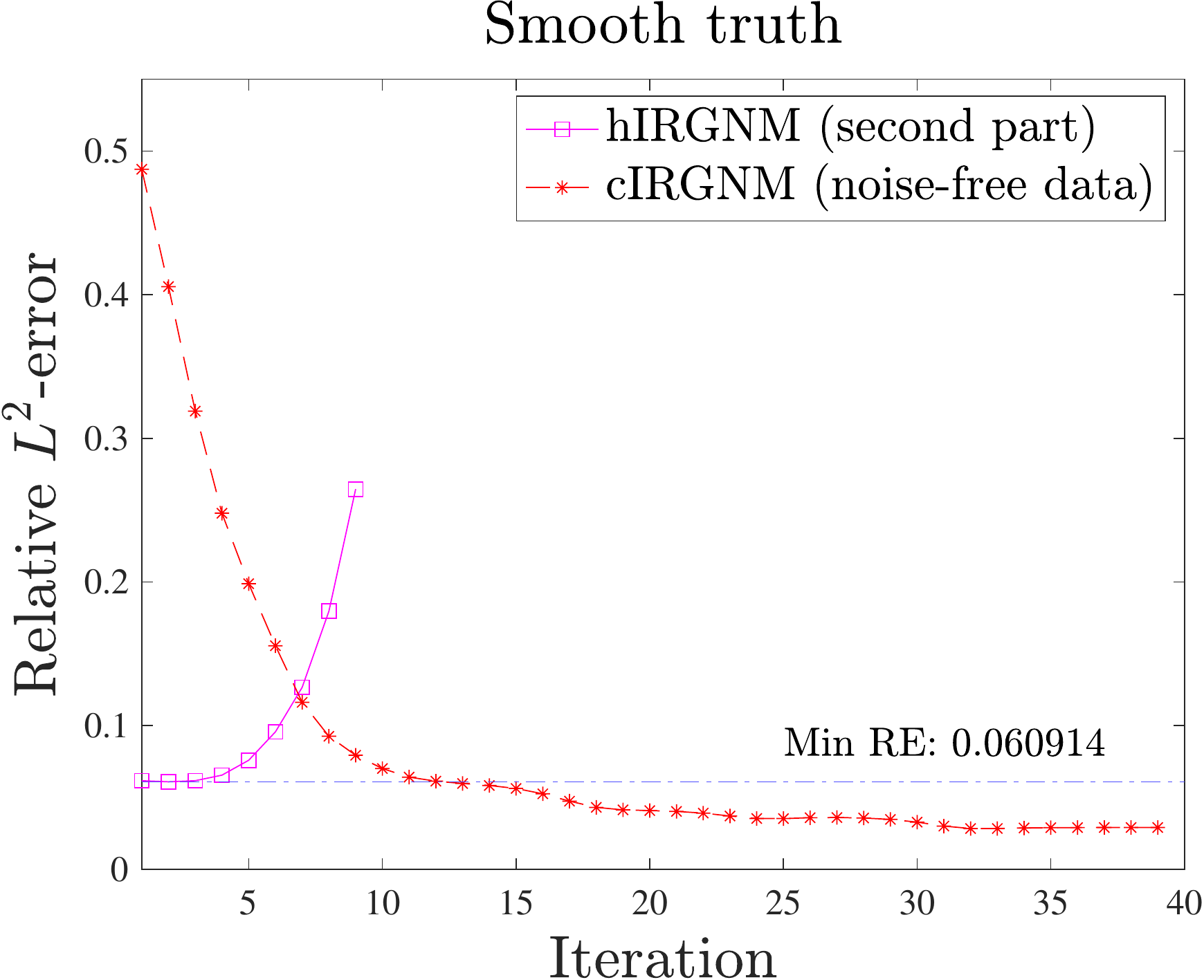}\\
\vspace{3mm}
\includegraphics[scale=0.34,trim=0 0 0 0]{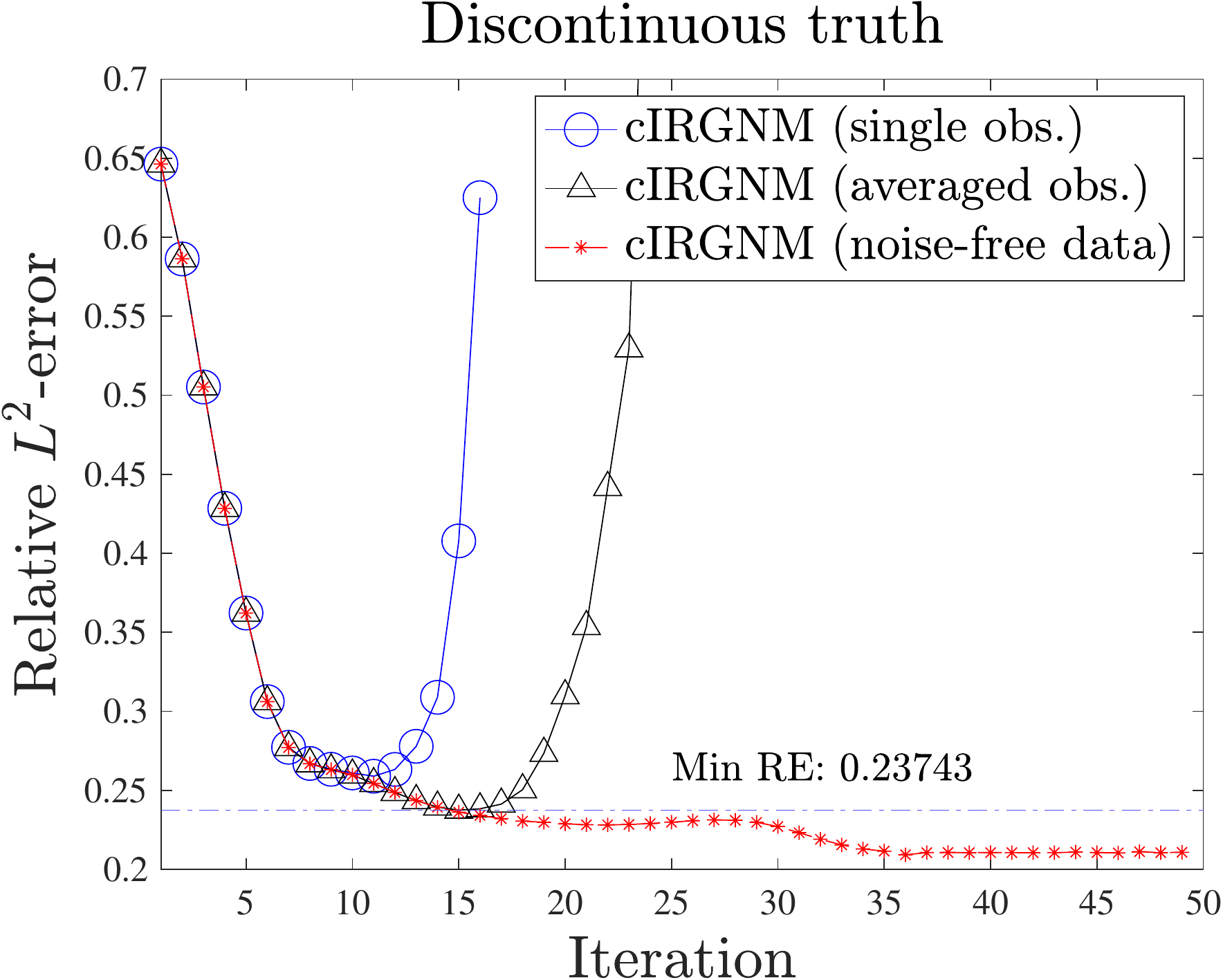}
\includegraphics[scale=0.34,trim=0 0 0 0]{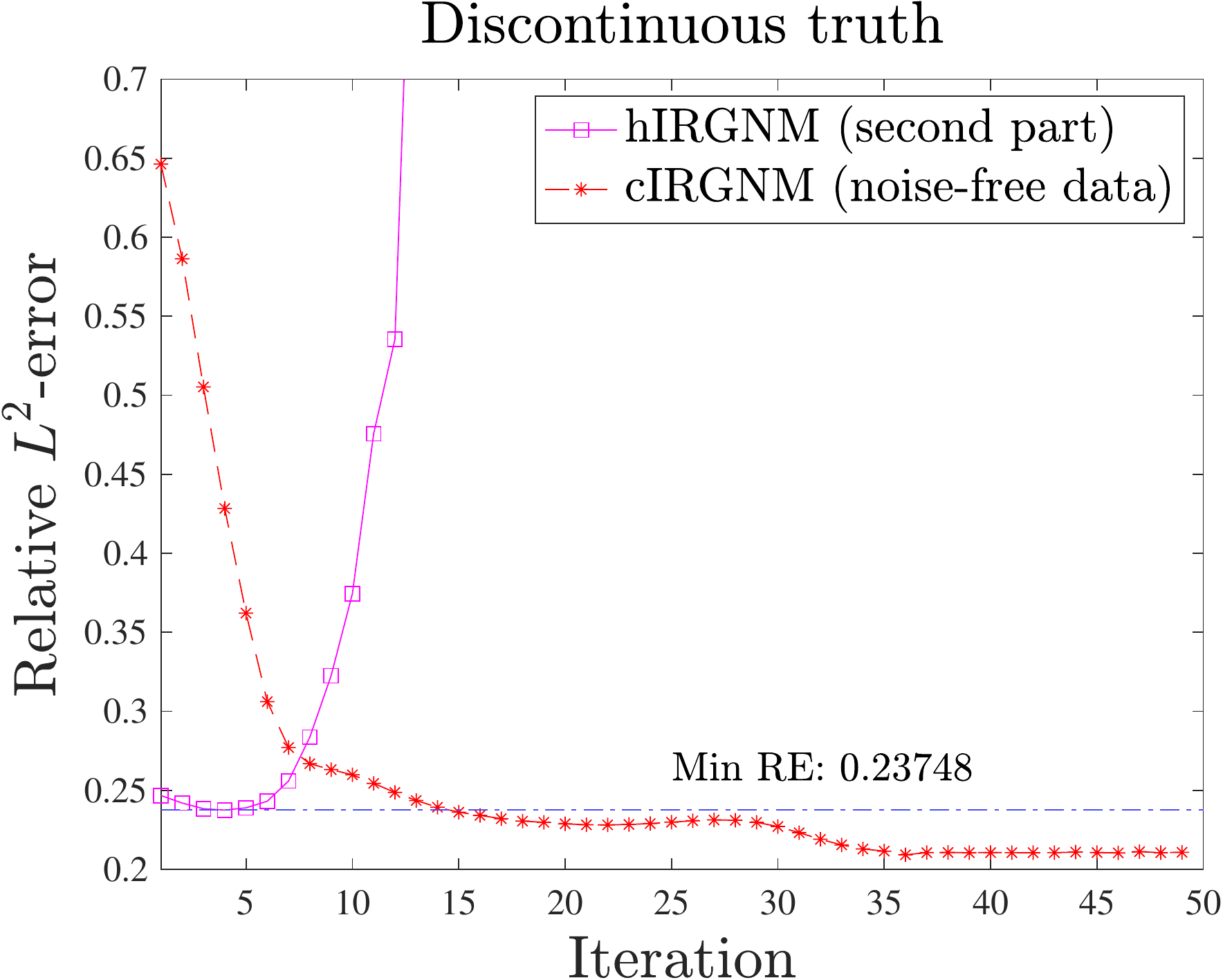}

 %\caption{Initial experiment for the iRGNM with the true underlying unknown, and the averaged noise after $n=25$ iterations.}
\caption{Example 3. Left: Relative $L^{2}$ errors for the case with the smooth (top) and discontinuous (bottom) truth obtained using the cIRGNM with three different observations: (i) noise-free, (ii) a single one and (iii) the averaged of $N=500$. Right: Relative $L^{2}$ errors obtained during the second part of the hIRGNM with the same $N=500$ observations. For comparisons the right panels also display the relative error obtained with the noise-free cIRGNM. The numerical values displayed on the left (resp. right) plots corresponds to the minimum relative error achieved via the cIRGNM with averaged measurements (resp. the second part of the hIRGNM).}
    \label{FigExp3_5}
\end{figure}

\begin{figure}[h!]
\centering
\includegraphics[scale=0.37,trim=0 0 0 0]{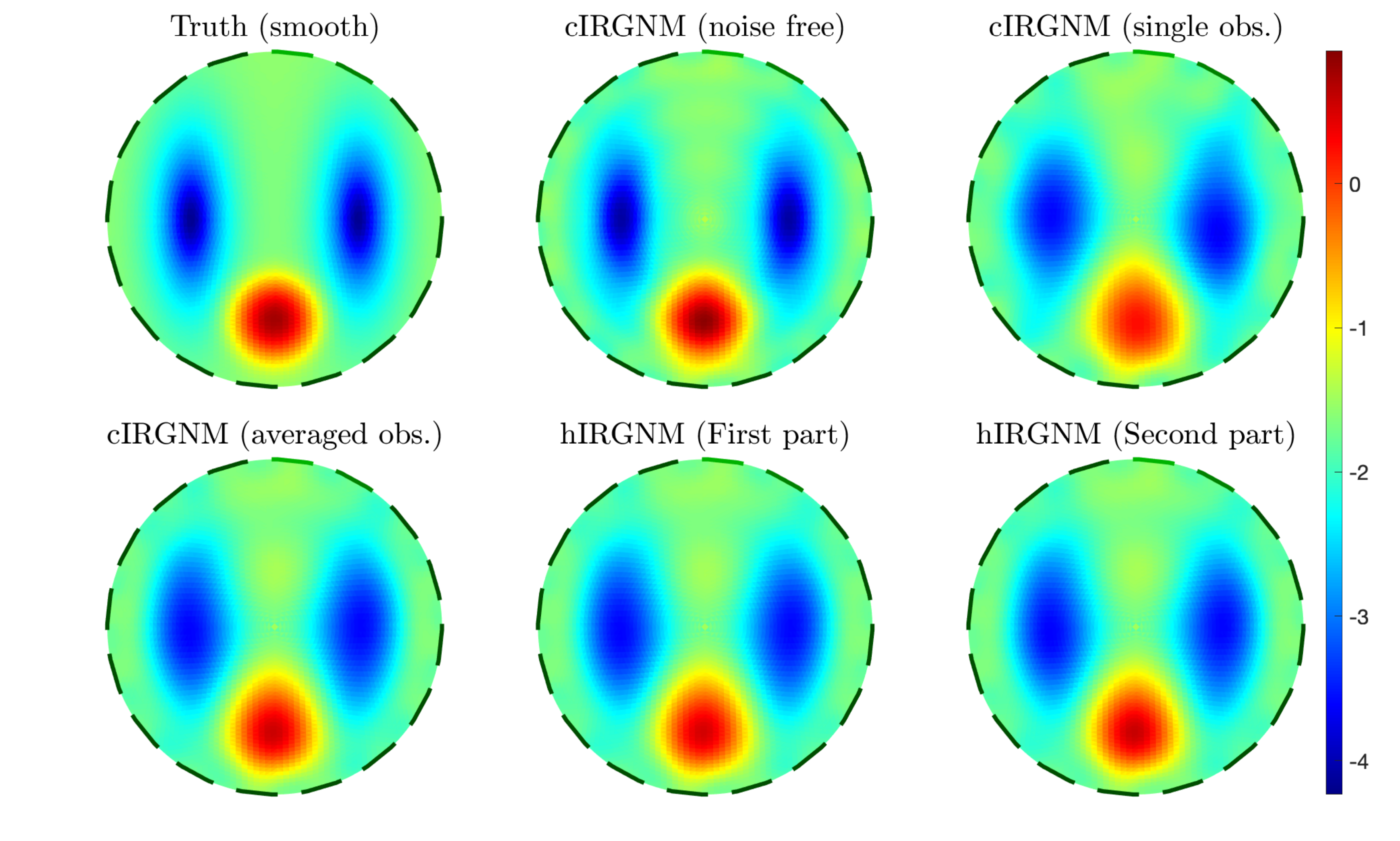}
 \caption{Example 3 (case with the smooth truth). Truth (top-left) and estimates of the unknown obtained with the cIGRNM with noise-free observations (top-middle), a single set of observations (top-right) and the average of $N=500$ observations (bottom-left). Bottom-middle and bottom-right panels show the estimates obtained from the first and second part of the hIRGNM using the same $N=500$ observations.}

    \label{FigExp3_2}
\end{figure}

\begin{figure}[h!]
\centering
\includegraphics[scale=0.37,trim=0 0 0 0]{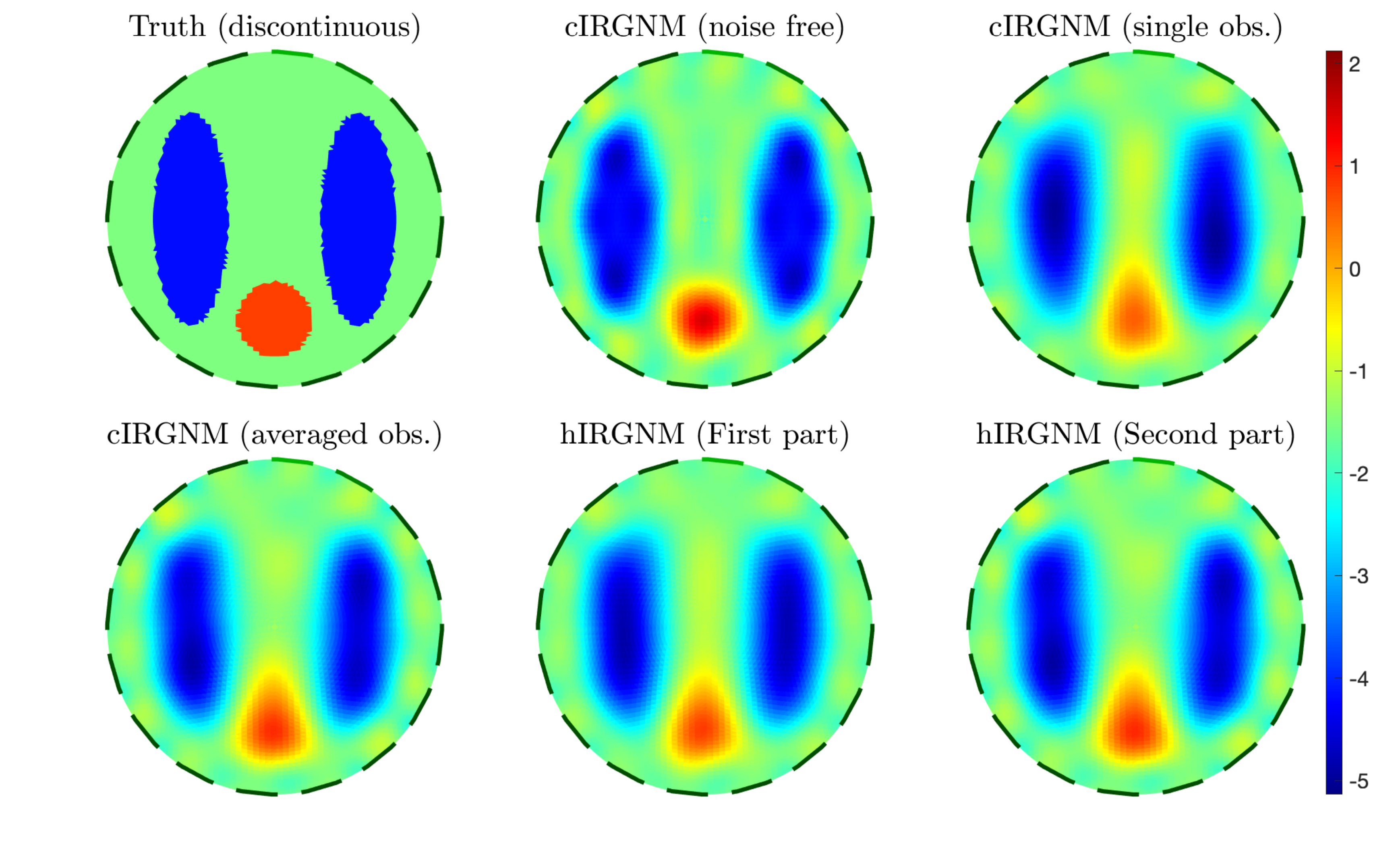}
 \caption{Example 3 (case with the discontinuous truth). Truth (top-left) and estimates of the unknown obtained with the cIGRNM with noise-free observations (top-middle), a single set of observations (top-right) and the average of $N=500$ observations (bottom-left). Bottom-middle and bottom-right panels show the estimates obtained from the first and second part of the hIRGNM using the same $N=500$ observations.}
    \label{FigExp3_3}
\end{figure}

To compare the performance of hIRGNM and cIRGNM, we show relative errors obtained by both algorithms in Figure \ref{FigExp3_5}. Algorithm \ref{alg:cIRGNM} is realized with the different type of observations (i.e. noise-free, single set and the average observation of $N=500$). As comparison, we apply the hIRGNM using the same $N=500$ observations and a value $\beta=2.0$ which, as mentioned earlier, produced stable results when using the dIRGNM with large $N$. The iterations achieved during the second part of the hIRGNM are shown in the right panels of Figure \ref{FigExp3_5}. The value displayed on these plots corresponds to the minimum relative error attained during the second part of hIRGNM. Similar to our previous experiment, we notice that while this value is approximately equal to the value obtained via the cIRGNM with averaged observations, the second part of the hIRGNM reached this minimum value in less number of iterations. In the case with the smooth truth only two iterations sufficed to attain such a minimum value while 14 iterations were required by the cIRGNM. 
Finally, the estimates obtained during the first and second part (when the minimum is attained) are shown in the bottom-middle and bottom-left panels of Figures \ref{FigExp3_2}-\ref{FigExp3_3}. We can see from these plots that the first part of the hGIRNM yields an estimate that is already very close to the truth.
As a comparison, the top-right and bottom-left panels of Figures \ref{FigExp3_2}-\ref{FigExp3_3} show the estimates of the unknown computed when the relative errors attain the minimum value. Again, a lower minimum (display on the plots) is achieved using the cIGRNM with average of all observations compared to the value when using a single set.
Top-middle panels of Figures \ref{FigExp3_2}-\ref{FigExp3_3} show the estimates obtained with the noise-free case.
%\todo[inline]{Discussion re. theory and explain here $beta$ can be even larger}

\section{Conclusion}
\label{sec:conc}
The purpose of this work was to investigate a dIRGNM (\ref{eq:IRGNM}) solving nonlinear inverse problems with
sequential observations. The idea behind our work is highly inspired by the artificial dynamic proposed in \cite{ILS13}
where we need to consider an augmented form, i.e. below
\begin{align}\label{eq:model_dyn2_modified}
\left\{
\begin{array}{l}
u_{n+1}  =u_{n}, \\
Y_{n+1} = \op\left(u_{n+1}\right) + \sigma \xi_{n+1}, \\
Z_{n+1} : =  \frac{1}{n+1}(nZ_{n}+Y_{n+1}) = \frac{1}{n+1}\sum_{i=1}^{n+1} Y_i,
\end{array}
\right.
\end{align}
with $u_0=\udag$.
Such an artificial dynamic contains a steady state equation associated with the unknown variable $\udag$ and two other observation equations with sequential observation $\{Y_n\}_{n=1,\ldots}$ and its average $\{Z_n\}_{n=1,\ldots}$. The proposed dIRGNM (\ref{eq:IRGNM}) is exactly an online filter algorithm towards the artificial dynamic \eqref{eq:model_dyn2_modified}. Systematic convergence analysis of this reconstruction algorithm has been provided in Sections \ref{sec:analysis}-\ref{sec:error} where the averaged observation $Z_n$ yields a vanishing asymptotical behavior if the regularization parameter is appropriately chosen. Such an observation verifies that the uncertainty of the nonlinear inverse problems has been dramatically weaken if the averaged observation is taken, i.e. $Z_n$ in (\ref{eq_DIRGNM}) or (\ref{eq:model_dyn2_modified}).
%The purpose of this work as the motivate and present nonlinear inverse problems as an artificial dynamical
%process. In particular we considered one well-known methodology, which is the iterative regularized Gauss-
%Newton method (IRGNM). As a result we incorporated new observations at each iteration. From this we provided a
%convergence analysis related to the mean square error, which was specific to case of the averaging of noise. Specifically,
%we were able to demonstrate a convergence rate of $\mathcal{O}(n^{-\frac{1}{2}})$, as $n \rightarrow \infty$.
%We also provided an extension of this to a particular choice o an a-posteriori stopping rule, which is the Lepskii
%principle.
Numerical evidence of our findings were presented through three inverse problems associating with elliptic partial differential equations. This was in terms of the rates attained, but also the numerical
performance of the dIRGNM compared to the cIRGNM.

For future work, there are various different avenues one can consider. Firstly as we considered the cIRGNM,
a natural direction would be other nonlinear methodologies such as the Levenberg--Marquardt method (LMM), which
is well known and has applications to geophysical sciences \cite{MH97,QJ10}.  We have not considered such
an analysis here, as the LMM commonly relies more on spectral methods, rather than a variational methods.
Another direction would be to consider other a-posteriori parameter choice rules for $\alpha_n$.
 Other common examples aside from Lepskii principle \cite{BH05,OVL90}, would include the empirical risk minimization.
 Finally given the results we have obtained, one could aim to characterize the ensemble Kalman filter \cite{GE09,GE94}, related
to inverse problems \cite{CIRS18,CCS21,CT21,ILS13}, in terms of convergence through asymptotic regularization \cite{LNW21}.
As of yet, this has only been achieved for linear filters.

\section*{Acknowledgments}
NKC was supported by KAUST baseline funding.
SL was NSFC (No.11925104),
Science and Technology Commission of Shanghai Municipality (19XD1420500, 21JC1400500). MI was supported by the Engineering and Physical Sciences Research Council, UK
[grant number EP/P006701/1]; through the EPSRC Future Composites Manufacturing Research Hub.

\end{document}